\documentclass{amsart}
\usepackage{amsfonts}
\usepackage{amssymb, amsmath, amsthm}
\usepackage{amscd}
\usepackage{verbatim}
\usepackage{enumerate}
\usepackage{setspace}

\allowdisplaybreaks
\newtheorem{theorem}{Theorem}[section]

\newtheorem{corollary}[theorem]{Corollary}

\newtheorem{definition}[theorem]{Definition}

\newtheorem{lemma}[theorem]{Lemma}

\newtheorem{proposition}[theorem]{Proposition}

\theoremstyle{remark}
\newtheorem{remark}{Remark}[section]
\newtheorem{example}[theorem]{Example}
\numberwithin{equation}{section}

\begin{document}
\title[Regularized family of Ericksen-Leslie system]{On a regularized family
of models for the full Ericksen-Leslie system}
\author{Ciprian G. Gal}
\address{Department of Mathematics $\&$ Statistics, Florida International
University, Miami, FL 33199, USA}
\email{cgal@fiu.edu}
\author{Louis Tebou}
\address{Department of Mathematics $\&$ Statistics, Florida International
University, Miami, FL 33199, USA}
\email{teboul@fiu.edu}

\begin{abstract}
We consider a general family of regularized systems for the full
Ericksen-Leslie model for the hydrodynamics of liquid crystals in $n$%
-dimensional compact Riemannian manifolds. The system we consider consists
of a regularized family of Navier-Stokes equations (including the
Navier-Stokes-$\alpha $-like equation, the Leray-$\alpha $ equation, the
Modified Leray-$\alpha $ equation, the Simplified Bardina model, the
Navier-Stokes-Voigt model and the Navier-Stokes equation) for the fluid
velocity $u$ suitably coupled with a parabolic equation for the director
field $d$. We establish existence, stability and regularity results for this
family. We also show the existence of a finite dimensional global attractor
for our general model, and then establish sufficiently general conditions
under which each trajectory converges to a single equilibrium by means of a
Lojasiewicz-Simon inequality.
\end{abstract}

\maketitle


\section{Introduction}

\label{s:intro}

A nematic liquid crystal is a phase of a material between the solid and
liquid phases, with the liquid phase having a certain degree of
orientational order. The flow in the liquid phase is described by a velocity
$u=(u_{1},...,u_{n})$ and by a director field $d=(d_{1},..,d_{n})$, which
stands for the averaged macroscopic/continuum orientation in $\mathbb{R}^{n}$
of the constituent molecules. One model that governs the flow of the nematic
liquid crystals is the general Ericksen-Leslie system (NS-EL)\ with
Ginzburg-Landau type approximation proposed in \cite{WXL12}. This system
consists of the Navier-Stokes equation for the fluid velocity coupled with
two additional anisotropic stress tensors, which are the elastic (Ericksen)
and the viscous (Leslie) stress tensors, respectively, and a parabolic
equation for the director field. Among the mathematical rigorous results for
the full (NS-EL) system one can barely find the references \cite{CRW, WXL12}
for incompressible fluid flows. These contributions are mainly concerned
with well-posedness and long-time behavior of solutions to the system under
suitable assumptions on the Leslie coefficients, ensuring that a certain
natural energy associated with the (NS-EL) system is dissipated. Especially
in \cite{CRW}, existence of a global-in-time weak solutions with \emph{%
finite energy} is proved as well as blow-up criterion is developed for the
existence of a globally-well defined classical solution of the 3D (NS-EL)
system with periodic boundary conditions. On the other hand, in \cite{WXL12}
global well-posedness of smooth solutions is established in certain special
cases and Lyapunov stability for this system near local energy minimizers is
shown. Due to the highly nonlinear and strong coupling of the (NS-EL) system
most of previous analytical studies were always restricted to some
simplified versions of the (NS-EL) system. Rather than giving a full account
of the literature, we refer the reader to \cite{CRW} where a complete
description of the most up-to-date analytical studies has been undertaken in
detail for these simplified models.

Regularized flow equations in hydrodynamics play a key role in understanding
turbulent phenomena in science. Given the nonlinear nature of turbulent
nematic liquid crystal flows and the ensuing multiscale interactions, direct
numerical simulations of the turbulent nematic liquid crystal flows is still
presently lacking apart from some investigations performed on simplified
systems which still retain the basic nonlinear structure and the essential
features of the full hydrodynamic (NS-EL) equations \cite{BFP, BGG, LLZ07}.
This is due mainly to two factors: (a) the numerical computation of the 3D
Navier-Stokes (NSE) equation with a high Reynolds number in regimes in which
the nonlinearities prevail is not possible at present \cite{DJT}\ and (b)
the strong coupling in the Ericksen-Leslie (NS-EL) equations make the
numerical approximation and computation of the solution quite expensive for
the nowadays computers, even for simplified versions of the original system
\cite{BGG}. Indeed, in turbulent flows most of the computational
difficulties lie in the understanding the dynamic interaction between small
and large scales of the flow \cite{OT}. Moreover, the number of degrees of
freedom needed to simulate the fluid flow increases quite drastically as a
function of the Reynolds number. In order to overcome these issues, in
recent years the approach of regularization modeling has been proposed and
tested successfully against experimental data for the 3D NSE equation. One
novelty of this approach is that the regularization models of the 3D
Navier-Stokes equation only modifies the spectral distribution of energy,
and the well-posedness (i.e., existence, uniqueness and stability with
respect to the initial data) of solutions can be rigorously proven unlike
for the 3D (NSE) equation \cite{HLT}. In order to handle these problems for
a simplified model of the original Ericksen-Leslie system, a general
three-parameter family of regularized equations has been proposed and
investigated in \cite{GM-JNS}\ for the purpose of direct numerical
simulations of turbulent incompressible flows of nematic liquid crystals.
Existence and uniqueness of smooth solutions can be rigorously proven for
the regularized family of \cite[Section 7]{GM-JNS}, as well as the existence
of finite dimensional global attractors and, under proper natural
conditions, the eventual asymptotic stabilization of the corresponding
solutions to single equilibria. The robust analytical properties of these
simplified Ericksen-Leslie models ensure computability of their solutions
and the stability of numerical schemes.

In this paper, our main goal is to investigate a wide range of regularized
models for the general Ericksen-Leslie system (NS-EL)\ with Ginzburg-Landau
type approximation proposed in \cite{WXL12}. As in \cite{GM-JNS}, we will
mainly be concerned with the same fundamental issues in the theory of
infinite-dimensional dynamical systems, that to give a unified analysis of
the entire family of regularized models and to establish existence,
stability and regularity results, and long-time results. As a representative
of a more general model, described in detail in the next section, the family
of regularized models associated with the original (NS-EL) system we wish to
consider formally reads%
\begin{equation}
\left\{
\begin{array}{l}
\partial _{t}u+A_{0}u+(Mu\cdot \nabla )(v)+\chi \nabla (Mu)^{T}\cdot
(v)+\nabla p=-\text{div}(\nabla d\odot \nabla d-\boldsymbol{\sigma }_{Q})+g,
\\
\partial _{t}d+v\cdot \nabla d-\omega _{Q}d+\frac{\lambda _{2}}{\lambda _{1}}%
A_{Q}d=\frac{1}{\lambda _{1}}\left( A_{1}d+\nabla _{d}W(d)\right) , \\
u=Q^{-1}v, \\
\text{div}\left( u\right) =0,\text{ div}\left( v\right) =0, \\
u\left( 0\right) =u_{0}, \\
d\left( 0\right) =d_{0}.%
\end{array}%
\right.  \label{e:rel}
\end{equation}%
Here, $A_{0}$, $A_{1}$, $M$, and $Q$ are linear operators having certain
mapping properties and $\chi $ is either $1$ or $0.$ The function $W\left(
d\right) =(\left\vert d\right\vert ^{2}-1)^{2}$ is used as a typical
approximation to penalize the deviation of the length $|d|$ from the value $%
1 $, under a generally accepted assumption that the liquid crystal molecules
are of similar size \cite{WXL12}. Following \cite{HLT} (cf. also \cite%
{GM-JNS}), there are three parameters which control the degree of smoothing
in the operators $A_{0}$, $M$ and $Q$, namely $\theta $, $\theta _{1}$ and $%
\theta _{2}$, while $A_{1}$ is a differential operator of \emph{second}
order. Some examples of operators $A_{0}$, $A_{1}$, $M$, and $Q$ which
satisfy the mapping assumptions imposed in this paper are
\begin{equation}
A_{0}=\mu _{4}(-\Delta )^{\theta },\quad A_{1}=-\Delta \text{, \ \ }%
M=(I-\alpha ^{2}\Delta )^{-\theta _{1}},\quad Q=(I-\alpha ^{2}\Delta
)^{-\theta _{2}},  \label{example}
\end{equation}%
for fixed positive real numbers $\alpha ,\mu _{4}$ and for specific choices
of the real parameters $\theta $, $\theta _{1}$, and $\theta _{2}$. The
notation $\nabla _{d}$ represents the gradient with respect to the variable $%
d$. Besides, the term $\nabla d\odot \nabla d$ denotes the $n\times n$
matrix whose $(i,j)$-th entry is given by $\nabla _{i}d\cdot \nabla _{j}d$,
for $1\leq i,j\leq n$, while for $v:=Qu,$%
\begin{equation}
A_{Q}=\frac{1}{2}(\nabla v+\nabla ^{T}v),\ \,\;\omega _{Q}=\frac{1}{2}%
(\nabla v-\nabla ^{T}v),  \label{v4}
\end{equation}%
represent the rate of strain tensor and the skew-symmetric part of the
strain rate, respectively. Moreover, as in \cite{CRW, WXL12} we denote by
\begin{equation}
\dot{d}=\partial _{t}d+v\cdot \nabla d,\ \,\;\mathcal{N}_{Q}=\dot{d}-\omega
_{Q}d  \label{def1}
\end{equation}%
the material derivative of $d$ and the rigid rotation part of the changing
rate of the director by fluid vorticity. The kinematic transport $\lambda
_{1}\mathcal{N}_{Q}+\lambda _{2}A_{Q}d$ represents the effect of the
macroscopic flow field on the microscopic structure such that the material
coefficients $\lambda _{1}$ and $\lambda _{2}$ reflect the molecular shape
and how slipper the particles are in the fluid, respectively. The Leslie
stress tensor $\boldsymbol{\sigma }_{Q}$ takes on the following general
form:
\begin{equation}
\boldsymbol{\sigma }_{Q}=\mu _{1}(d^{T}A_{Q}d)d\otimes d+\mu _{2}\mathcal{N}%
_{Q}\otimes d+\mu _{3}d\otimes \mathcal{N}_{Q}+\mu _{5}(A_{Q}d)\otimes d+\mu
_{6}d\otimes (A_{Q}d),  \label{v5}
\end{equation}%
where $\otimes $ stands for the usual Kronecker product, i.e., $(a\otimes
b)_{ij}:=a_{i}b_{j}$, for $1\leq i,j\leq n$. The six independent
coefficients $\mu _{1},...,\mu _{6}$ from (\ref{example}) and (\ref{v5}),
are called Leslie coefficients. Finally, $g=g\left( t\right) $ is an
external body force acting on the fluid.

It is rather clear that one recovers the original (NS-EL) system of \cite%
{CRW, WXL12} by setting $\theta =1,$ $\theta _{1}=\theta _{2}=0$ and $\chi
=0 $ in (\ref{e:rel}). We recall that some theoretical aspects (i.e.,
existence of globally-defined weak solutions and blow-up criteria for smooth
solutions) have been recently developed in \cite{CRW}. Beyond \cite{CRW} and
\cite{WXL12}, not much else seems to be known in terms of analytical and
numerical results for this system to the best of our knowledge. It is worth
noting that among the models considered in (\ref{e:rel})-(\ref{example}),
when restricted only to the equation for the fluid velocity, one can also
find the globally well-posed 3D Leray-$\alpha $ equations, the modified 3D
Leray (ML) equations, the simplified 3D Bardina (SBM) models, the 3D
Navier-Stokes-Voigt (NSV) equations, and their inviscid counterparts. The
corresponding parameter values of $\left( \theta ,\theta _{1},\theta
_{2}\right) $ and operators (\ref{example}) associated with these models are
described in detail in Section \ref{s:prelim}. Inspired by work in \cite%
{GM-JNS} performed on a simplified family of Ericken-Leslie models, we
proceed to develop a complete theory for the whole family of (\ref{e:rel}).
First, we develop well-posedness and long-time dynamics results for the
entire three-parameter family of models (\ref{e:rel}), and then subsequently
recover results of this type for the specific regularization models that
have not been previously studied in the literature, including results for
the original (NS-EL) system.

The main novelties of the present paper with respect to previous results on
the original (NS-EL) model are the following:

(i) The existence result of the globally-defined weak solutions with \emph{%
finite} energy\ are extended to the general family of models (\ref{e:rel})
on an $n$-dimensional compact Riemannian manifold with or without boundary.
We also address both cases of dimension when $n=2,3.$ Furthermore, our
setting allows for the treatment of all kinds of boundary conditions (i.e.,
periodic, no-slip, no-flux, etc)\ for $\left( u,d\right) $; they will be
incorporated in the weak formulation for the problem (\ref{e:rel}) and the
information associated with the dissipation and smoothing operators from (%
\ref{example}).

(ii) We establish general results on regularity, uniqueness and continuous
dependence with respect to initial data for the family (\ref{e:rel}) in the
general case when $\lambda _{2}\neq 0$ and $\mu _{1}\geq 0$.

(iii) We prove results on the existence of finite-dimensional of global
attractors, and existence of exponential attractors (also known as inertial
sets) for the entire three-parameter family (\ref{e:rel}) in the general
case of (ii) and when $\theta >0$. Due to loss of compactness of the
semigroup associated with problem (\ref{e:rel}), the proofs require a
completely different argument than the one given in \cite{GM-JNS}, for a
simplified Ericksen-Leslie family, and is based on a short trajectory type
technique (see \cite{EZ}).

(iv) We discuss the convergence, as time goes to infinity, of solutions of (%
\ref{e:rel}) to single equilibria. More precisely, by the Lojasiewicz--Simon
technique we establish the convergence of any globally-defined weak solution
of (\ref{e:rel}) with finite energy to a single steady state, regardless of
whether uniqueness is known or not for (\ref{e:rel}), provided that the
time-dependent body force $g$ is asymptotically decaying in a precise way,
i.e.,%
\begin{equation*}
\int_{t}^{\infty }\left\Vert g\left( s\right) \right\Vert _{H^{-\theta
-\theta _{2}}}^{2}ds\lesssim \left( 1+t\right) ^{-\left( 1+\delta \right) },%
\text{ for all }t\geq 0,
\end{equation*}%
for some $\delta \in \left( 0,1\right) $. In particular, we show for any
fixed initial datum $\left( u_{0},d_{0}\right) \in H^{-\theta _{2}}\times
H^{1},$ the corresponding trajectory $\left( u\left( t\right) ,d\left(
t\right) \right) $ satisfies%
\begin{equation}
u\left( t\right) \rightarrow 0\text{ weakly in }H^{-\theta _{2}}\text{ and }%
d(t)\rightarrow d_{\ast }\text{ strongly in }H^{0},  \label{weak-conv}
\end{equation}%
as $t$ tends to $\infty $, where $d_{\ast }$ is a steady-state of $%
A_{1}d_{\ast }+f\left( d_{\ast }\right) =0$. We emphasize that (\ref%
{weak-conv}) holds for all weak solutions satisfying a suitable energy
inequality, and so it holds in particular for the limit points of
approximate solutions constructed within a numerical scheme. This result is
also valid for the original (NS-EL) model (\ref{e:rel}) with $\left( \theta
,\theta _{1},\theta _{2}\right) =\left( 1,0,0\right) ,$ $\chi =0$ and
extends a result obtained for a simplified version of the (NS-EL) model
analyzed in \cite{prslong}. Finally, we also give sufficient conditions for
the model (\ref{e:rel}) in order to have a stronger convergence result in (%
\ref{weak-conv}). More precisely, we show that%
\begin{equation*}
u\left( t\right) \rightarrow 0\text{ strongly in }H^{-\theta _{2}}\text{ and
}d(t)\rightarrow d_{\ast }\text{ strongly in }H^{1},
\end{equation*}%
provided that%
\begin{equation*}
\theta +\theta _{2}\geq 1\text{ and }d\text{ belongs to }L^{\infty }\left(
\mathbb{R}_{+};L^{\infty }\left( \Omega \right) \right) .
\end{equation*}

(v) Exploiting the framework of \cite{HLT} which is also extended in \cite%
{GM-JNS}, the abstract mapping assumptions we employ for (\ref{e:rel}) are
more general, and as a result do not require any specific form of the
parametrizations of $A_{0}$, $M$, and $Q,$ as in (\ref{example}). As a
consequence, our framework allows~us to derive new results for a much larger
three-parameter family of models that have not been explicitly studied
elsewhere in detail. Finally, in Section \ref{s:cr} we give some conclusions
about the abstract model and its connection to the standard models as 
introduced in Table \ref{t:spec}. 

The remainder of the paper is structured as follows. In~Section \ref%
{s:prelim}, we establish our notation and give some basic preliminary
results for the operators appearing in the general regularized model.
In~Section \ref{s:well}, we build some well-posedness results for the
general model; in particular, we establish existence results (Section \ref%
{ss:exist}), regularity results (Section \ref{ss:reg}), and uniqueness and
continuous dependence results (Section \ref{ss:stab}). In~Section \ref%
{ss:global}, we show existence of a finite-dimensional global attractor for
the general model by employing the approach from \cite{EZ}. In~Section \ref%
{s:convss}, we establish the eventual asymptotic stabilization as time goes
to infinity of solutions to our regularized models, with the help from a
Lojasiewicz--Simon technique. To make the paper sufficiently self-contained,
our final Section \ref{ss:app} contains supporting material on Sobolev and Gr%
\"{o}nwall-type inequalities, and several other abstract results which are
needed to prove our main results.

\section{Preliminary material}

\label{s:prelim}

\subsection{The functional framework}

\label{ss:ff}

We follow the same framework and notation as in \cite{HLT} (cf. also \cite%
{GM-JNS}). To this end, let $\Omega $ be an $n$-dimensional smooth compact
manifold with or without boundary and equipped with a volume form, and let $%
E\rightarrow \Omega $ be a vector bundle over $\Omega $ endowed with a
Riemannian metric $h=\left( h_{ij}\right) _{n\times n}$. With $C^{\infty
}(E) $ denoting the space of smooth sections of $E$, let $\mathcal{V}%
\subseteq C^{\infty }(E)$ be a linear subspace, let $A_{0}:\mathcal{V}%
\rightarrow \mathcal{V}$ be a linear operator, and let $B_{0}:\mathcal{V}%
\times \mathcal{V}\rightarrow \mathcal{V}$ be a bilinear map. At this point $%
\mathcal{V}$ is conceived to be an arbitrary linear subspace of $C^{\infty
}(E)$; however, later on, we will impose some explicit restrictions on $%
\mathcal{V}$ (see below). Furthermore, we let $\mathcal{W}\subseteq
C^{\infty }(E)$ be a linear subspace and let $A_{1}:\mathcal{W}\rightarrow
\mathcal{W}$ be a linear operator satisfying various assumptions below. In
order to define the variational setting for the phase-field component we
also need to introduce the bilinear operators $R_{0}:\mathcal{W}\times
\mathcal{W\rightarrow V}$, $B_{1}:\mathcal{V}\times \mathcal{W}\rightarrow
\mathcal{W}$, as follows:%
\begin{equation}
B_{1}\left( u\left( x\right) ,d\left( x\right) \right) :=Qu\left( x\right)
\cdot \nabla d\left( x\right) ,\text{ }R_{0}\left( \psi \left( x\right)
,d\left( x\right) \right) :=\psi \left( x\right) \cdot \nabla d\left(
x\right) .  \label{b00}
\end{equation}%
Given the initial data $u_{0}\in \mathcal{V}$, $d_{0}\in \mathcal{W}$ and
forcing term $g\in C^{\infty }(0,T;\mathcal{V})$ with $T>0$, consider the
following system%
\begin{equation}
\left\{
\begin{array}{l}
\partial _{t}u+A_{0}u+B_{0}(u,u)=R_{0}\left( A_{1}d,d\right) +\text{div}%
\left( \sigma _{Q}\right) +g, \\
\partial _{t}d+B_{1}\left( u,d\right) -\omega _{Q}d+\frac{\lambda _{2}}{%
\lambda _{1}}A_{Q}d=\frac{1}{\lambda _{1}}\left( A_{1}d+\nabla
_{d}W(d)\right) , \\
u\left( 0\right) =u_{0},d\left( 0\right) =d_{0},%
\end{array}%
\right.  \label{e:op}
\end{equation}%
on the time interval $[0,T]$. Bearing in mind the model (\ref{e:rel}), we
are mainly interested in bilinear maps of the form
\begin{equation}
B_{0}(v,w)=\bar{B}_{0}(Mv,Qw),  \label{e:b-def}
\end{equation}%
where $M$ and $Q$ are linear operators in $\mathcal{V}$ that are relatively
flexible, and $\bar{B}_{0}$ is a bilinear map fixing the underlying
nonlinear structure of the fluid equation. In the following, denote $%
P:C^{\infty }(E)\rightarrow \mathcal{V}$ as the $L^{2}$-orthogonal projector
onto $\mathcal{V}$. When $\sigma _{Q}\equiv 0,$ $\omega _{Q}\equiv 0$ and $%
\lambda _{2}=0$, the system (\ref{e:op}) corresponds to a simplified
(regularized) Ericksen-Leslie system that was fully investigated in \cite%
{GM-JNS}.

We will study the regularized system (\ref{e:op}) by extending it to
function spaces that have weaker differentiability properties. To this end,
we interpret (\ref{e:op}) in a distributional sense, and need to
continuously extend $A_{0},$ $A_{1}$ and $B_{0},B_{1}$ and $R_{0}$ to
appropriate smoothness spaces. Namely, we employ the spaces $V^{s}=\mathrm{%
clos}_{H^{s}}\mathcal{V}$, $W^{s}=\mathrm{clos}_{H^{s}}\mathcal{W}$ ($H^{s}$
denotes the Sobolev space of order $s$), which will informally be called
Sobolev spaces in the sequel. The pair of spaces $V^{s}$ and $V^{-s}$ are
equipped with the duality pairing $\left\langle \cdot ,\cdot \right\rangle $%
, that is, the continuous extension of the $L^{2}$-inner product to $V^{0}$.
Same applies to the triplet $W^{s}\subset W^{0}=\left( W^{0}\right) ^{\ast
}\subset W^{-s}.$ Moreover, we assume that there are self-adjoint \emph{%
positive} operators $\Lambda $ and $A_{1},$ respectively, such that $\Lambda
^{s}:V^{s}\rightarrow V^{0},$ $A_{1}^{s/2}:W^{s}\rightarrow W^{0}$ are
isometries for any $s\in \mathbb{R}$, and $\Lambda ^{-1}$, $\left(
A_{1}\right) ^{-1}$ are compact operators. For arbitrary real $s$, assume
that $A_{0}$, $A_{1}$, $M$, and $Q$ can be continuously extended so that
\begin{equation}
A_{0}:V^{s}\rightarrow V^{s-2\theta },\quad A_{1}:W^{s}\rightarrow
W^{s-2},\quad M:V^{s}\rightarrow V^{s+2\theta _{1}},\quad \text{and}\quad
Q:V^{s}\rightarrow V^{s+2\theta _{2}},  \label{e:bdd-amn}
\end{equation}%
are bounded operators. Again, we emphasize that the assumptions we will need
for $A_{0}$, $M$, and $Q$ are more general, and do not require this
particular form of the parametrization (see (\ref{e:coercive-a})-(\ref%
{e:coercive-abis}) below). We will assume $\theta ,\theta _{2}\geq 0$ and no
\emph{a priori} sign restriction on $\theta _{1}$. The canonical norm in the
Hilbert spaces $V^{s}$ and $W^{s},$ respectively, will be denoted by the
same quantity $\left\Vert \cdot \right\Vert _{s}$ whenever no further
confusion arises, while we will use the notation $\left\Vert \cdot
\right\Vert _{L^{p}}$ for the $L^{p}$-norm. Furthermore, we assume that $%
A_{0}$ and $Q$ are both self-adjoint, and coercive in the sense that for $%
\beta \in \mathbb{R}$,
\begin{equation}
\left\langle A_{0}w,\Lambda ^{2\beta }w\right\rangle \geq c_{A_{0}}\Vert
w\Vert _{\theta +\beta }^{2}-C_{A_{0}}\Vert w\Vert _{\beta }^{2},\qquad w\in
V^{\theta +\beta },  \label{e:coercive-a}
\end{equation}%
with $c_{A_{0}}=c_{A_{0}}(\beta )>0$, and $C_{A_{0}}=C_{A_{0}}(\beta )\geq 0$%
, and that
\begin{equation}
\left\langle Qw,w\right\rangle \geq c_{Q}\Vert w\Vert _{-\theta
_{2}}^{2},\qquad w\in V^{-\theta _{2}},  \label{e:coercive-n}
\end{equation}%
with $c_{Q}>0$. We also assume that%
\begin{equation}
\left\langle A_{0}w,Qw\right\rangle \geq c_{A_{0}}\Vert w\Vert _{\theta
-\theta _{2}}^{2},\qquad w\in V^{\theta -\theta _{2}},
\label{e:coercive-abis}
\end{equation}%
Note that if $\theta =0$, (\ref{e:coercive-a}) is strictly speaking not
coercivity and follows from the boundedness of $A_{0}$, and note also that (%
\ref{e:coercive-n}) implies the invertibility of $Q$.

One may typically consider the following examples of operators occurring in
various combinations in (\ref{e:op}).

\begin{example}
\label{x:spaces}

(a) When $\Omega $ is a closed Riemannian manifold, and $E=T\Omega $ the
tangent bundle, an example of $\mathcal{V}$ is $\mathcal{V}_{\mathrm{per}%
}\subseteq \{u\in C^{\infty }(T\Omega ):\mathrm{div}\,\left( u\right) =0\}$,
a subspace of the divergence-free functions. The space of periodic functions
with vanishing mean on a torus $\mathbb{T}^{n}$ is a special case of this
example. In this case, one typically has $A_{0}=(-\Delta )^{\theta }$, $%
M=(I-\alpha ^{2}\Delta )^{-\theta _{1}}$, $Q=(I-\alpha ^{2}\Delta )^{-\theta
_{2}}$ and $A_{1}=-\Delta ,$ as operators that satisfy (\ref{e:bdd-amn}),
cf. \cite[Example 2.1, (a)]{HLT}.

(b) When $\Omega $ is a compact Riemannian manifold with boundary $\Gamma $
and again $E=T\Omega $ the tangent bundle, a typical example of $\mathcal{V}$
is $\mathcal{V}_{\mathrm{hom}}=\{u\in C_{0}^{\infty }(T\Omega ):\mathrm{div}%
\left( \,u\right) =0\}$. In this case, one may consider the choices $%
A_{0}=(-P\Delta )^{\theta }$, $A_{1}=-\Delta $, $M=(I-\alpha ^{2}P\Delta
)^{-\theta _{1}}$, and $Q=(I-\alpha ^{2}P\Delta )^{-\theta _{2}}$,
respectively, as operators satisfying (\ref{e:bdd-amn}), cf. \cite[Example
2.1, (b)]{HLT}.

(c) Let $\Omega $ be connected Riemannian $n$-dimensional manifold with
non-empty (sufficiently smooth) boundary $\partial \Omega .$ Define $%
A_{1}=-\Delta $, as the Laplacian of the metric $h$, acting on%
\begin{equation*}
D\left( A_{1}\right) =\left\{ \phi \in W^{2}:\phi =0\text{ on }\partial
\Omega \right\} ,
\end{equation*}%
where in local coordinates $\left\{ x_{i}\right\} _{i=1}^{n},$ the Laplacian
reads%
\begin{equation*}
\Delta \left( \cdot \right) =\frac{1}{\sqrt{\det \left( h\right) }}%
\sum\nolimits_{i,j=1}^{n}\partial _{x_{j}}\left( h^{ij}\sqrt{\det \left(
h\right) }\partial _{x_{i}}\left( \cdot \right) \right) ;
\end{equation*}%
the matrix $\left( h^{ij}\right) $ denotes the inverse of $h$. We have that $%
A_{1}$ is a positive self-adjoint operator on $W^{0}$.
\end{example}

\begin{example}
\label{y:spaces}In Example~\ref{x:spaces} above, the bilinear map $\bar{B}%
_{0}$ can be taken to be
\begin{equation}
\bar{B}_{0\chi }(v,w)=P[(v\cdot \nabla )w+\chi (\nabla w^{T})v],
\label{b00bis}
\end{equation}%
which correspond to the models with $\chi \in \left\{ 0,1\right\} $ as
introduced in the system (\ref{e:rel}).
\end{example}

To refer to the above examples, let us further introduce the shorthand
notation:%
\begin{equation}
B_{0\chi }\left( v,w\right) =\overline{B}_{0\chi }\left( Mv,Qw\right) ,\text{
}\chi \in \left\{ 0,1\right\} .  \label{b01}
\end{equation}%
For clarity, we list in Table \ref{t:spec2} the corresponding values of the
parameters and bilinear maps discussed above for special cases as given by (%
\ref{example}).{\small
\begin{table}[th]
\caption{Values of the parameters $\protect\theta ,$ $\protect\theta _{1}$
and $\protect\theta _{2}$, and the particular form of the bilinear map $%
B_{0} $ for some special cases of the model (\protect\ref{e:op}). (The
bilinear maps $B_{00}$ and $B_{01}$ are as in (\protect\ref{b01})). More
precisely, it allows us to include the special cases when $u$ satisfies the
Navier-Stokes equation (NSE), the Leray-$\protect\alpha $-equation, the
modified Leray-$\protect\alpha $-equation (ML), the simplified Bardina
model (SBM), the Navier-Stokes-Voigt equation (NSV)\ and the Lagrangian
averaged Navier-Stokes-$\protect\alpha $ model (NS-$\protect\alpha $-model).}
\label{t:spec2}
\begin{center}
{\small $%
\begin{tabular}{|c|c|c|c|c|c|c|}
\hline\hline
Model & NSE-EL & Leray-EL-$\alpha $ & ML--EL-$\alpha $ & SBM-EL & NSV-EL &
NS-EL-$\alpha $ \\ \hline
$\theta $ & 1 & 1 & 1 & 1 & 0 & 1 \\
$\theta _{1}$ & 0 & 1 & 0 & 1 & 1 & 0 \\
$\theta _{2}$ & 0 & 0 & 1 & 1 & 1 & 1 \\
$B_{0}$ & $B_{00}$ & $B_{00}$ & $B_{00}$ & $B_{00}$ & $B_{00}$ & $B_{01}$ \\
\hline\hline
\end{tabular}%
$  }
\end{center}
\end{table}
}

Next, we denote the trilinear forms
\begin{equation}
b_{0}(u,v,w)=\langle B_{0}(u,v),w\rangle ,\text{ }b_{1}\left( u,d,\psi
\right) =\left\langle B_{1}\left( u,d\right) ,\psi \right\rangle ,
\label{b01bis}
\end{equation}%
and similarly the forms $\bar{b}_{0\chi }$ and $b_{0\chi }$, following (\ref%
{b00}), (\ref{b00bis}) and (\ref{b01}). Then our notion of \emph{weak
solution} for problem (\ref{e:op}) can be formulated as follows.

\begin{definition}
\label{weak}\textit{Let }$g\left( t\right) \in L^{2}(0,T;V^{-s}),$ for some $%
s\in \mathbb{R}$ and%
\begin{equation*}
\left( u_{0},d_{0}\right) \in \mathcal{Y}_{\theta _{2}}:=\left\{
\begin{array}{ll}
V^{-\theta _{2}}\times W^{1}, & \text{if }\lambda _{2}\neq 0, \\
V^{-\theta _{2}}\times \left( W^{1}\cap \left\{ d_{0}\in L^{\infty }\left(
\Omega \right) \right\} :\left\Vert d_{0}\right\Vert _{L^{\infty }}\leq
1\right) , & \text{if }\lambda _{2}=0.%
\end{array}%
\right.
\end{equation*}%
\textit{Find a pair of functions}
\begin{equation}
\left( u,d\right) \in L^{\infty }\left( 0,T;\mathcal{Y}_{\theta _{2}}\right)
\cap L^{2}\left( 0,T;V^{\theta -\theta _{2}}\times W^{2}\right)  \label{1.8}
\end{equation}%
\textit{such that}%
\begin{equation}
\partial _{t}u\in L^{p}\left( 0,T;V^{-\gamma }\right) ,\text{ }\partial
_{t}d\in L^{2}\left( 0,T;W^{-2}\right)  \label{1.8bis}
\end{equation}%
for some $p>1$ and $\gamma \geq 0$, \textit{such that }$\left( u,d\right) $%
\textit{\ fulfills }$u\left( 0\right) =u_{0},$ $d\left( 0\right) =d_{0}$
\textit{and satisfies}%
\begin{align}
& \int\nolimits_{0}^{T}\left( -\left\langle u\left( t\right) ,w^{^{\prime
}}\left( t\right) \right\rangle +\left\langle A_{0}u\left( t\right) ,w\left(
t\right) \right\rangle +b_{0}\left( u\left( t\right) ,u\left( t\right)
,w\left( t\right) \right) \right) dt  \label{weak1} \\
& =\int\nolimits_{0}^{T}\left( \left\langle g\left( t\right) ,w\left(
t\right) \right\rangle +\left\langle R_{0}\left( A_{1}d\left( t\right)
,d\left( t\right) \right) ,w\left( t\right) \right\rangle -\left\langle
\sigma _{Q},\nabla w\left( t\right) \right\rangle \right) dt,  \notag
\end{align}%
\begin{align}
& \int\nolimits_{0}^{T}\left( -\left\langle d\left( t\right) ,\psi
^{^{\prime }}\left( t\right) \right\rangle +\left\langle \mu \left( t\right)
,\psi \left( t\right) \right\rangle +b_{1}\left( u\left( t\right) ,d\left(
t\right) ,\psi \left( t\right) \right) \right) dt  \label{weak2} \\
& =\int\nolimits_{0}^{T}\left( \left\langle \omega _{Q}d,\psi \left(
t\right) \right\rangle -\frac{\lambda _{2}}{\lambda _{1}}\left\langle
A_{Q}d,\psi \left( t\right) \right\rangle \right) dt,  \notag
\end{align}%
for any $\left( w,\psi \right) \in C_{0}^{\infty }\left( 0,T;\mathcal{V}%
\times \mathcal{W}\right) $, with $\mu \left( t\right) :=-\lambda
_{1}^{-1}\left( A_{1}d\left( t\right) +\nabla _{d}W(d\left( t\right)
)\right) $ a.e. on $\Omega \times \left( 0,T\right) .$
\end{definition}

\begin{remark}
As far as the interpretation of the initial conditions $u\left( 0\right)
=u_{0},$ $d\left( 0\right) =d_{0}$ is concerned, note that properties (\ref%
{1.8})-(\ref{1.8bis}) imply that $u\in C(0,T;V^{-\gamma })$ and $d\in
C(0,T;W^{0})$. Thus, the initial conditions are satisfied in a weak sense.
All kinds of boundary conditions (i.e., periodic, no-slip, no-flux, etc)\
for $\left( u,d\right) $ can be treated and will be included in our
analysis; they will be incorporated in the weak formulation for the problem (%
\ref{e:op}). On the other hand, for those values of $\left( \theta ,\theta
_{1},\theta _{2}\right) $ from Table \ref{t:spec2} we recover some specific
regularization models given by (\ref{e:rel}) for the particular choices of
the operators $A_{0},M,Q$ and $\chi $ in (\ref{example}), as listed in Table~%
\ref{t:spec2}.
\end{remark}

{\small
\begin{table}[th]
\caption{Some special cases of the model (\protect\ref{e:rel}) with $\protect%
\alpha >0$, and with $\Pi =(I-\protect\alpha ^{2}\Delta )^{-1}$.}
\label{t:spec}
\begin{center}
{\small $%
\begin{tabular}{|l||llllll|}
\hline
Model & \multicolumn{1}{||l|}{NSE-EL} & \multicolumn{1}{l|}{Leray-EL-$\alpha
$} & \multicolumn{1}{l|}{ML-EL-$\alpha $} & \multicolumn{1}{l|}{SBM-EL} &
\multicolumn{1}{l|}{NSV-EL} & NS-EL-$\alpha $ \\ \hline
$A_{0}$ & \multicolumn{1}{||l|}{$-\mu _{4}\Delta $} & \multicolumn{1}{l|}{$%
-\mu _{4}\Delta $} & \multicolumn{1}{l|}{$-\mu _{4}\Delta $} &
\multicolumn{1}{l|}{$-\mu _{4}\Delta $} & \multicolumn{1}{l|}{$-\mu
_{4}\Delta \Pi $} & $-\mu _{4}\Delta $ \\ \hline
$M$ & \multicolumn{1}{||l|}{$I$} & \multicolumn{1}{l|}{$\Pi $} &
\multicolumn{1}{l|}{$I$} & \multicolumn{1}{l|}{$\Pi $} & \multicolumn{1}{l|}{%
$\Pi $} & $\Pi $ \\ \hline
$Q$ & \multicolumn{1}{||l|}{$I$} & \multicolumn{1}{l|}{$I$} &
\multicolumn{1}{l|}{$\Pi $} & \multicolumn{1}{l|}{$\Pi $} &
\multicolumn{1}{l|}{$\Pi $} & $I$ \\ \hline
$\chi $ & \multicolumn{1}{||l|}{$0$} & \multicolumn{1}{l|}{$0$} &
\multicolumn{1}{l|}{$0$} & \multicolumn{1}{l|}{$0$} & \multicolumn{1}{l|}{$0$%
} & $1$ \\ \hline
\end{tabular}%
$  }
\end{center}
\end{table}
}

Throughout the paper, $C\geq 0$ will denote a \emph{generic} constant whose
further dependence on certain quantities will be specified on occurrence.
The value of the constant can change even within the same line. Furthermore,
we introduce the notation $a\lesssim b$ to mean that there exists a constant
$C>0$ such that $a\leq Cb.$ This notation will be used when the explicit
value of $C$ is irrelevant or tedious to write down.

\subsection{Energy estimates and solutions}

\label{ss:ee}

The system (\ref{e:op}) admits a \emph{total regularized energy}, consisting
of kinetic and potential energies, given by%
\begin{equation}
E_{Q}\left( t\right) =\frac{1}{2}\left[ \left\langle u\left( t\right)
,Qu\left( t\right) \right\rangle +||A_{1}^{1/2}d\left( t\right)
||_{L^{2}}^{2}\right] +\int_{\Omega }W\left( d\right) dx.
\label{total energy of the system}
\end{equation}%
In particular, for the smoothed systems introduced in (\ref{e:rel}), the
total energy $E_{Q}$ can be identified with the energy of the original
NSE-EL system under suitable boundary conditions. Furthermore, in the case
of the $\alpha $-models from Table \ref{t:spec2}, the invariant $E_{Q}$
reduces, as $\alpha \rightarrow 0$, to the dissipated energy $E_{I}$ of the
NSE-EL system. In order to show that $E_{Q}$ is an ideal invariant for the
system (\ref{e:op}), we need to perform some basic energy estimates and
computations. In what follows, we will always force the following

\noindent \textbf{Assumption on} $\mathcal{V}$: For a given smooth tensor $%
\Xi =\Xi \left( x\right) \in \mathbb{R}^{n\times n}$, we require that the
following identity holds:%
\begin{equation}
\left\langle div\left( \Xi \right) ,v\right\rangle +\left\langle \Xi ,\nabla
v\right\rangle =0,  \label{select-p}
\end{equation}%
for any $v=Qu\in \mathcal{V}$. In particular, such an assumption always
holds provided that $\mathcal{V}=\mathcal{V}_{per}$ is the space of periodic
(divergence-free) functions with vanishing mean on a torus $\Omega =\mathbb{T%
}^{n}$, see Example \ref{x:spaces}, (a). Clearly, (\ref{select-p}) will also
hold in function spaces $V^{s}=\mathrm{clos}_{H^{s}}\mathcal{V}$, $s\geq 1,$%
\ that have weaker differentiability properties. For more details on the
nature of this assumption, we refer the reader to Section \ref{s:cr}.

\noindent \textbf{Energy estimates}: In order to deduce a particular energy
identity, we will also assume that $b_{0}(u,u,Qu)=0,$ for any $Qu\in
\mathcal{V}$; pairing the first equation of (\ref{e:op}) with $Qu$ and the
second equation with $A_{1}d+\nabla _{d}W(d)$, respectively, by virtue of (%
\ref{select-p}) we deduce%
\begin{align}
& \frac{d}{dt}E_{Q}\left( t\right) +\left\langle A_{0}u,Qu\right\rangle -%
\frac{1}{\lambda _{1}}\left\Vert A_{1}d+\nabla _{d}W(d)\right\Vert
_{L^{2}}^{2}  \label{en-1} \\
& =-\left\langle \sigma _{Q},\nabla \left( Qu\right) \right\rangle
+\left\langle \omega _{Q}d,A_{1}d+\nabla _{d}W(d)\right\rangle -\frac{%
\lambda _{2}}{\lambda _{1}}\left\langle A_{Q}d,A_{1}d+\nabla
_{d}W(d)\right\rangle  \notag \\
& +\left\langle g\left( t\right) ,Qu\right\rangle ,  \notag
\end{align}%
for $t\in \left( 0,T\right) ,$ for any fixed but otherwise arbitrary $T>0$.
In order to simplify this identity further, we assume as in \cite{WXL12}
that the coefficients $\lambda _{1},\lambda _{2},\mu _{1},...,\mu _{6}$ obey
certain constraints:
\begin{eqnarray}
&&\lambda _{1}<0,  \label{lama1a} \\
&&\mu _{1}\geq 0,\quad \mu _{4}>0,  \label{mu14} \\
&&\mu _{5}+\mu _{6}\geq 0,  \label{mu56} \\
&&\lambda _{1}=\mu _{2}-\mu _{3},\ \ \ \lambda _{2}=\mu _{5}-\mu _{6}.
\label{lam2}
\end{eqnarray}%
Then, we insert the expression for the Leslie stress tensor $\sigma _{Q}$
from (\ref{v5}) and perform analogous computations as in \cite{LL01, WXL12},
relying on the symmetric properties of $A_{Q}$ and anti-symmetric properties
of $\omega _{Q}.$ We obtain after some lengthy but standard transformations
that%
\begin{align}
& \frac{d}{dt}E_{Q}\left( t\right) +\left\langle A_{0}u,Qu\right\rangle -%
\frac{1}{\lambda _{1}}\left\Vert A_{1}d+\nabla _{d}W(d)\right\Vert
_{L^{2}}^{2}+\mu _{1}\left\Vert d^{T}A_{Q}d\right\Vert _{L^{2}}^{2}
\label{en-2} \\
& =-\left( \mu _{2}+\mu _{3}\right) \left\langle d\otimes \mathcal{N}%
_{Q},A_{Q}\right\rangle -\left( \mu _{5}+\mu _{6}\right) \left\Vert
A_{Q}d\right\Vert _{L^{2}}^{2}  \notag \\
& -\frac{\lambda _{2}}{\lambda _{1}}\left\langle A_{Q}d,\lambda _{1}\mathcal{%
N}_{Q}+\lambda _{2}A_{Q}d\right\rangle +\left\langle g\left( t\right)
,Qu\right\rangle .  \notag
\end{align}%
We recall that (\ref{lama1a})-(\ref{lam2}) are always necessary in order for
the energy $E_{Q},$ $Q=I,$ of the system (\ref{e:rel}) to be nonincreasing
in the absence of external forces (cf. \cite{WXL12}). In addition, we'll
also assume two different sets of hypotheses on the coefficients according
to \cite{WXL12}.

\begin{itemize}
\item \textbf{Case 1} (with Parodi's relation). Suppose that \eqref{lama1a}--%
\eqref{lam2} are satisfied. Moreover, we enforce the following Parodi's
relation $\mu _{2}+\mu _{3}=\mu _{6}-\mu _{5}$ and%
\begin{equation}
\frac{(\lambda _{2})^{2}}{-\lambda _{1}}\leq \mu _{5}+\mu _{6}.
\label{critical point of lambda 2}
\end{equation}

\item \textbf{Case 2} (without Parodi's relation). Suppose that %
\eqref{lama1a}--\eqref{lam2} are satisfied. Moreover, we assume
\begin{equation}
|\lambda _{2}-\mu _{2}-\mu _{3}|<2\sqrt{-\lambda _{1}}\sqrt{\mu _{5}+\mu _{6}%
}.  \label{noPa1}
\end{equation}
\end{itemize}

In \textbf{Case 1}, it turns out that the regularized energy $E_{Q}$
satisfies for smooth solutions the identity%
\begin{align}
& \frac{d}{dt}E_{Q}\left( t\right) +\left\langle A_{0}u,Qu\right\rangle -%
\frac{1}{\lambda _{1}}\left\Vert A_{1}d+\nabla _{d}W(d)\right\Vert
_{L^{2}}^{2}+\mu _{1}\left\Vert d^{T}A_{Q}d\right\Vert _{L^{2}}^{2}
\label{en-case1} \\
& =-\left( \mu _{5}+\mu _{6}+\frac{\lambda _{2}^{2}}{\lambda _{1}}\right)
\left\Vert A_{Q}d\right\Vert _{L^{2}}^{2}+\left\langle g\left( t\right)
,Qu\right\rangle ,  \notag
\end{align}%
while in the \textbf{Case 2}, it obeys%
\begin{align}
& \frac{d}{dt}E_{Q}\left( t\right) +\left\langle A_{0}u,Qu\right\rangle +\mu
_{1}\left\Vert d^{T}A_{Q}d\right\Vert _{L^{2}}^{2}-\frac{3\lambda _{1}}{4}%
\left\Vert \mathcal{N}_{Q}\right\Vert _{L^{2}}^{2}  \label{en-case2} \\
& \leq -\left( \mu _{5}+\mu _{6}+\frac{\left( \lambda _{2}-\left( \mu
_{2}+\mu _{3}\right) \right) ^{2}}{\lambda _{1}}\right) \left\Vert
A_{Q}d\right\Vert _{L^{2}}^{2}+\left\langle g\left( t\right)
,Qu\right\rangle .  \notag
\end{align}%
Next, for every $\varepsilon >0$ we have%
\begin{equation*}
\langle g\left( t\right) ,Qu\rangle \leq \varepsilon ^{-1}\Vert g\left(
t\right) \Vert _{-\theta -\theta _{2}}^{2}+\varepsilon \Vert Q\Vert
_{-\theta _{2};\theta _{2}}^{2}\Vert u\Vert _{\theta -\theta _{2}}^{2}.
\end{equation*}%
Employing now the condition (\ref{e:coercive-abis}) for the operator $A_{0}$%
, we can absorb this term on the right-hand side of (\ref{en-case1})-(\ref%
{en-case2}). In either case, by integrating the resulting relation on the
time interval $\left( s,t\right) ,$ we easily derive that $E_{Q}$ satisfies
an energy inequality. More precisely, owing to (\ref{e:coercive-abis}) there
holds a.e. $t>0,$%
\begin{align}
\frac{d}{dt}& E_{Q}\left( t\right) +\frac{c_{A_{0}}}{2}\left\Vert u\left(
t\right) \right\Vert _{\theta -\theta _{2}}^{2}-\frac{1}{\lambda _{1}}%
\left\Vert A_{1}d\left( t\right) +\nabla _{d}W\left( d\left( t\right)
\right) \right\Vert _{L^{2}}^{2}+\mu _{1}\left\Vert \left(
d^{T}A_{Q}d\right) \left( t\right) \right\Vert _{L^{2}}^{2}  \label{en-both}
\\
& \leq \frac{2\Vert Q\Vert _{-\theta _{2};\theta _{2}}^{2}}{c_{A_{0}}}%
\left\Vert g\left( t\right) \right\Vert _{-\theta -\theta _{2}}^{2}.  \notag
\end{align}%
It follows from (\ref{en-both}) by integration over $\left( 0,t\right) $
that $\left( u\left( t\right) ,d\left( t\right) \right) $ belongs to the
functional class (\ref{1.8}) given $g\in L^{2}\left( 0,T;V^{-\theta -\theta
_{2}}\right) $. We note that $E_{Q}\left( 0\right) <\infty $ is equivalent
to having $\left( u_{0},d_{0}\right) \in \mathcal{Y}_{\theta _{2}}$.

We now introduce another notion of weak solutions which is also essential in
our subsequent study.

\begin{definition}
\label{weak-energy}Let $\lambda _{1},\lambda _{2},\mu _{1},...,\mu _{6}$
satisfy the above assumptions according to the \textbf{Cases 1-2}. By an
\emph{energy solution} we will mean a weak distributional solution $\left(
u,d\right) $, satisfying the weak formulation (\ref{weak1})-(\ref{weak2})
and obeying the energy inequality (\ref{en-both}) according to the Cases 1
and 2, respectively.
\end{definition}

It is worth pointing out that, by virtue of (\ref{en-case1})-(\ref{en-case2}%
), energy solutions of the regularized Ericksen-Leslie system (\ref{e:op})
satisfy:%
\begin{equation}
\left\{
\begin{array}{l}
A_{1}d+\nabla _{d}W(d)\in L^{2}\left( 0,T;L^{2}\left( \Omega \right) \right)
,\text{ }A_{Q}d\in L^{2}\left( 0,T;L^{2}\left( \Omega \right) \right) , \\
\mathcal{N}_{Q}\in L^{2}\left( 0,T;L^{2}\left( \Omega \right) \right) ,\text{
}d^{T}A_{Q}d\in L^{2}\left( 0,T;L^{2}\left( \Omega \right) \right) .%
\end{array}%
\right.  \label{apriori-tensor}
\end{equation}%
Indeed, the energy dissipations provided by the inequalities in (\ref%
{en-case1}) and (\ref{en-case2}) are equivalent because of definitions (\ref%
{v4}), (\ref{def1}) and the equation for the director field $d$ from (\ref%
{e:op}) (see also \cite{WXL12}). Such knowledge will also become important
in the study of global regularity.

\section{Well-posedness results}

\label{s:well}

Analogous to the theory of regularized flows we have developed for a
simplified Ericksen-Leslie model in \cite{GM-JNS}, we begin to devise a
solution theory for the general three-parameter family of regularized models
from (\ref{e:rel}). We begin to establish existence and regularity results,
and under appropriate assumptions uniqueness and stability. At the end of
the proof of each theorem, we give the corresponding conditions for $\left(
\theta ,\theta _{1},\theta _{2}\right) $ which allow us not only to
establish old results but also new results in the literature, especially for
the cases listed in Table \ref{t:spec2}. The analysis in this section is
divided mainly into two parts according to whether $\lambda _{2}=0$ or $%
\lambda _{2}\neq 0$.

\subsection{Existence of weak solutions}

\label{ss:exist}

In this subsection, we establish sufficient conditions for the existence of
energy solutions to the problem (\ref{e:op}) (cf. Definition \ref%
{weak-energy}). As noted previously, in the case when $\lambda _{2}=0,$ a
maximum principle holds for the director field $d$ of any weak solution.

\begin{proposition}
\label{maxp}Suppose that $b_{1}\left( v,\psi ,\psi \right) =0$, for any $%
v\in V^{\theta -\theta _{2}},$ $\psi \in W^{1}$. Let $d_{0}\in L^{\infty
}\left( \Omega \right) $ such that $\left\Vert d_{0}\right\Vert _{L^{\infty
}\left( \Omega \right) }\leq 1.$ Then, for any weak solution $\left(
u,d\right) $ to problem (\ref{e:op}) in the sense of Definition \ref{weak},
we have $d\in L^{\infty }\left( 0,T;L^{\infty }\left( \Omega \right) \right)
$ and%
\begin{equation}
\left\vert d\left( x,t\right) \right\vert \leq \left\Vert d_{0}\right\Vert
_{L^{\infty }\left( \Omega \right) }\text{, a.e. on }\Omega \times \left(
0,T\right) .  \label{6.3}
\end{equation}
\end{proposition}

\begin{proof}
The inequality in (\ref{6.3}) follows from a straightforward application of
the weak maximum principle, since $\lambda _{2}=0$ and the tensor $\omega
_{Q}$ is skew-symmetric. A Moser type of iteration argument also gives the
desired regularity $d\in L^{\infty }\left( 0,T;L^{\infty }\left( \Omega
\right) \right) $ (see, e.g., \cite[Lemma 9.3.1]{CD}).
\end{proof}

It is worth emphasizing that when $\lambda _{2}\neq 0$, the inequality (\ref%
{6.3}) is generally not expected to hold.

\begin{theorem}
\label{t:exist} Assume Proposition \ref{maxp} \emph{only} when $\lambda
_{2}=0$ and let the following conditions hold.

\begin{itemize}
\item[i)] $\left( u_{0},d_{0}\right) \in \mathcal{Y}_{\theta _{2}}$ with any
$\theta _{2}\geq 0$ and $g\in L^{2}(0,T;V^{-\theta -\theta _{2}})$, $T>0$.

\item[ii)] $b_{0}(v,v,Qv)=0,$ for any $v\in V^{\theta -\theta _{2}}$;

\item[iii)] $b_{0}:V^{\bar{\sigma}_{1}}\times V^{\bar{\sigma}_{2}}\times V^{%
\bar{\gamma}}\rightarrow \mathbb{R}$ is bounded for some $\bar{\sigma}%
_{i}<\theta -\theta _{2}$, $i=1,2$, and $\bar{\gamma}\geq \gamma $;

\item[iv)] $b_{0}:V^{\sigma _{1}}\times V^{\sigma _{2}}\times V^{\gamma
}\rightarrow \mathbb{R}$ is bounded for some $\sigma _{i}\in \lbrack -\theta
_{2},\theta -\theta _{2}]$, $i=1,2$, and $\gamma \in \lbrack \theta +\theta
_{2},\infty )\cap (\theta _{2},\infty )\cap (1+n/6,\infty )$;

Then, there exists at least one energy solution $\left( u,d\right) $
satisfying (\ref{1.8})-(\ref{1.8bis}), (\ref{apriori-tensor}) such that%
\begin{equation}
p=\left\{
\begin{array}{ll}
\min \{2,\frac{2\theta }{\sigma _{1}+\sigma _{2}+2\theta _{2}},\frac{4\left(
6-n\right) }{12-n}\}, & \text{if }\theta >0, \\
\frac{4\left( 6-n\right) }{12-n}, & \text{if }\theta =0.%
\end{array}%
\right.  \label{thep}
\end{equation}
\end{itemize}
\end{theorem}

\begin{proof}
We rely on a Galerkin approximation scheme by borrowing ideas from \cite{HLT}%
. To this end, let $\{V_{m}:m\in \mathbb{N}\}\subset V^{\theta -\theta
_{2}}, $ $\{W_{m}:m\in \mathbb{N}\}\subset D\left( A_{1}\right) \cap
L^{\infty }\left( \Omega \right) $ be sequences of finite dimensional
(smooth)\ subspaces of $V^{\theta -\theta _{2}}$ and $D\left( A_{1}\right) ,$
respectively, such that

\begin{enumerate}
\item $V_{m}\subset V_{m+1},$ $W_{m}\subset W_{m+1},$ for all $m\in \mathbb{N%
}$;

\item $\cup _{m\in \mathbb{N}}V_{m}$ is dense in $V^{\theta -\theta _{2}}$,
and $\cup _{m\in \mathbb{N}}W_{m}$ is dense in $D\left( A_{1}\right) ;$

\item For $m\in \mathbb{N}$, with $\widetilde{V}_{m}=QV_{m}\subset V^{\theta
+\theta _{2}}$, the projectors $P_{m}:V^{\theta -\theta _{2}}\rightarrow
V_{m},$ $S_{m}:D\left( A_{1}\right) \rightarrow W^{m}$, defined by
\begin{align*}
\left\langle P_{m}v,w_{m}\right\rangle & =\left\langle v,w_{m}\right\rangle
,\qquad w_{m}\in \widetilde{V}_{m},\,v\in V^{\theta -\theta _{2}}, \\
\left\langle S_{m}d,\psi _{m}\right\rangle & =\left\langle d,\psi
_{m}\right\rangle ,\qquad \psi _{m}\in W_{m},\,d\in D\left( A_{1}\right) ,
\end{align*}%
are uniformly bounded as maps from $V^{-\gamma }\rightarrow V^{-\gamma }$
and $W^{-2}\rightarrow W^{-2}$, respectively.
\end{enumerate}

Such sequences can be constructed e.g., by using the eigenfunctions of the
isometries $\Lambda ^{1+\theta }:V^{1+\theta -\theta _{2}}\rightarrow
V^{-\theta _{2}}$, $A_{1}:D\left( A_{1}\right) \rightarrow W^{0}.$ Consider
the problem of finding $\left( u_{m},d_{m}\right) \in C^{1}(0,T;V_{m}\times
W_{m})$ such that for all $\left( w_{m},\psi _{m}\right) \in \widetilde{V}%
_{m}\times W_{m}$,%
\begin{equation}
\left\{
\begin{array}{l}
\langle \partial _{t}u_{m},w_{m}\rangle +\langle A_{0}u_{m},w_{m}\rangle
+b_{0}(u_{m},u_{m},w_{m}) \\
=\langle g,w_{m}\rangle +\left\langle R_{0}\left( A_{1}d_{m},d_{m}\right)
+div\left( \sigma _{Q}^{m}\right) ,w_{m}\right\rangle , \\
\left\langle \partial _{t}d_{m},\psi _{m}\right\rangle +b_{1}\left(
u_{m},d_{m},\psi _{m}\right) +\left\langle \omega _{Q}^{m}d_{m}+\frac{%
\lambda _{2}}{\lambda _{1}}A_{Q}^{m}d_{m},\psi _{m}\right\rangle \\
=\frac{1}{\lambda _{1}}\left\langle A_{1}d_{m}+\nabla _{d}W\left(
d_{m}\right) ,\psi _{m}\right\rangle , \\
\langle u_{m}(0),w_{m}\rangle =\langle u_{0},w_{m}\rangle , \\
\left\langle d_{m}\left( 0\right) ,\psi _{m}\right\rangle =\left\langle
d_{0},\psi _{m}\right\rangle ,%
\end{array}%
\right.  \label{e:galerkin}
\end{equation}%
where%
\begin{align*}
A_{Q}^{m}& =\frac{1}{2}(\nabla \left( Qu_{m}\right) +\nabla ^{T}\left(
Qu_{m}\right) ),\ \,\;\omega _{Q}^{m}=\frac{1}{2}(\nabla \left(
Qu_{m}\right) -\nabla ^{T}\left( Qu_{m}\right) ), \\
\dot{d}_{m}& =\partial _{t}d_{m}+Qu_{m}\cdot \nabla d_{m},\ \,\;\mathcal{N}%
_{Q}^{m}=\dot{d}_{m}-\omega _{Q}^{m}d_{m}
\end{align*}%
and%
\begin{equation*}
\boldsymbol{\sigma }_{Q}^{m}=\mu _{1}(d_{m}^{T}A_{Q}^{m}d_{m})d_{m}\otimes
d_{m}+\mu _{2}\mathcal{N}_{Q}^{m}\otimes d_{m}+\mu _{3}d_{m}\otimes \mathcal{%
N}_{Q}^{m}+\mu _{5}(A_{Q}^{m}d_{m})\otimes d_{m}+\mu _{6}d_{m}\otimes
(A_{Q}^{m}d_{m}).
\end{equation*}%
Choosing a basis for $V_{m}\times W_{m}$, one sees that the system (\ref%
{e:galerkin}) is an initial value problem for a system of ODE's. By
definition, the operator $Q$ is invertible so that the standard ODE theory
gives a unique solution to (\ref{e:galerkin}), which is locally-defined in
time. Using the definition of $Q$ once more, one checks that%
\begin{equation*}
c_{Q}\Vert u_{m}(0)\Vert _{-\theta _{2}}^{2}\leq \langle
u_{m}(0),Qu_{m}(0)\rangle =\langle u(0),Qu_{m}(0)\rangle \leq \Vert
u(0)\Vert _{-\theta _{2}}\Vert Qu_{m}(0)\Vert _{\theta _{2}},
\end{equation*}%
so that $\Vert u_{m}(0)\Vert _{-\theta _{2}}$ is uniformly bounded.
Similarly, one shows that $\left\Vert d_{m}\left( 0\right) \right\Vert _{1}$
is uniformly bounded.

Now in the first and second equalities of (\ref{e:galerkin}), taking $%
w_{m}=Qu_{m}$ and $\psi _{m}=A_{1}d_{m}+\nabla _{d}W\left( d_{m}\right) $,
respectively, and using the a priori estimates established earlier in (\ref%
{en-case1})-(\ref{apriori-tensor}), one derives that the solution $\left(
u_{m},d_{m}\right) $ is uniformly bounded in $L^{\infty }(0,T;\mathcal{Y}%
_{\theta _{2}})\cap L^{2}(0,T;V^{\theta -\theta _{2}}\times D\left(
A_{1}\right) )$. Moreover, the terms $A_{1}d_{m}+\nabla _{d}W(d_{m})$, $%
A_{Q}^{m}d_{m},$ $\mathcal{N}_{Q}^{m}$ and $d_{m}^{T}A_{Q}^{m}d_{m}$ are
also uniformly bounded in $L^{2}\left( 0,T;L^{2}\left( \Omega \right)
\right) $. Thus, passing to a subsequence, one has%
\begin{equation}
\left\{
\begin{array}{l}
\left( u_{m},d_{m}\right) \rightarrow \left( u,d\right)
\mbox{
weak-star in }L^{\infty }(0,T;\mathcal{Y}_{\theta _{2}}), \\
\left( u_{m},d_{m}\right) \rightarrow \left( u,d\right) \mbox{
weakly in }L^{2}(0,T;V^{\theta -\theta _{2}}\times D\left( A_{1}\right) ).%
\end{array}%
\right.  \label{e:exist-weak-conv}
\end{equation}

Passing to the limit as $m\rightarrow \infty $ in (\ref{e:galerkin})
requires the use of compactness arguments. To this end, we start by
estimating $\left\Vert \partial _{t}u_{m}\right\Vert _{-\gamma }$ and $%
\left\Vert \partial _{t}d_{m}\right\Vert _{-2}$, respectively. The first
equation in (\ref{e:galerkin}) may be recast as:%
\begin{equation}
\partial _{t}u_{m}+P_{m}A_{0}u_{m}+P_{m}B_{0}(u_{m},u_{m})=P_{m}\left(
g+R_{0}\left( A_{1}d_{m},d_{m}\right) +div(\sigma _{Q}^{m})\right) .
\label{e:galerkin-abs}
\end{equation}%
Consequently,%
\begin{align}
\Vert \partial _{t}u_{m}\Vert _{-\gamma }& \leq \Vert P_{m}A_{0}u_{m}\Vert
_{-\gamma }+\Vert P_{m}B_{0}(u_{m},u_{m})\Vert _{-\gamma }+\Vert P_{m}g\Vert
_{-\gamma }  \label{exist-1} \\
& +\Vert P_{m}R_{0}\left( A_{1}d_{m},d_{m}\right) \Vert _{-\gamma }+\Vert
P_{m}div(\sigma _{Q}^{m})\Vert _{-\gamma }  \notag \\
& \lesssim \Vert u_{m}\Vert _{\theta -\theta _{2}}+\Vert
P_{m}B_{0}(u_{m},u_{m})\Vert _{-\gamma }+\Vert g\Vert _{-\theta -\theta _{2}}
\notag \\
& +\Vert P_{m}R_{0}\left( A_{1}d_{m},d_{m}\right) \Vert _{-\gamma }+\Vert
P_{m}div(\sigma _{Q}^{m})\Vert _{-\gamma }.  \notag
\end{align}%
Now, thanks to the boundedness of $B_{0}$ (see (iv)), it follows as in \cite[%
Theorem 3.1]{HLT} that%
\begin{equation}
\Vert P_{m}B_{0}(u_{m},u_{m})\Vert _{-\gamma }\lesssim \Vert u_{m}\Vert
_{\sigma _{1}}\Vert u_{m}\Vert _{\sigma _{2}}.  \label{exist0}
\end{equation}%
If $\theta =0$, then the norms in the right hand side are the $V^{-\theta
_{2}}$-norm which is uniformly bounded. On the other hand, if $\theta >0$,
then by interpolation, one gets
\begin{equation}
\Vert u_{m}\Vert _{\sigma _{i}}\lesssim \Vert u_{m}\Vert _{-\theta
_{2}}^{1-\lambda _{i}}\Vert u_{m}\Vert _{\theta -\theta _{2}}^{\lambda
_{i}},\qquad \lambda _{i}=\frac{\sigma _{i}+\theta _{2}}{\theta },\quad
i=1,2,  \label{exist1}
\end{equation}%
so that

\begin{equation}
\Vert P_{m}B_{0}(u_{m},u_{m})\Vert _{-\gamma }\lesssim \Vert u_{m}\Vert
_{-\theta _{2}}^{2-\lambda _{1}-\lambda _{2}}\Vert u_{m}\Vert _{\theta
-\theta _{2}}^{\lambda _{1}+\lambda _{2}}\lesssim \Vert u_{m}\Vert _{\theta
-\theta _{2}}^{\lambda _{1}+\lambda _{2}}.  \label{exist2}
\end{equation}%
Hence, with $\lambda :=\lambda _{1}+\lambda _{2}=\frac{\sigma _{1}+\sigma
_{2}+2\theta _{2}}{\theta }$ if $\theta >0$, and with $\lambda =1$ if $%
\theta =0$, we get
\begin{equation}
\int_0^T \Vert P_{m}B_{0}(u_{m},u_{m})\Vert _{-\gamma }^{p}\,dt\lesssim
\Vert u_{m}\Vert _{L^{p}(V^{\theta -\theta _{2}})}^{p}+\Vert u_{m}\Vert
_{L^{p\lambda }(V^{\theta -\theta _{2}})}^{p}.  \label{exist2bis}
\end{equation}%
The first term on the right-hand side is bounded uniformly when $p\leq 2$.
The second term is bounded if $p\lambda \leq 2$, that is $p\leq 2/\lambda $.
We conclude that $P_{m}B_{0}(u_{m},u_{m})$ is uniformly bounded in $%
L^{p}\left( V^{-\gamma }\right) $, with $p=\min \{2,2/\lambda \}$.
Concerning a uniform bound for $P_{m}R_{0}\left( A_{1}d_{m},d_{m}\right) $
in $L^{2}\left( V^{-\gamma }\right) $, we argue as in \cite[Theorem 3.2]%
{GM-JNS} to derive that%
\begin{equation}
\Vert P_{m}R_{0}\left( A_{1}d_{m},d_{m}\right) \Vert _{-\gamma }^{2}\lesssim
\left\Vert A_{1}d_{m}\right\Vert _{L^{2}}^{2}\left\Vert d_{m}\right\Vert
_{1}^{2}  \label{exist2biss}
\end{equation}%
provided that $\gamma >1+\frac{n}{6}\geq \frac{n}{2}$. As for the remaining
term in (\ref{exist-1}), one has:%
\begin{align}
\Vert P_{m}div(\sigma _{Q}^{m})\Vert _{-\gamma }& \lesssim \mu
_{1}\left\Vert div\left( d_{m}^{T}A_{Q}^{m}d_{m})d_{m}\otimes d_{m}\right)
\right\Vert _{-\gamma }  \label{exist-3} \\
& +\left\Vert div\left( \mu _{2}\mathcal{N}_{Q}^{m}\otimes d_{m}+\mu
_{3}d_{m}\otimes \mathcal{N}_{Q}^{m}\right) \right\Vert _{-\gamma }  \notag
\\
& +\left\Vert div\left( \mu _{5}(A_{Q}^{m}d_{m})\otimes d_{m}+\mu
_{6}d_{m}\otimes (A_{Q}^{m}d_{m})\right) \right\Vert _{-\gamma }  \notag \\
& =:I_{1}+I_{2}+I_{3}.  \notag
\end{align}%
We'll just estimate $I_{1}$; estimating $I_{2}$ and $I_{3}$ follows suit. To
this end, let $\varphi _{m}\in V^{\gamma }\subset W^{1,3}$ with $\left\Vert
\varphi _{m}\right\Vert _{\gamma }=1$ and $\gamma >\frac{n}{6}+1$. Then%
\begin{align}
\left\vert \left\langle div\left( d_{m}^{T}A_{Q}^{m}d_{m})d_{m}\otimes
d_{m}\right) ,\varphi _{m}\right\rangle \right\vert & =\left\vert
\left\langle \left( d_{m}^{T}A_{Q}^{m}d_{m}\right) d_{m}\otimes d_{m},\nabla
\varphi _{m}\right\rangle \right\vert  \label{exist-3bis} \\
& \lesssim \left\Vert d_{m}^{T}A_{Q}^{m}d_{m}\right\Vert _{L^{2}}\left\Vert
A_{1}^{1/2}d_{m}\right\Vert _{L^{2}}^{2\left( 1-\delta _{n}\right)
}\left\Vert A_{1}d_{m}\right\Vert _{L^{2}}^{2\delta _{n}}  \notag
\end{align}%
with $\delta _{n}=\frac{n}{4\left( 6-n\right) }$. Similarly, we have%
\begin{equation}
I_{2}\lesssim \left\Vert \mathcal{N}_{Q}^{m}\right\Vert _{L^{2}}\left\Vert
d_{m}\right\Vert _{1}\text{ and }I_{3}\lesssim \left\Vert
A_{Q}^{m}d_{m}\right\Vert _{L^{2}}\left\Vert d_{m}\right\Vert _{1}.
\label{exist-3biss}
\end{equation}%
Consequently, on account of the estimates (\ref{exist2bis})-(\ref%
{exist-3biss}) and thanks to H\"{o}lder's inequality, it follows from (\ref%
{exist-1}) that $\partial _{t}u_{m}$ is uniformly bounded in $L^{p}\left(
0,T;V^{-\gamma }\right) $, provided that $p\leq 2,$ $p\leq \frac{4\left(
6-n\right) }{12-n}$ and $p\leq 2\lambda $ if $\theta >0$, and $p\leq \frac{%
4\left( 6-n\right) }{12-n}$ for $\theta =0$. To estimate $\partial _{t}\phi
_{m},$ we recast the second equation in (\ref{e:galerkin}) as%
\begin{align*}
& \partial _{t}d_{m}+S_{m}B_{1}\left( u_{m},d_{m}\right) +S_{m}\omega
_{Q}^{m}d_{m}+\frac{\lambda _{2}}{\lambda _{1}}S_{m}A_{Q}^{m}d_{m} \\
& =\frac{1}{\lambda _{1}}S_{m}A_{1}d_{m}+S_{m}\nabla _{d}W\left(
d_{m}\right) .
\end{align*}%
It follows from the uniform boundedness of $S_{m}$ that%
\begin{align}
\left\Vert \partial _{t}d_{m}\right\Vert _{-2}& \lesssim \left\Vert
B_{1}\left( u_{m},d_{m}\right) \right\Vert _{-2}+\left\Vert \omega
_{Q}^{m}d_{m}\right\Vert _{-2}+\left\Vert A_{Q}^{m}d_{m}\right\Vert _{-2}
\label{exist-4} \\
& +\left\Vert A_{1}d_{m}+\nabla _{d}W\left( d_{m}\right) \right\Vert _{-2}.
\notag
\end{align}%
Thanks to the Hahn-Banach theorem, H\"{o}lder's inequality and a proper
Sobolev embedding theorem, we can argue as in \cite[Theorem 3.2]{GM-JNS} to
get the following estimates:%
\begin{align}
\left\Vert B_{1}\left( u_{m},d_{m}\right) \right\Vert _{-2}& =\left\langle
Qu_{m}\cdot \nabla d_{m},\varphi _{m}\right\rangle \text{, }\varphi _{m}\in
V^{2},\text{ }\left\Vert \varphi _{m}\right\Vert _{2}=1  \label{exist-4bis}
\\
& \lesssim \left\Vert Qu_{m}\right\Vert _{\theta _{2}}\left\Vert \nabla
d_{m}\right\Vert _{1}\left\Vert \varphi _{m}\right\Vert _{2}  \notag \\
& \lesssim \left\Vert u_{m}\right\Vert _{-\theta _{2}}\left\Vert
A_{1}d_{m}\right\Vert _{L^{2}},\text{ since }\theta _{2}\geq 0,  \notag
\end{align}%
and, using Einstein's summation convention,%
\begin{align}
\left\Vert \omega _{Q}^{m}d_{m}\right\Vert _{-2}& =\left\langle \omega
_{Q}^{m}d_{m},\psi _{m}\right\rangle ,\text{ }\psi _{m}\in V^{2},\text{ }%
\left\Vert \psi _{m}\right\Vert _{2}=1  \label{exist-4tris} \\
& =\left\langle \left( \omega _{Q}^{m}\right) _{ij}d_{mj},\psi
_{mi}\right\rangle =\left\langle \left( \omega _{Q}^{m}\right)
_{ij},d_{mj}\psi _{mi}\right\rangle  \notag \\
& \lesssim \left\Vert \left( \omega _{Q}^{m}\right) _{ij}\right\Vert
_{-1}\left\Vert d_{mj}\psi _{mi}\right\Vert _{1}  \notag \\
& \lesssim \left\Vert \left( \omega _{Q}^{m}\right) _{ij}\right\Vert
_{-1}\left\Vert d_{m}\right\Vert _{2}.  \notag
\end{align}%
Now, since $\theta _{2}\geq 0$, it holds%
\begin{equation*}
\left\Vert \left( \omega _{Q}^{m}\right) _{ij}\right\Vert _{-1}\lesssim
\left\Vert \left( \nabla Qu_{m}\right) _{ij}\right\Vert _{-1}\lesssim
\left\Vert Qu_{m}\right\Vert _{0}\lesssim \left\Vert u_{m}\right\Vert
_{-\theta _{2}}.
\end{equation*}%
Substituting this bound into (\ref{exist-4tris}), we easily arrive at the
bound%
\begin{equation}
\left\Vert \omega _{Q}^{m}d_{m}\right\Vert _{-2}\lesssim \left\Vert
u_{m}\right\Vert _{-\theta _{2}}\left\Vert A_{1}d_{m}\right\Vert _{L^{2}},
\label{exist-4quad}
\end{equation}%
It follows from (\ref{exist-4})-(\ref{exist-4quad}) and earlier estimates
that $\partial _{t}d_{m}$ is uniformly bounded in $L^{2}\left(
0,T;V^{-2}\right) $. With the estimates for $\partial _{t}u_{m}$ and $%
\partial _{t}d_{m}$, we now have the required ingredients for the
application of the Aubin-Lions-Simon compactness theorem (see, e.g., \cite[%
Appendix]{HLT}). In particular, we can infer the existence of a limit couple%
\begin{equation*}
\left( u,d\right) \in C(0,T;V^{-\gamma }\times W^{0})\cap L^{\infty }(0,T;%
\mathcal{Y}_{\theta _{2}})
\end{equation*}%
such that, in addition to (\ref{e:exist-weak-conv}), we also have%
\begin{equation}
\left\{
\begin{array}{l}
\left( u_{m},d_{m}\right) \rightarrow \left( u,d \right)
\mbox{
strongly in }L^{2}(0,T;V^{s}\times W^{2-}) \\
d_{m}\rightarrow d\mbox{ strongly in }C(0,T;W^{1-}),%
\end{array}%
\right.  \label{e:exist-conv-u}
\end{equation}%
for any $s<\theta -\theta _{2},$ where $W^{s-}$ denotes $W^{s-\delta },$ for
some sufficiently small $\delta \in (0,s]$.

We are now able to pass to the limit in all the nonlinear terms of (\ref%
{e:galerkin}) so that this limit couple $\left( u,d\right) $ indeed
satisfies the weak formulation (\ref{weak1})-(\ref{weak2}) of Definition \ref%
{weak}. This is standard procedure and so we leave the details to the
interested reader. However, we refer the reader to \cite[Theorem 3.1]{HLT}
for passage to the limit in the equation for the velocity and to \cite[%
Theorem 3.2]{GM-JNS} for passage to the limit in the elastic (Ericksen)
stress tensor $R_{0}$. The proof of the theorem is now finished.
\end{proof}

Our theorem covers the following special cases listed in Table \ref{t:spec2}.

\begin{remark}
\label{rem:all-mod}Let $\theta +\theta _{1}>\frac{1}{2}$ and recall that $%
\theta ,\theta _{2}\geq 0$. By \cite[Proposition 2.5]{HLT}, the trilinear
form $b_{00},$ defined by (\ref{b01})-(\ref{b01bis}), fulfills the
hypotheses (ii)-(iv) of Theorem \ref{t:exist} for $-\gamma \leq \theta
-\theta _{2}-1$ with $-\gamma <\min \{2\theta +2\theta _{1}-\frac{n+2}{2}%
,\theta -\theta _{2}+2\theta _{1},\theta +\theta _{2}-1\}$. Similarly, the
trilinear form $b_{01}$ satisfies (ii)-(iv) for $-\gamma \leq \theta -\theta
_{2}-1$ with $-\gamma <\min \{2\theta +2\theta _{1}-\frac{n+2}{2},\theta
-\theta _{2}+2\theta _{1}-1,\theta +\theta _{2}\}$. In particular, our
result yields the global existence of a weak energy solution for both the
inviscid and viscous Leray-EL-$\alpha $ models in three space dimensions,
and for all the other regularized models listed in Table \ref{t:spec2}.\ As
far as we know, except for the 3D NSE-EL system reported in \cite{CRW}, none
of these results have been reported previously.
\end{remark}

\begin{remark}
As in \cite[Section 4]{GM-JNS} for the simplified Ericksen-Leslie model ($%
\sigma _{Q}\equiv 0$, $\omega _{Q}\equiv 0,$ $\lambda _{2}=0$), it is also
possible to consider the situation where the operators $A_{0}$ and $B_{0}$
in the general three-parameter family of regularized models represented by
problem (\ref{e:op}) have values from a convergent (in a certain sense)
sequence, and study the limiting behavior of the corresponding sequence of
energy solutions. As a special case this includes the $\alpha \rightarrow
0^{+}$ limits in the $\alpha $-models (\ref{e:rel}). We leave the details
for future contributions.
\end{remark}

\subsection{Regularity of weak solutions}

\label{ss:reg}

In this subsection, we develop a regularity result for the energy solutions
of the general family of regularized models constructed in Section \ref%
{ss:exist}. Recall that $\theta ,\theta _{2}\geq 0$, $\theta _{1}\in \mathbb{%
R}$ and that, in general, $\lambda _{2}\neq 0$.

\begin{theorem}
\label{t:reg} Let%
\begin{equation*}
\left( u,d\right) \in L^{\infty }\left( 0,T;\mathcal{Y}_{\theta _{2}}\right)
\cap L^{2}\left( 0,T;V^{\theta -\theta _{2}}\times D\left( A_{1}\right)
\right)
\end{equation*}%
be an energy solution in the sense of Definition \ref{weak}. Let $s\in
\left( \frac{n}{4},1\right] ,$ $n=2,3$ and consider the following nonempty
interval%
\begin{equation*}
J_{n}:=\left( -\theta _{2},\theta -\frac{n}{2}\right)\cap \left[ 4s-\theta
-3\theta _{2},+\infty \right) .
\end{equation*}%
For $\beta \in J_{n}\neq \varnothing $, let the following conditions hold.

(i) $b_{0}:V^{\alpha }\times V^{\alpha }\times V^{\theta -\beta }\rightarrow
\mathbb{R}$ is bounded, where $\alpha =\min \{\beta ,\theta -\theta _{2}\}$;

(ii) $b_{0}(v,w,Qw)=0$ for any $v,w\in \mathcal{V}$;

(iii) $u_{0}\in V^{\beta }$, $d_{0}\in D\left( A_{1}^{s}\right) $, and $g\in
L^{2}(0,T;V^{\beta -\theta })$.

Then we have
\begin{equation}
\left( u,d\right) \in L^{\infty }(0,T;V^{\beta }\times D\left(
A_{1}^{s}\right) )\cap L^{2}(0,T;V^{\beta +\theta }\times D(A_{1}^{\left(
2s+1\right) /2}))  \label{e:reg}
\end{equation}%
and%
\begin{align}
& \left\Vert u\left( t\right) \right\Vert _{\beta }^{2}+\left\Vert d\left(
t\right) \right\Vert _{2s}^{2}+\int_{0}^{t}\left( \left\Vert u\left(
s\right) \right\Vert _{\theta +\beta }^{2}+\left\Vert d\left( s\right)
\right\Vert _{2s+1}^{2}\right) ds  \label{claim-est} \\
& \leq \varphi \left( t\right) \left( \left\Vert u_{0}\right\Vert _{\beta
}^{2}+\left\Vert d_{0}\right\Vert _{2s}^{2}+\left\Vert g\right\Vert
_{L^{2}\left( 0,T;V^{\beta -\theta }\right) }^{2}\right) ,  \notag
\end{align}%
for some positive function $\varphi $ which depends on time, the norm of the
initial data $\left( u_{0},d_{0}\right) $ in $V^{\beta }\times D\left(
A_{1}^{s}\right) $ and on $g$.
\end{theorem}

\begin{proof}
The following estimates can be rigorously justified working with a
sufficiently smooth approximating solution, see Theorem \ref{t:exist}. We
will proceed formally. Pairing the first equation of (\ref{e:op}) with $%
\Lambda ^{2\beta }u$ yields%
\begin{align}
\frac{1}{2}& \frac{d}{dt}\left\langle u,\Lambda ^{2\beta }u\right\rangle
+\left\langle A_{0}u,\Lambda ^{2\beta }u\right\rangle +b_{0}\left(
u,u,\Lambda ^{2\beta }u\right)  \label{e:innpro} \\
& =\left\langle R_{0}\left( A_{1}d,d\right) ,\Lambda ^{2\beta
}u\right\rangle +\left\langle g,\Lambda ^{2\beta }u\right\rangle
+\left\langle \text{div}\left( \sigma _{Q}\right) ,\Lambda ^{2\beta
}u\right\rangle .  \notag
\end{align}%
Similarly, taking the inner product of the second equation of (\ref{e:op})
with $A_{1}^{2s}d$ gives%
\begin{align}
& \frac{1}{2}\frac{d}{dt}\left\Vert A_{1}^{s}d\right\Vert _{L^{2}}^{2}-\frac{%
1}{\lambda _{1}}\left\Vert A_{1}^{\left( 2s+1\right) /2}d\right\Vert
_{L^{2}}^{2}  \label{e:innpro2} \\
& =-\left\langle B_{1}\left( u,d\right) ,A_{1}^{2s}d\right\rangle
+\left\langle \omega _{Q}d,A_{1}^{2s}d\right\rangle +\frac{1}{\lambda _{1}}%
\left\langle f\left( d\right) ,A_{1}^{2s}d\right\rangle  \notag \\
& -\frac{\lambda _{2}}{\lambda _{1}}\left\langle
A_{Q}d,A_{1}^{2s}d\right\rangle .  \notag
\end{align}%
First, we are going to estimate the $b_{0}$-term as well as all the other
terms on the right-hand side of (\ref{e:innpro}), then we'll estimate the
terms on the right-hand side of (\ref{e:innpro2}). Combining the boundedness
of $b_{0}$ (see (i)) with the definition of $\Lambda ^{2\beta }u$\ and
Young's inequality, we find%
\begin{equation}
b_{0}(u,u,\Lambda ^{2\beta }u)\lesssim \delta ^{-1}\Vert u\Vert _{\theta
-\theta _{2}}^{2}\Vert u\Vert _{\beta }^{2}+\delta \left\Vert u\right\Vert
_{\beta +\theta }^{2},\qquad a.e.\ \text{ in }\left( 0,T\right) ,
\label{est10}
\end{equation}%
for any $\delta >0$; clearly, we also have%
\begin{equation}
\left\langle g,\Lambda ^{2\beta }u\right\rangle \lesssim \delta \left\Vert
u\right\Vert _{\beta +\theta }^{2}+C_{\delta }\left\Vert g\right\Vert
_{\beta -\theta }^{2}.  \label{est11}
\end{equation}%
Using a duality argument, we get%
\begin{align*}
\left\vert \left\langle R_{0}\left( A_{1}d,d\right) ,\Lambda ^{2\beta
}u\right\rangle \right\vert & \lesssim \left\Vert R_{0}\left(
A_{1}d,d\right) \right\Vert _{\beta -\theta }\left\Vert \Lambda ^{2\beta
}u\right\Vert _{\theta -\beta } \\
& \lesssim \left\Vert R_{0}\left( A_{1}d,d\right) \right\Vert
_{-1}\left\Vert u\right\Vert _{\theta +\beta }
\end{align*}%
since $\beta -\theta <-1$ (in all space dimensions, for $\beta \in J_{n}$).
Now, by Hahn-Banach theorem and H\"{o}lder's inequality,%
\begin{align*}
\left\Vert R_{0}\left( A_{1}d,d\right) \right\Vert _{-1}& =\sup_{\varphi
}\left\langle A_{1}d\cdot \nabla d,\varphi \right\rangle ,\text{ }\varphi
\in W^{1},\left\Vert \varphi \right\Vert _{1}=1, \\
& =\sup_{\varphi }\left\langle A_{1}^{s}d,A_{1}^{1-s}\left( \nabla d\varphi
\right) \right\rangle \\
& \leq \sup_{\varphi }\left\Vert A_{1}^{s}d\right\Vert _{L^{2}}\left\Vert
A_{1}^{1-s}\left( \nabla d\cdot \varphi \right) \right\Vert _{L^{2}} \\
& \leq \left\Vert A_{1}^{s}d\right\Vert _{L^{2}}\left\Vert d\right\Vert _{2}
\end{align*}%
so that the preceding inequality becomes%
\begin{equation}
\left\vert \left\langle R_{0}\left( A_{1}d,d\right) ,\Lambda ^{2\beta
}u\right\rangle \right\vert \lesssim \delta \left\Vert u\right\Vert _{\theta
+\beta }^{2}+C_{\delta }\left\Vert A_{1}^{s}d\right\Vert
_{L^{2}}^{2}\left\Vert A_{1}d\right\Vert _{L^{2}}^{2}.  \label{est11bis}
\end{equation}%
It remains to estimate the term involving $\sigma _{Q}$. To this end, we
note the identity%
\begin{align}
\left\langle \text{div}\left( \sigma _{Q}\right) ,\Lambda ^{2\beta
}u\right\rangle & =\mu _{1}\left\langle \text{div}\left( \left(
d^{T}A_{Q}d\right) \left( d\otimes d\right) \right) ,\Lambda ^{2\beta
}u\right\rangle  \label{3.10} \\
& +\left\langle \text{div}\left( \mu _{2}\mathcal{N}_{Q}\otimes d+\mu
_{3}d\otimes \mathcal{N}_{Q}\right) ,\Lambda ^{2\beta }u\right\rangle  \notag
\\
& +\left\langle \text{div}\left( \mu _{5}A_{Q}d\otimes d+\mu _{6}d\otimes
A_{Q}d\right) ,\Lambda ^{2\beta }u\right\rangle  \notag \\
& =:I_{1}+I_{2}+I_{3}.  \notag
\end{align}%
For these nonlinear terms, bounds are derived employing Lemma \ref{l:hole}
(Appendix), as follows:

\begin{itemize}
\item[(a)] The terms $\left( d^{T}A_{Q}d\right) _{ij}\left( d\otimes
d\right) _{jk}$ are a product of functions in $L^{2}$ and $H^{1}$ and
therefore bounded in $H^{\beta -\theta +1}$ since $\beta <\theta -n/2$.
Moreover, the terms $d_{i}d_{j}$ are a product of functions in $H^{2s}$ and $%
H^{1},$ and therefore bounded in $H^{1}$ since, by assumption, $2s>n/2$.

\item[(b)] To estimate all the nonlinear terms $I_{2},I_{3}$, we have to
estimate terms of the form $\left( A_{Q}d\right) _{i}d_{j}$, $d_{k}\left(
A_{Q}d\right) _{l}$, which are a product of functions in $L^{2}$ and $H^{2s}$
(respectively, in $H^{2s}$ and $L^{2}$), and therefore bounded in $H^{\beta
-\theta +1}$ provided that $\beta <2s+\theta -1-n/2$. On the other hand, we
have to estimate terms of the form $\left( \mathcal{N}_{Q}\right) _{i}d_{j}$%
, $d_{k}\left( \mathcal{N}_{Q}\right) _{l}$, which are a product of
functions in $L^{2}$ and $H^{2s}$ (respectively, in $H^{2s}$ and $L^{2}$)
and therefore are also bounded in $H^{\beta -\theta +1}$ for $\beta \in
J_{n}.$
\end{itemize}

By a duality argument, (a) and Young's inequality, it follows%
\begin{eqnarray}
\left\vert I_{1}\right\vert &\lesssim &\mu _{1}\left\Vert \left(
d^{T}A_{Q}d\right) \left( d\otimes d\right) \right\Vert _{\beta -\theta
+1}\left\Vert u\right\Vert _{\beta +\theta }  \label{est13} \\
&\lesssim &\mu _{1}\left\Vert d^{T}A_{Q}d\right\Vert _{L^{2}}\left\Vert
d\otimes d\right\Vert _{1}\left\Vert u\right\Vert _{\beta +\theta }  \notag
\\
&\lesssim &C_{\delta }\left( \mu _{1}\left\Vert d^{T}A_{Q}d\right\Vert
_{L^{2}}\right) ^{2}\left( \left\Vert d\right\Vert _{2s}\left\Vert
d\right\Vert _{1}\right) ^{2}+\delta \left\Vert u\right\Vert _{\beta +\theta
}^{2}.  \notag
\end{eqnarray}%
Using a duality argument once more and exploiting (b), we immediately get%
\begin{equation}
\left\vert I_{2}\right\vert +\left\vert I_{3}\right\vert \lesssim C_{\delta
}\left( \left\Vert \mathcal{N}_{Q}\right\Vert _{L^{2}}^{2}+\left\Vert
A_{Q}d\right\Vert _{L^{2}}^{2}\right) \left\Vert d\right\Vert
_{2s}^{2}+\delta \left\Vert u\right\Vert _{\beta +\theta }^{2}.
\label{estimate14}
\end{equation}%
Inserting (\ref{est10})-(\ref{est11bis}) and (\ref{est13})-(\ref{est14})
into (\ref{e:innpro}), then using the coercitivity of $A_{0}$, we derive%
\begin{align}
& \frac{1}{2}\frac{d}{dt}\left\langle u,\Lambda ^{2\beta }u\right\rangle
+c_{A_{0}}\left\Vert u\right\Vert _{\theta +\beta }^{2}  \label{innbis} \\
& \leq 5\delta \left\Vert u\right\Vert _{\beta +\theta }^{2}+C_{\delta
}\left( \Vert u\Vert _{\theta -\theta _{2}}^{2}\Vert u\Vert _{\beta
}^{2}+\left\Vert g\right\Vert _{\beta -\theta }^{2}\right) +C_{\delta
}\left\Vert A_{1}^{s}d\right\Vert _{L^{2}}^{2}\left\Vert A_{1}d\right\Vert
_{L^{2}}^{2}  \notag \\
& +C_{\delta }\left( \mu _{1}\left\Vert d^{T}A_{Q}d\right\Vert
_{L^{2}}^{2}\right) \left\Vert d\right\Vert _{1}^{2}\left\Vert d\right\Vert
_{2s}^{2}+C_{\delta }\left( \left\Vert \mathcal{N}_{Q}\right\Vert
_{L^{2}}^{2}+\left\Vert A_{Q}d\right\Vert _{L^{2}}^{2}\right) \left\Vert
d\right\Vert _{2s}^{2}.  \notag
\end{align}

We now turn to estimating the right-hand side of (\ref{e:innpro2}). First,
we note the identity%
\begin{align*}
& -\left\langle B_{1}\left( u,d\right) ,A_{1}^{2s}d\right\rangle
+\left\langle \omega _{Q}d,A_{1}^{2s}d\right\rangle +\frac{1}{\lambda _{1}}%
\left\langle f\left( d\right) ,A_{1}^{2s}d\right\rangle -\frac{\lambda _{2}}{%
\lambda _{1}}\left\langle A_{Q}d,A_{1}^{2s}d\right\rangle \\
& =-\left\langle A_{1}^{\left( 2s-1\right) /2}B_{1}\left( u,d\right)
,A_{1}^{\left( 2s+1\right) /2}d\right\rangle +\left\langle A_{1}^{\left(
2s-1\right) /2}\left( \omega _{Q}d\right) ,A_{1}^{\left( 2s+1\right)
/2}d\right\rangle \\
& +\frac{1}{\lambda _{1}}\left\langle A_{1}^{\left( 2s-1\right) /2}f\left(
d\right) ,A_{1}^{\left( 2s+1\right) /2}d\right\rangle -\frac{\lambda _{2}}{%
\lambda _{1}}\left\langle A_{1}^{\left( 2s-1\right) /2}\left( A_{Q}d\right)
,A_{1}^{\left( 2s+1\right) /2}d\right\rangle \\
& =J_{1}+...+J_{4}.
\end{align*}%
For these terms, bounds are derived employing Lemma \ref{l:hole} (Appendix),
as follows:

\begin{itemize}
\item[(c)] The terms $\left( Qu\right) _{i}\partial _{i}d_{j}$ are a product
of functions in $H^{\left( \beta +\theta +3\theta _{2}\right) /2}$ and $%
H^{s} $ and therefore bounded in $H^{2s-1}$ provided that $\beta \geq
4s-2-\theta -3\theta _{2}$ and $\beta >n+2s-2-\theta -3\theta _{2}$, which
are satisfied if $\beta \in J_{n}$.

\item[(d)] Finally, we have to estimate terms of the form $\left( \nabla
Qu\right) _{ij}d_{j}$, which are a product of functions in $H^{\left( \beta
+\theta +3\theta _{2}\right) /2-1}$ and $H^{s+1},$ and therefore bounded in $%
H^{2s-1}$ provided that $\beta \geq 4s-\theta -3\theta _{2}$ and $\beta
>n+2s-2-\theta -3\theta _{2}$, which once again holds for $\beta \in J_{n}$.
\end{itemize}

We begin with an easy bound on $J_{3}$ since $f\left( d\right) =\left(
\left\vert d\right\vert ^{2}-1\right) d$. We have%
\begin{eqnarray}
\left\vert J_{3}\right\vert &\lesssim &||A_{1}^{1/2}\left( f\left( d\right)
\right) ||_{L^{2}}||A_{1}^{\left( 2s+1\right) /2}d||_{L^{2}}  \label{est14}
\\
&\leq &\delta \left\Vert d\right\Vert _{2s+1}^{2}+C_{\delta }\left\Vert
A_{1}d\right\Vert _{L^{2}}^{2}\left\Vert A_{1}^{s}d\right\Vert _{L^{2}}^{2}.
\notag
\end{eqnarray}%
By an interpolation inequality in the triple $W^{2s+1}\subset W^{s+1}\subset
W^{1}$, Holder and Young inequalities, in view of (c) we find%
\begin{align}
\left\vert J_{1}\right\vert & \lesssim \left\Vert B_{1}\left( u,d\right)
\right\Vert _{2s-1}\left\Vert d\right\Vert _{2s+1}  \label{est15} \\
& \lesssim \left\Vert Qu\right\Vert _{\left( \beta +\theta +3\theta
_{2}\right) /2}\left\Vert d\right\Vert _{s+1}\left\Vert d\right\Vert _{2s+1}
\notag \\
& \lesssim \left\Vert u\right\Vert _{\left( \beta +\theta -\theta
_{2}\right) /2}\left\Vert d\right\Vert _{1}^{1/2}\left\Vert d\right\Vert
_{2s+1}^{3/2}  \notag \\
& \leq \delta \left\Vert d\right\Vert _{2s+1}^{2}+C_{\delta }\left\Vert
u\right\Vert _{\left( \beta +\theta -\theta _{2}\right) /2}^{4}\left\Vert
d\right\Vert _{1}^{2}  \notag \\
& \leq \delta \left\Vert d\right\Vert _{2s+1}^{2}+C_{\delta }\left\Vert
u\right\Vert _{\beta }^{2}\left\Vert u\right\Vert _{\theta -\theta
_{2}}^{2}\left\Vert d\right\Vert _{1}^{2}.  \notag
\end{align}%
Similar to the bound for $J_{1}$, using (d) one deduces%
\begin{align}
\left\vert J_{2}\right\vert & \lesssim \left\Vert \omega _{Q}d\right\Vert
_{2s-1}\left\Vert d\right\Vert _{2s+1}  \label{est15bis} \\
& \lesssim \left\Vert Qu\right\Vert _{\left( \beta +\theta +3\theta
_{2}\right) /2}\left\Vert d\right\Vert _{s+1}\left\Vert d\right\Vert _{2s+1}
\notag \\
& \leq \delta \left\Vert d\right\Vert _{2s+1}^{2}+C_{\delta }\left\Vert
u\right\Vert _{\beta }^{2}\left\Vert u\right\Vert _{\theta -\theta
_{2}}^{2}\left\Vert d\right\Vert _{1}^{2}.  \notag
\end{align}%
As a result of (\ref{est15})-(\ref{est15bis}), we also find the same bound (%
\ref{est15bis}) for $J_{4}$. Putting all the above estimates (\ref{est14})-(%
\ref{est15bis}) together with (\ref{e:innpro2}), we arrive at the inequality%
\begin{align}
& \frac{1}{2}\frac{d}{dt}\left\Vert A_{1}^{s}d\right\Vert _{L^{2}}^{2}-\frac{%
1}{\lambda _{1}}\left\Vert A_{1}^{\left( 2s+1\right) /2}d\right\Vert
_{L^{2}}^{2}  \label{inn2} \\
& \leq C_{\delta }\left\Vert u\right\Vert _{\beta }^{2}\left\Vert
u\right\Vert _{\theta -\theta _{2}}^{2}\left\Vert d\right\Vert
_{1}^{2}+4\delta \left\Vert d\right\Vert _{2s+1}^{2}+C_{\delta }\left\Vert
A_{1}d\right\Vert _{L^{2}}^{2}\left\Vert A_{1}^{s}d\right\Vert _{L^{2}}^{2}.
\notag
\end{align}

Finally, combining (\ref{innbis}) and (\ref{inn2}), we infer%
\begin{align}
& \frac{1}{2}\frac{d}{dt}\left[ \left\Vert A_{1}^{s}d\right\Vert
_{L^{2}}^{2}+\left\langle u,\Lambda ^{2\beta }u\right\rangle \right] -\frac{1%
}{\lambda _{1}}\left\Vert A_{1}^{\left( 2s+1\right) /2}d\right\Vert
_{L^{2}}^{2}+c_{A_{0}}\left\Vert u\right\Vert _{\theta +\beta }^{2}
\label{inn3} \\
& \leq 5\delta \left\Vert u\right\Vert _{\beta +\theta }^{2}+4\delta
\left\Vert d\right\Vert _{2s+1}^{2}+C_{\delta }\left\Vert u\right\Vert
_{\beta }^{2}\left\Vert u\right\Vert _{\theta -\theta _{2}}^{2}\left\Vert
d\right\Vert _{1}^{2}+C_{\delta }\left\Vert A_{1}d\right\Vert
_{L^{2}}^{2}\left\Vert A_{1}^{s}d\right\Vert _{L^{2}}^{2}  \notag \\
& +C_{\delta }\left( \Vert u\Vert _{\theta -\theta _{2}}^{2}\Vert u\Vert
_{\beta }^{2}+\left\Vert g\right\Vert _{\beta -\theta }^{2}\right)
+C_{\delta }\left( \left\Vert \mathcal{N}_{Q}\right\Vert
_{L^{2}}^{2}+\left\Vert A_{Q}d\right\Vert _{L^{2}}^{2}\right) \left\Vert
d\right\Vert _{2s}^{2}  \notag \\
& +C_{\delta }\left( \mu _{1}\left\Vert d^{T}A_{Q}d\right\Vert
_{L^{2}}\right) ^{2}\left\Vert d\right\Vert _{1}^{2}\left\Vert d\right\Vert
_{2s}^{2}.  \notag
\end{align}%
Thus, choosing a sufficiently small $\delta \sim \min \left(
c_{A_{0}},-\lambda _{1}^{-1}\right) >0$ in (\ref{inn3}), by Gronwall's
inequality we conclude (\ref{claim-est}). The proof of the theorem is
finished.
\end{proof}

To clarify the previous result in the case of specific models, the
corresponding conditions (in particular, (i)-(ii)) when $b_{0}$ is either $%
b_{00}$ or $b_{01}$, as given by Example \ref{y:spaces}, are listed below in
the following remarks. We note that this procedure always produces a new
interval $\mathcal{Y}_{n}$ for the parameter $\beta $ so that one must
ensure that $J_{n}\cap \mathcal{Y}_{n}$ stays nonempty.

\begin{remark}
\label{reg-rem}Let $4\theta +4\theta _{1}+2\theta _{2}>n+2$, $2\theta
+2\theta _{1}\geq 1-k$, $\theta +2\theta _{2}\geq 1$, $3\theta +4\theta
_{1}\geq 1$, $\theta +2\theta _{1}\geq \ell $, and $3\theta +2\theta
_{1}+2\theta _{2}\geq 2-\ell $, for some $k,\ell \in \{0,1\}$. For%
\begin{equation*}
\beta \in (\frac{n+2}{2}-2(\theta _{1}+\theta _{2})-\theta ,3\theta +2\theta
_{1}-\frac{n+2}{2})\cap \lbrack \frac{1-\ell }{2}-\theta _{1}-\theta
_{2},\min \{2\theta +\theta _{2}-1,2\theta -\theta _{2}+2\theta _{1}-k\}].
\end{equation*}%
from \cite[Proposition 2.5]{HLT} we infer that the trilinear form $b_{00}$
satisfies the hypotheses (i)-(ii) of the above theorem.
\end{remark}

\begin{remark}
Let $4\theta +4\theta _{1}+2\theta _{2}>n+2$, $\theta +2\theta _{2}\geq 0$,
and $\theta +2\theta _{1}\geq 1$. For%
\begin{equation*}
\beta \in (\frac{n+2}{2}-2(\theta _{1}+\theta _{2})-\theta ,3\theta +2\theta
_{1}-\frac{n+2}{2})\cap \lbrack \frac{1}{2}-\theta _{1}-\theta _{2},\min
\{2\theta +\theta _{2},2\theta -\theta _{2}+2\theta _{1}-1\}].
\end{equation*}%
from \cite[Proposition 2.5]{HLT} it follows that the trilinear form $b_{01}$
satisfies the hypotheses (i)-(ii) of the above theorem.
\end{remark}

\begin{remark}
\label{reg-gen}We note that the interval $J_{n}$ is exactly the same for any
fixed values of $\theta ,\theta _{2}$. This is the case, for instance, when $%
\theta =\theta _{2}=1$, refer to Table \ref{t:spec2}. We also observe that
for $J_{n}\neq \varnothing $, we must always ask that $\theta +\theta
_{2}>n/2$.
\end{remark}

\begin{remark}
For any $s\in \left( 0.75,0.875\right) $ such that $\beta \in \left[
4s-4,-0.5\right) $, Theorem \ref{t:reg} implies global regularity of the
energy solutions of Definition \ref{weak} for the modified Leray-EL-$\alpha $
(ML-EL-$\alpha $) model, the SBM-EL model and the NS-EL-$\alpha $ system in
three space dimensions. We emphasize that $J_{3}=\varnothing $ when $s=1$
for all these models, and that these results are valid without any
restrictions on the physical parameters $\lambda _{1},\lambda _{2},$ $\mu
_{1},...,\mu _{6},$ other than what was already assumed in Section \ref%
{s:prelim} (cf. (\ref{lama1a})-(\ref{lam2}) and \textbf{Cases 1-2}). On the
other hand, global regularity of the energy solutions for the 3D Leray-EL-$%
\alpha $ ($\theta =1,\theta _{2}=0$) model and the 3D NSV-EL model ($\theta
=0,$ $\theta _{2}=1$) are not covered here since both models fail to satisfy
the condition $\theta +\theta _{2}>n/2$.
\end{remark}

\begin{remark}
\label{rem-bd}Any regular weak solution $\left( u,d\right) $, as given by
Theorem \ref{t:reg}, satisfies
\begin{equation*}
d\in L^{\infty }\left( 0,T;L^{\infty }\left( \Omega \right) \right)
\end{equation*}%
due to the Sobolev embedding $W^{2s}\subset L^{\infty }$, as $2s>n/2$.
\end{remark}

\subsection{Uniqueness and stability}

\label{ss:stab}

Now we shall provide sufficient conditions for uniqueness and continuous
dependence on the initial data for any weak solutions of the general
three-parameter family of regularized models. Recall that $\theta _{1}\in
\mathbb{R}$ and $\theta ,\theta _{2}\geq 0.$

Our first result is concerned with the case when a maximum principle applies
to the director field $d$ (i.e., when $\lambda _{2}=0,$ such that any
stretching of the crystal molecules is ignored).

\begin{theorem}
\label{t:stab} Let $\left( u_{i},d_{i}\right) \in L^{\infty }(0,T;\mathcal{Y}%
_{\theta _{2}}),$ $i=1,2,$ be two energy solutions in the sense of
Definition \ref{weak}, corresponding to the initial conditions $\left(
u_{i}(0),d_{i}\left( 0\right) \right) \in \mathcal{Y}_{\theta _{2}}$, $%
i=1,2. $ Assume Proposition \ref{maxp} and the following conditions.

(i) $b_{0}:V^{\sigma _{1}}\times V^{\theta -\theta _{2}}\times V^{\sigma
_{2}}\rightarrow \mathbb{R}$ is bounded for some $\sigma _{1}\leq \theta
-\theta _{2}$ and $\sigma _{2}\leq \theta +\theta _{2}$ with $\sigma
_{1}+\sigma _{2}\leq \theta $.

(ii) $b_{0}(v,w,Qw)=0$ for any $v\in V^{\sigma _{1}}$ and $w\in V^{\sigma
_{2}}.$

Further suppose that $\theta _{2}\geq 1$. Then the following estimate holds%
\begin{align}
& \Vert u_{1}(t)-u_{2}(t)\Vert _{-\theta _{2}}^{2}+\left\Vert d_{1}\left(
t\right) -d_{2}\left( t\right) \right\Vert _{1}^{2}  \label{uniq_stab} \\
& +\int_{0}^{t}\left( \Vert u_{1}\left( s\right) -u_{2}\left( s\right) \Vert
_{\theta -\theta _{2}}^{2}+\left\Vert A_{1}\left( d_{1}\left( s\right)
-d_{2}\left( s\right) \right) \right\Vert _{L^{2}}^{2}\right) ds  \notag \\
& \leq \varrho \left( t\right) \left( \Vert u_{1}(0)-u_{2}(0)\Vert _{-\theta
_{2}}^{2}+\left\Vert d_{1}\left( 0\right) -d_{2}\left( 0\right) \right\Vert
_{1}^{2}\right) ,  \notag
\end{align}%
for $t\in \lbrack 0,T],$ for some positive continuous function $\varrho :%
\mathbb{R}_{+}\rightarrow \mathbb{R}_{+},$ $\varrho \left( 0\right) >0$,
which depends only on the initial data $\left( u_{i}\left( 0\right)
,d_{i}\left( 0\right) \right) $ in $\mathcal{Y}_{\theta _{2}}$-norm.
\end{theorem}

\begin{proof}
First of all, when $\lambda _{2}=0$ by Proposition \ref{maxp} there exists a
constant $M>0$ such that $\left\Vert d_{i}\right\Vert _{L_{t}^{\infty
}L_{x}^{\infty }}\leq M,$ $i=1,2.$ Set $\nabla _{d}W\left( d\right) =f\left(
d\right) $ and let $u=u_{1}-u_{2}$, $d=d_{1}-d_{2}$. Then subtracting the
equations for $\left( u_{1},d_{1}\right) $ and $\left( u_{2},d_{2}\right) $
we have%
\begin{align}
& \left\langle \partial _{t}u,w\right\rangle +\left\langle
A_{0}u,w\right\rangle +\left\langle B_{0}\left( u,u_{1}\right)
,w\right\rangle +\left\langle B_{0}\left( u_{2},u\right) ,w\right\rangle
\label{diffuniq1} \\
& =\left\langle R_{0}\left( A_{1}d_{2},d\right) ,w\right\rangle
+\left\langle R_{0}\left( A_{1}d,d_{1}\right) ,w\right\rangle +\left\langle
\sigma _{Qu_{1}}-\sigma _{Qu_{2}},\nabla w\right\rangle ,  \notag
\end{align}%
and%
\begin{align}
& \left\langle \partial _{t}d,\eta \right\rangle -\frac{1}{\lambda _{1}}%
\left\langle A_{1}d,\eta \right\rangle +\left\langle B_{1}\left(
u,d_{1}\right) ,\eta \right\rangle +\left\langle B_{1}\left( u_{2},d\right)
,\eta \right\rangle  \label{diffuniq2} \\
& =\frac{1}{\lambda _{1}}\left\langle f\left( d_{1}\right) -f\left(
d_{2}\right) ,\eta \right\rangle +\left\langle \omega _{Qu}d_{1}+\omega
_{Qu_{2}}d,\eta \right\rangle .  \notag
\end{align}%
Here, we denote $\sigma _{Qu_{i}},$ $A_{Qu_{i}}$, $\omega _{Qu_{i}}$ to be
exactly the same stress tensors from (\ref{v5}), (\ref{v4}) associated with
a given weak solution $u=u_{i}$, $i=1,2$, while $\omega _{Qu}:=\omega
_{Qu_{1}}-\omega _{Qu_{2}}$. First, observe that by the assumptions on $%
\theta ,\theta _{2},$ according to the estimates that we will perform below,
the weak solution $\left( u_{i},d_{i}\right) $\ of (\ref{e:op}) enjoys in
fact \emph{additional} regularity. Indeed, these subsequent estimates yield
that $R_{0}\left( A_{1}d_{i},d_{i}\right) $, $B_{0}\left( u_{i},u_{i}\right)
\in L^{2}\left( 0,T;V^{-\theta -\theta _{2}}\right) $ and $\sigma
_{Qu_{i}}\in L^{2}\left( 0,T;V^{-\theta -\theta _{2}+1}\right) $, and $%
B_{1}\left( u_{i},d_{i}\right) $, $\omega _{Qu_{i}}d_{i}\in L^{2}\left(
0,T;L^{2}\left( \Omega \right) \right) $. This regularity effectively
translates to regularity of the time derivatives $\partial _{t}u\in
L^{2}\left( 0,T;V^{-\theta -\theta _{2}}\right) $ and $\partial _{t}d\in
L^{2}\left( 0,T;L^{2}\left( \Omega \right) \right) $ such that each of the
corresponding functional pairings $\left\langle \partial
_{t}u,w\right\rangle $ and $\left\langle \partial _{t}d,\eta \right\rangle $
is integrable for $w\in L^{2}\left( 0,T;V^{\theta +\theta _{2}}\right) $ and
$\eta \in L^{2}\left( 0,T;L^{2}\left( \Omega \right) \right) $,
respectively. Thus, in what follows we can take $w=Qu$ and $\eta =A_{1}d$
into (\ref{diffuniq1})-(\ref{diffuniq2}) to infer%
\begin{align}
& \frac{d}{dt}\left( \Vert u\Vert _{-\theta _{2}}^{2}+\left\Vert
A_{1}^{1/2}d\right\Vert _{L^{2}}^{2}\right) +2c_{A_{0}}\Vert u\Vert _{\theta
-\theta _{2}}^{2}-\frac{2}{\lambda _{1}}\left\Vert A_{1}d\right\Vert
_{L^{2}}^{2}  \label{uniqest1} \\
& \leq 2\left\vert b_{0}(u,u_{1},Qu)\right\vert +2\left\vert b_{1}\left(
u,d,A_{1}d_{2}\right) \right\vert +2\left\vert b_{1}\left(
u_{2},d,A_{1}d\right) \right\vert +\frac{2}{\lambda _{1}}\left\vert
\left\langle f\left( d_{1}\right) -f\left( d_{2}\right) ,A_{1}d\right\rangle
\right\vert  \notag \\
& +2\left\vert \left\langle \omega _{Qu}d_{1}+\omega
_{Qu_{2}}d,A_{1}d\right\rangle \right\vert +2\left\vert \left\langle \sigma
_{Qu_{1}}-\sigma _{Qu_{2}},\nabla \left( Qu\right) \right\rangle \right\vert
\notag \\
& =:I_{1}+I_{2}+...+I_{6}.  \notag
\end{align}%
All the terms $I_{1}-I_{4}$ on the right-hand side of (\ref{uniqest1}) were
estimated in \cite[Theorem 3.4]{GM-JNS} for the corresponding regularized
simplified Ericksen-Leslie system (\ref{e:op}) when $\sigma _{Q}\equiv 0$
and $\omega _{Q}\equiv 0$. The bounds\footnote{%
The estimates from (\ref{uniqest2}) performed in \cite{GM-JNS} required that
$\theta +\theta _{2}\geq 1,$ $\theta _{2}\geq 0$ in 2D and $\theta _{2}\geq
1 $ in 3D. Alternatively, one can replace these conditions by $\theta
+\theta _{2}>\frac{n}{2}$ which is complementary.} for these nonlinear terms
read as follows:%
\begin{equation}
\left\{
\begin{array}{l}
I_{1}\lesssim C_{\delta }\Vert u\Vert _{-\theta _{2}}^{2}\Vert u_{1}\Vert
_{\theta -\theta _{2}}^{2}+\delta \Vert u\Vert _{\theta -\theta _{2}}^{2},
\\
I_{2}\lesssim \delta \left\Vert A_{1}d\right\Vert _{L^{2}}^{2}+C_{\delta
}\left( \left\Vert d\right\Vert _{1}^{2}+\left\Vert u\right\Vert _{-\theta
_{2}}^{2}\left\Vert A_{1}d_{2}\right\Vert _{L^{2}}^{2}\right) , \\
I_{3}\lesssim \delta \left\Vert A_{1}d\right\Vert _{L^{2}}^{2}+C_{\delta
}\left\Vert u_{2}\right\Vert _{\theta -\theta _{2}}^{2}\left\Vert
u_{2}\right\Vert _{-\theta _{2}}^{\kappa }\left\Vert d\right\Vert _{1}^{2},
\\
I_{4}\lesssim \delta \left\Vert A_{1}d\right\Vert _{L^{2}}^{2}+C_{\delta
}\left\Vert d\right\Vert _{L^{2}}^{2},%
\end{array}%
\right.  \label{uniqest2}
\end{equation}%
for any $\delta >0$, for some $\kappa =\kappa \left( n\right) \geq 2$ and $%
C_{\delta }>0$ sufficiently large. Now we proceed to estimate $I_{5}$ and $%
I_{6}$. This is mainly the place where the main condition $\theta _{2}\geq 1$
(in all space dimensions) must be enforced. We begin with the most
challenging term $I_{6}$. With the following definitions%
\begin{align*}
\widehat{\mathcal{N}}& :=\frac{1}{\lambda _{1}}\left( A_{1}d+\left( f\left(
d_{1}\right) -f\left( d_{2}\right) \right) d_{1}+f\left( d_{2}\right)
d\right) , \\
\mathcal{N}_{i}& :=\dot{d}_{i}-\omega _{Qu_{i}}d_{i},
\end{align*}%
we note the following identity%
\begin{align}
\left\langle \sigma _{Qu_{1}}-\sigma _{Qu_{2}},\nabla \left( Qu\right)
\right\rangle & =\mu _{1}\left\langle \left( d_{1}^{T}A_{Qu_{1}}d_{1}\right)
\left( d_{1}\otimes d_{1}\right) -\left( d_{2}^{T}A_{Qu_{2}}d_{2}\right)
\left( d_{2}\otimes d_{2}\right) ,\nabla Qu\right\rangle  \label{diff-Leslie}
\\
& +\mu _{2}\left\langle \mathcal{N}_{1}\otimes d_{1}-\mathcal{N}_{2}\otimes
d_{2},\nabla Qu\right\rangle +\mu _{3}\left\langle d_{1}\otimes \mathcal{N}%
_{1}-d_{2}\otimes \mathcal{N}_{2},\nabla Qu\right\rangle  \notag \\
& +\mu _{5}\left\langle \left( A_{Qu_{1}}d_{1}\right) \otimes d_{1}-\left(
A_{Qu_{2}}d_{2}\right) \otimes d_{2},\nabla Qu\right\rangle  \notag \\
& +\mu _{6}\left\langle d_{1}\otimes \left( A_{Qu_{1}}d_{1}\right)
-d_{2}\otimes \left( A_{Qu_{2}}d_{2}\right) ,\nabla Qu\right\rangle  \notag
\\
& =\mu _{1}\left\langle \left(
d^{T}A_{Qu_{1}}d_{1}+d_{2}^{T}A_{Qu}d_{1}+d_{2}^{T}A_{Qu_{2}}d\right) \left(
d_{1}\otimes d_{1}\right) ,\nabla Qu\right\rangle  \notag \\
& +\mu _{1}\left\langle \left( d_{2}^{T}A_{Qu_{2}}d_{2}\right) \left(
d\otimes d_{1}+d_{2}\otimes d\right) ,\nabla Qu\right\rangle  \notag \\
& +\mu _{2}\left\langle \widehat{\mathcal{N}}\otimes d_{1}+\mathcal{N}%
_{2}\otimes d,\nabla Qu\right\rangle +\mu _{3}\left\langle d_{1}\otimes
\widehat{\mathcal{N}}+d\otimes \mathcal{N}_{2},\nabla Qu\right\rangle  \notag
\\
& +\mu _{5}\left\langle \left( A_{Qu}d_{1}+A_{Qu_{2}}d\right) \otimes
d_{1}+\left( A_{Qu_{2}}d_{2}\right) \otimes d,\nabla Qu\right\rangle  \notag
\\
& +\mu _{6}\left\langle d_{1}\otimes \left( A_{Qu}d_{1}+A_{Qu_{2}}d\right)
+d\otimes \left( A_{Qu_{2}}d_{2}\right) ,\nabla Qu\right\rangle  \notag \\
& =:I_{61}+I_{62}+...+I_{66}.  \notag
\end{align}%
We shall estimate $I_{61}-I_{66}$ now. By the Holder inequality and proper
Sobolev embedding theorems (e.g., $V^{-\theta _{2}}\subseteq V^{1-2\theta
_{2}}$ and $W^{2}\subset L^{\infty }$), we have%
\begin{align}
\left\vert I_{61}\right\vert & \lesssim \left( \left\Vert d\right\Vert
_{L^{\infty }}\left\Vert \nabla Qu_{1}\right\Vert _{0}\left\Vert
d_{1}\right\Vert _{L^{\infty }}+\left\Vert d_{2}\right\Vert _{L^{\infty
}}\left\Vert \nabla Qu\right\Vert _{0}\left\Vert d_{1}\right\Vert
_{L^{\infty }}+\left\Vert d_{2}\right\Vert _{L^{\infty }}\left\Vert \nabla
Qu_{2}\right\Vert _{0}\left\Vert d\right\Vert _{L^{\infty }}\right)
\label{uniqest61} \\
& \times \left\Vert d_{1}\right\Vert _{L^{\infty }}^{2}\left\Vert \nabla
Qu\right\Vert _{0}  \notag \\
& \leq C_{M}\left( \left\Vert d\right\Vert _{L^{\infty }}\left\Vert
u_{1}\right\Vert _{-\theta _{2}}\left\Vert u\right\Vert _{-\theta
_{2}}+\left\Vert u\right\Vert _{-\theta _{2}}^{2}+\left\Vert
u_{2}\right\Vert _{-\theta _{2}}\left\Vert d\right\Vert _{L^{\infty
}}\left\Vert u\right\Vert _{-\theta _{2}}\right)  \notag \\
& \leq C_{M}\left\Vert u\right\Vert _{-\theta _{2}}^{2}+\delta \left\Vert
A_{1}d\right\Vert _{L^{2}}^{2}+C_{M,\delta }\left( \left\Vert
u_{1}\right\Vert _{-\theta _{2}}^{2}+\left\Vert u_{2}\right\Vert _{-\theta
_{2}}^{2}\right) \left\Vert u\right\Vert _{-\theta _{2}}^{2}  \notag
\end{align}%
for any $\delta >0,$ since $\left\Vert d_{i}\right\Vert _{L_{t}^{\infty
}L_{x}^{\infty }}\leq M,$ $i=1,2$. Similarly, we have%
\begin{align}
\left\vert I_{62}\right\vert & \lesssim \left\Vert d\right\Vert _{L^{\infty
}}\left( \left\Vert d_{2}\right\Vert _{L^{\infty }}^{2}\left\Vert
d_{1}\right\Vert _{L^{\infty }}+\left\Vert d_{2}\right\Vert _{L^{\infty
}}^{3}\right) \left\Vert \nabla Qu_{2}\right\Vert _{0}\left\Vert \nabla
Qu\right\Vert _{0}  \label{uniqest62} \\
& \leq C_{M,\delta }\left\Vert u_{2}\right\Vert _{-\theta
_{2}}^{2}\left\Vert u\right\Vert _{-\theta _{2}}^{2}+\delta \left\Vert
A_{1}d\right\Vert _{L^{2}}^{2}.  \notag
\end{align}%
Moreover, recalling that $\mathcal{N}_{i}\in L^{2}\left( 0,T;L^{2}\left(
\Omega \right) \right) ,$ $i=1,2$, for every weak energy solution $\left(
u_{i},d_{i}\right) $ of Definition \ref{weak}, we estimate%
\begin{align}
\left\vert I_{63}\right\vert & \leq C_{M}\left( \left\Vert A_{1}d\right\Vert
_{L^{2}}+\left\Vert d\right\Vert _{L^{2}}\right) \left\Vert \nabla
Qu\right\Vert _{0}+\left\Vert \mathcal{N}_{2}\right\Vert _{L^{2}}\left\Vert
d\right\Vert _{L^{\infty }}\left\Vert \nabla Qu\right\Vert _{0}
\label{uniqest63} \\
& \leq \delta \left\Vert A_{1}d\right\Vert _{L^{2}}^{2}+C_{M,\delta
}\left\Vert u\right\Vert _{-\theta _{2}}^{2}\left( 1+\left\Vert \mathcal{N}%
_{2}\right\Vert _{L^{2}}^{2}\right) .  \notag
\end{align}%
The bound for $I_{64}$ is exactly the same as in (\ref{uniqest63}). On the
other hand, the bound for the last integrals $I_{65},I_{66}$ can be obtained
as follows:%
\begin{align}
\max \left\{ \left\vert I_{65}\right\vert ,\left\vert I_{66}\right\vert
\right\} & \leq C_{M}\left( \left\Vert \nabla Qu\right\Vert
_{0}^{2}+\left\Vert \nabla Qu_{2}\right\Vert _{0}\left\Vert d\right\Vert
_{L^{\infty }}\left\Vert \nabla Qu\right\Vert _{0}\right)  \label{uniqest66}
\\
& \leq \delta \left\Vert A_{1}d\right\Vert _{L^{2}}^{2}+C_{M,\delta
}\left\Vert u\right\Vert _{-\theta _{2}}^{2}\left( 1+\left\Vert
u_{2}\right\Vert _{-\theta _{2}}^{2}\right) .  \notag
\end{align}%
To estimate $I_{5}$, we start with the preliminary estimate%
\begin{align*}
\left\vert \left\langle \omega _{Qu}d_{1},A_{1}d\right\rangle \right\vert &
\lesssim \left\Vert d_{1}\right\Vert _{L^{\infty }}\left\Vert
A_{1}d\right\Vert _{L^{2}}\left\Vert \nabla Qu\right\Vert _{0} \\
& \leq \delta \left\Vert A_{1}d\right\Vert _{L^{2}}^{2}+C_{M,\delta
}\left\Vert u\right\Vert _{-\theta _{2}}^{2}.
\end{align*}%
Next, by Agmon inequalities we have%
\begin{align*}
\left\vert \left\langle \omega _{Qu_{2}}d,A_{1}d\right\rangle \right\vert &
\lesssim \left\Vert d\right\Vert _{L^{\infty }}\left\Vert \nabla
Qu_{2}\right\Vert _{0}\left\Vert A_{1}d\right\Vert _{L^{2}} \\
& \lesssim \left\{
\begin{array}{ll}
\left\Vert d\right\Vert _{1}^{1/2}\left\Vert A_{1}d\right\Vert
_{L^{2}}^{3/2}\left\Vert u_2\right\Vert _{-\theta _{2}}, & \text{if }n=2, \\
\left\Vert d\right\Vert _{1}^{1/4}\left\Vert A_{1}d\right\Vert
_{L^{2}}^{7/4}\left\Vert u_2\right\Vert _{-\theta _{2}}, & \text{if }n=3.%
\end{array}%
\right.
\end{align*}%
Therefore, using Young's inequality, it follows%
\begin{equation}
\left\vert I_{5}\right\vert \leq \delta \left\Vert A_{1}d\right\Vert
_{L^{2}}^{2}+C_{\delta }\left[ \left\Vert u\right\Vert _{-\theta
_{2}}^{2}+\left( \left\Vert u_{2}\right\Vert _{-\theta _{2}}^{4}+\left\Vert
u_{2}\right\Vert _{-\theta _{2}}^{8}\right) \left\Vert d\right\Vert _{1}^{2}%
\right] .  \label{uniqest5}
\end{equation}%
Combining (\ref{uniqest2}), (\ref{uniqest61})-(\ref{uniqest5}), then
choosing a sufficiently small $\delta \sim \min \left( c_{A_{0}},-\lambda
_{1}^{-1}\right) >0$ into (\ref{uniqest1}), by application of Gronwall's
inequality, one finds%
\begin{equation*}
\Vert u(t)\Vert _{-\theta _{2}}^{2}+\left\Vert d\left( t\right) \right\Vert
_{1}^{2}\leq \left( \Vert u(0)\Vert _{-\theta _{2}}^{2}+\left\Vert d\left(
0\right) \right\Vert _{1}^{2}\right) \exp \int_{0}^{t}\Theta \left( s\right)
ds,
\end{equation*}%
for a suitable function $\Theta \in L^{1}\left( 0,T\right) .$ Integrating (%
\ref{uniqest1}) once more over $\left( 0,t\right) $ yields the desired
inequality (\ref{uniq_stab}). The proof is finished.
\end{proof}

To clarify these stability results at least in the case of the specific
models listed in Table \ref{t:spec2}, the corresponding conditions and
stability results derived from Theorem \ref{t:stab} are given below.

\begin{remark}
Exploiting \cite[Proposition 2.5]{HLT}, the trilinear form $b_{00}$
satisfies the hypotheses of Theorem \ref{t:stab} provided $\theta +\theta
_{1}\geq \frac{1-k}{2}$, $\theta +2\theta _{1}\geq k$, $\theta +\theta
_{2}\geq \frac{1}{2}$, $2\theta +2\theta _{1}+\theta _{2}>\frac{n+2}{2}$,
and $3\theta +2\theta _{1}+2\theta _{2}\geq 2-k$, for some $k\in \{0,1\}$.
The trilinear form $b_{01}$ satisfies the hypotheses of Theorem \ref{t:stab}
for $\theta +2\theta _{1}\geq 1$, $\theta +\theta _{1}\geq \frac{1}{2},$ $%
\theta +\theta _{2}\geq 0$, $2\theta +2\theta _{1}+\theta _{2}>\frac{n+2}{2}$%
, and $3\theta +2\theta _{1}+2\theta _{2}\geq 1$. Together with $\theta
_{2}\geq 1$ and $\lambda _{2}=0$, these assumptions allow us to recover the
stability and uniqueness of energy solutions for the 3D NS--EL-$\alpha $%
-model, the 3D NSV--EL-$\alpha $-model, the 3D ML--EL-$\alpha $-model and
the 3D SBM--EL model (see Table \ref{t:spec2}). These results were not
reported anywhere else.
\end{remark}

We conclude the section with a result that handles the general case when $%
\lambda _{2}\neq 0$.

\begin{theorem}
\label{t:stab2} Let%
\begin{equation*}
\left( u_{i},d_{i}\right) \in L^{\infty }(0,T;V^{\beta }\times D\left(
A_{1}^{s}\right) )\cap L^{2}(0,T;V^{\beta +\theta }\times D(A_{1}^{\left(
2s+1\right) /2}))
\end{equation*}%
be two energy solutions that satisfy the assumptions of Theorem \ref{t:reg}.
For $\theta _{2}\geq 1$, the estimate%
\begin{align}
& \Vert u_{1}(t)-u_{2}(t)\Vert _{-\theta _{2}}^{2}+\left\Vert d_{1}\left(
t\right) -d_{2}\left( t\right) \right\Vert _{1}^{2}  \label{c-dep} \\
& +\int_{0}^{t}\left( \Vert u_{1}\left( s\right) -u_{2}\left( s\right) \Vert
_{\theta -\theta _{2}}^{2}+\left\Vert A_{1}\left( d_{1}\left( s\right)
-d_{2}\left( s\right) \right) \right\Vert _{L^{2}}^{2}\right) ds  \notag \\
& \leq \varrho \left( t\right) \left( \Vert u_{1}(0)-u_{2}(0)\Vert _{-\theta
_{2}}^{2}+\left\Vert d_{1}\left( 0\right) -d_{2}\left( 0\right) \right\Vert
_{1}^{2}\right) ,  \notag
\end{align}%
holds for $t\in \lbrack 0,T],$ for some positive continuous function $%
\varrho :\mathbb{R}_{+}\rightarrow \mathbb{R}_{+},$ $\varrho \left( 0\right)
>0$, which depends on the initial data $\left( u_{i}\left( 0\right)
,d_{i}\left( 0\right) \right) $ in $V^{\beta }\times W^{2s}$-norm.
\end{theorem}

\begin{proof}
Indeed, the main ingredient in the proof of Theorem \ref{t:stab} was the
fact that $d_{i}\in L_{t}^{\infty }\left( L_{x}^{\infty }\right) $ which is
now provided by Remark \ref{rem-bd}. It is worth pointing out that in the
general case when $\lambda _{2}\neq 0$, the inequality (\ref{uniqest1})
reads as follows:%
\begin{align}
& \frac{d}{dt}\left( \Vert u\Vert _{-\theta _{2}}^{2}+\left\Vert
A_{1}^{1/2}d\right\Vert _{L^{2}}^{2}\right) +2c_{A_{0}}\Vert u\Vert _{\theta
-\theta _{2}}^{2}-\frac{2}{\lambda _{1}}\left\Vert A_{1}d\right\Vert
_{L^{2}}^{2}  \label{uniqest1bis} \\
& \leq 2\left\vert b_{0}(u,u_{1},Qu)\right\vert +2\left\vert b_{1}\left(
u,d,A_{1}d_{2}\right) \right\vert +2\left\vert b_{1}\left(
u_{2},d,A_{1}d\right) \right\vert +\frac{2}{\lambda _{1}}\left\vert
\left\langle f\left( d_{1}\right) -f\left( d_{2}\right) ,A_{1}d\right\rangle
\right\vert  \notag \\
& +2\left\vert \left\langle \omega _{Qu}d_{1}+\omega
_{Qu_{2}}d,A_{1}d\right\rangle \right\vert +2\left\vert \left\langle \sigma
_{Qu_{1}}-\sigma _{Qu_{2}},\nabla \left( Qu\right) \right\rangle \right\vert
-2\lambda _{2}\lambda _{1}^{-1}\left\vert \left\langle
A_{Qu}d_{1}+A_{Qu_{2}}d,A_{1}d\right\rangle \right\vert  \notag \\
& =:I_{1}+I_{2}+...+I_{6}+I_{7}.  \notag
\end{align}%
More precisely, with respect to (\ref{uniqest1}), there is one additional
term $I_{7}$ on the right-hand side. Bounds on the first six terms $I_{1}$-$%
I_{6}$ are already provided by (\ref{uniqest2})-(\ref{uniqest5}). To find a
proper bound for the final term $I_{7}$ one may proceed \emph{verbatim} as
in getting estimate (\ref{uniqest5}) for the term $I_{5}$; indeed, note that
$A_{Q}$ and $\omega _{Q}$ are in fact the symmetric and the skew-symmetric
parts of the strain rate, respectively (cf. (\ref{v4})). Hence, the proof of
(\ref{c-dep}) follows from that of (\ref{uniq_stab}) with some minor
modifications.
\end{proof}

\section{Finite dimensional global attractors}

\label{ss:global}

In this section we establish the existence of (smooth)\ finite dimensional
global attractors for the general three-parameter family of regularized
models (\ref{e:op}). For the sake of reference below, recall the following
definition for the space of translation bounded functions%
\begin{equation*}
L_{tb}^{2}\left( \mathbb{R}_{+};X\right) :=\left\{ g\in L_{loc}^{2}\left(
\mathbb{R}_{+};X\right) :\left\Vert g\right\Vert _{L_{tb}^{2}\left( \mathbb{R%
}_{+};X\right) }^{2}:=\sup_{t\geq 0}\int_{t}^{t+1}\left\Vert g\left(
s\right) \right\Vert _{X}^{2}ds<\infty \right\} ,
\end{equation*}%
where $X$ is a given Banach space.

We begin with a first basic dissipative inequality which is satisfied by any
weak energy solution of problem (\ref{e:op}). The following result holds for
any $\theta ,\theta _{2}\geq 0$.

\begin{proposition}
\label{t:attr-exist}Let $\left( u,d\right) \in L_{\mathrm{loc}}^{\infty
}(0,\infty ;\mathcal{Y}_{\theta _{2}})\cap L_{\mathrm{loc}}^{2}(0,\infty
;V^{\theta -\theta _{2}}\times D\left( A_{1}\right) )$ be any energy
solution in the sense of Definition \ref{weak-energy} with $\left( u\left(
0\right) ,d\left( 0\right) \right) \in \mathcal{Y}_{\theta _{2}}$. Let the
following conditions hold.

(i) $\left\langle A_{0}v,Qv\right\rangle \geq c_{A_{0}}\Vert v\Vert _{\theta
-\theta _{2}}^{2}$ for any $v\in V^{\theta -\theta _{2}}$, with a constant $%
c_{A_{0}}>0$;

(ii) $g\in L_{tb}^{2}\left( \mathbb{R}_{+};V^{-\theta -\theta _{2}}\right) $;

Then for some constant $\kappa >0$ independent of time and the initial
condition, we have
\begin{align}
& \Vert u(t)\Vert _{-\theta _{2}}^{2}+\left\Vert d\left( t\right)
\right\Vert _{1}^{2}+\Vert \left( u,d\right) \Vert _{L^{2}(t,t+1;V^{\theta
-\theta _{2}}\times D\left( A_{1}\right) )}^{2}  \label{energy_fbis} \\
& \lesssim e^{-\kappa t}\left( \Vert u(0)\Vert _{-\theta _{2}}^{2}+\mathcal{L%
}\left( \left\Vert d\left( 0\right) \right\Vert _{1}^{2}\right) \right)
+C_{\ast },  \notag
\end{align}%
for all $t\geq 0,$ for some constant $C_{\ast }>0$ and a function $\mathcal{L%
}>0$ independent of time and the initial data.
\end{proposition}

\begin{proof}
The proof of estimate (\ref{energy_fbis}) follows the line of arguments
given in \cite[Proposition 5.1]{GM-JNS}. However the arguments in the
present case are simpler since $f$ is precisely given. For completeness
sake, we include a short proof of this dissipative estimate. As usual, one
proves the claim for smooth approximate solutions and then one passes to the
limit in the end result. First, we observe that we can find a positive
function $\mathcal{L}$ such that
\begin{equation}
\Vert u(t)\Vert _{-\theta _{2}}^{2}+\left\Vert d\left( t\right) \right\Vert
_{1}^{2}\lesssim E_{Q}\left( t\right) \lesssim \Vert u(t)\Vert _{-\theta
_{2}}^{2}+\mathcal{L}\left( \left\Vert d \left( t\right) \right\Vert
_{1}^{2}\right) ,  \label{3.9}
\end{equation}%
owing to the definition of $Q$, the fact that $W\left( d\right) =\frac{1}{4}%
\left( \left\vert d\right\vert ^{2}-1\right) ^{2}$ and standard Sobolev
inequalities. Next, let us set $\rho :=A_{1}d+f\left( d\right) $ and note
that%
\begin{equation*}
\left\Vert A_{1}^{1/2}d\right\Vert _{L^{2}}^{2}+\left\langle f\left(
d\right) d,1\right\rangle =\left\Vert A_{1}^{1/2}d\right\Vert
_{L^{2}}^{2}+\int_{\Omega }\left\vert d\right\vert ^{4}-\int_{\Omega
}\left\vert d\right\vert ^{2}=\left\langle \rho ,d\right\rangle
\end{equation*}%
since $f\left( d\right) =\left( \left\vert d\right\vert ^{2}-1\right) d$. By
simple manipulation of Young and Cauchy-Schwarz inequalities, for any $%
\delta \geq 1/8$, it easily follows%
\begin{align*}
\frac{1}{2}\left\Vert A_{1}^{1/2}d\right\Vert _{L^{2}}^{2}+2\left\langle
W\left( d\right) ,1\right\rangle & =\frac{1}{2}\left\langle \rho
,d\right\rangle -\frac{1}{2}\left\langle \left\vert d\right\vert
^{2}-1,1\right\rangle \\
& \leq \delta \left\Vert \rho \right\Vert _{L^{2}}^{2}+\frac{1}{2}\left(
\frac{1}{8\delta }-1\right) \left\Vert d\right\Vert _{L^{2}}^{2}+\frac{1}{2}%
\left\vert \Omega \right\vert \\
& \leq \delta \left\Vert \rho \right\Vert _{L^{2}}^{2}+\frac{1}{2}\left\vert
\Omega \right\vert .
\end{align*}%
Thus, we obtain%
\begin{equation*}
\frac{1}{2}\left\Vert A_{1}^{1/2}d\right\Vert _{L^{2}}^{2}+\left\langle
W\left( d\right) ,1\right\rangle \leq \frac{3\delta }{2}\left\Vert \rho
\right\Vert _{L^{2}}^{2}+C.
\end{equation*}%
This relation together with either one of the energy identities (\ref%
{en-case1})-(\ref{en-case2}) and the assumption (i), yields for $\kappa \in
\left( 0,1\right) $, the following inequality%
\begin{equation}
\frac{d}{dt}E_{Q}\left( t\right) +\kappa E_{Q}\left( t\right) +\mu
_{1}\left\Vert \left( d^{T}A_{Q}d\right) \left( t\right) \right\Vert
_{L^{2}}^{2}\leq \Theta \left( t\right) ,  \label{diss-ener}
\end{equation}%
where we have set%
\begin{align*}
\Theta \left( t\right) :& =-c_{A_{0}}\left\Vert u\left( t\right) \right\Vert
_{\theta -\theta _{2}}^{2}+\kappa c_{Q}\left\Vert u\left( t\right)
\right\Vert _{-\theta _{2}}^{2}+\frac{1}{\lambda _{1}}\left\Vert \rho \left(
t\right) \right\Vert _{L^{2}}^{2}+\widetilde{\delta }^{-1}\Vert g\left(
t\right) \Vert _{-\theta -\theta _{2}}^{2} \\
& +\widetilde{\delta }\Vert Q\Vert _{-\theta _{2};\theta _{2}}^{2}\Vert
u\left( t\right) \Vert _{\theta -\theta _{2}}^{2}+\frac{3\kappa \delta }{2}%
\left\Vert \rho \left( t\right) \right\Vert _{L^{2}}^{2},
\end{align*}%
for any $\widetilde{\delta }>0$. Setting now $\delta =2/3$ and $\widetilde{%
\delta }>0$ in such a way that $\widetilde{\delta }\Vert Q\Vert _{-\theta
_{2};\theta _{2}}^{2}=c_{A_{0}}/2,$ and exploiting the Sobolev embedding $%
V^{\theta -\theta _{2}}\subseteq V^{-\theta _{2}}$ (i.e., $\left\Vert
u\right\Vert _{-\theta _{2}}^{2}\leq C_{\theta _{2}}\left\Vert u\right\Vert
_{\theta -\theta _{2}}^{2}$ for $\theta \geq 0$), it follows%
\begin{equation*}
\Theta \left( t\right) \leq -\frac{c_{A_{0}}}{2}\left( 1-\frac{2}{c_{A_{0}}}%
\kappa c_{Q}C_{\theta _{2}}\right) \left\Vert u\left( t\right) \right\Vert
_{\theta -\theta _{2}}^{2}-\left( -\lambda _{1}^{-1}-\kappa \right)
\left\Vert \rho \left( t\right) \right\Vert _{L^{2}}^{2}+C\left( 1+\Vert
g\left( t\right) \Vert _{-\theta -\theta _{2}}^{2}\right) ,
\end{equation*}%
for all $t\geq 0$. Adjusting a sufficiently small constant $\kappa \in
\left( 0,\min \left( -\lambda _{1}^{-1}/2,c_{A_{0}}/\left( 4c_{Q}C_{\theta
_{2}}\right) \right) \right) ,$ from (\ref{diss-ener}) we infer%
\begin{align}
& \frac{d}{dt}E_{Q}\left( t\right) +\kappa E_{Q}\left( t\right) +\mu
_{1}\left\Vert \left( d^{T}A_{Q}d\right) \left( t\right) \right\Vert
_{L^{2}}^{2}+C_{\kappa }\left( \left\Vert \rho \left( t\right) \right\Vert
_{L^{2}}^{2}+\left\Vert u\left( t\right) \right\Vert _{\theta -\theta
_{2}}^{2}\right)  \label{diss} \\
& \lesssim 1+\Vert g\left( t\right) \Vert _{-\theta -\theta _{2}}^{2}.
\notag
\end{align}%
The application of Gronwall's inequality (see Appendix, Lemma \ref{Gineq})
in (\ref{diss}) allows us to deduce%
\begin{align}
& E_{Q}\left( t\right) +C_{\kappa }\int_{t}^{t+1}\left( \left\Vert u\left(
s\right) \right\Vert _{\theta -\theta _{2}}^{2}+\left\Vert \rho \left(
s\right) \right\Vert _{L^{2}}^{2}\right) ds  \label{3.9bis} \\
& \leq E_{Q}\left( 0\right) e^{-\kappa t}+C(1+\left\Vert g\right\Vert
_{L_{tb}^{2}\left( \mathbb{R}_{+};V^{-\theta -\theta _{2}}\right) }^{2}),
\notag
\end{align}%
for all $t\geq 0$, for some positive constants $C_{\kappa },C$ independent
of time and the initial data. Reporting (\ref{3.9}) in (\ref{3.9bis}), we
easily arrive at the dissipative estimate (\ref{energy_fbis}). This
completes the proof.
\end{proof}

Next, we recall that by Theorem \ref{t:stab2}, there exists a unique energy
solution%
\begin{equation*}
\left( u,d\right) \in L_{\text{loc}}^{\infty }(0,\infty ;\Upsilon _{\beta
,s})\cap L_{\text{loc}}^{2}(0,\infty ;V^{\beta +\theta }\times
D(A_{1}^{\left( 2s+1\right) /2})),
\end{equation*}%
satisfying (\ref{weak1})-(\ref{weak2}) with any given initial data $\left(
u_{0},d_{0}\right) \in \Upsilon _{\beta ,s}:=V^{\beta }\times D\left(
A_{1}^{s}\right) $. Thus, when the body force $g$ is time independent we can
define a dynamical system for these regular energy solutions. Indeed, system
(\ref{e:op}) generates a semigroup $\left\{ S_{\theta _{2}}(t)\right\}
_{t\geq 0}$ of \emph{closed} operators on the Hilbert space $\Upsilon
_{\beta ,s}$ (when endowed with the metric of $V^{-\theta _{2}}\times W^{1}$%
), given by%
\begin{equation}
\begin{array}{ll}
S_{\theta _{2}}(t):\Upsilon _{\beta ,s}\rightarrow \Upsilon _{\beta ,s},%
\text{ }t\geq 0, &  \\
\left( u_{0},d_{0}\right) \mapsto \left( u\left( t\right) ,d\left( t\right)
\right) . &
\end{array}
\label{semigroup}
\end{equation}

\begin{remark}
\label{disc2}In the case $\lambda _{2}=0,$ $\mu _{1}\geq 0,$ by Theorem \ref%
{t:stab} one can also define the dynamical system $\left( S_{\theta _{2}},%
\mathcal{Y}_{\theta _{2}}\right) $ for problem (\ref{e:op}). In this
instance in (\ref{semigroup}), $\left( u\left( t\right) ,d\left( t\right)
\right) $ is the (unique) energy solution associated with a given initial
datum $\left( u_{0},d_{0}\right) $ in the space $\mathcal{Y}_{\theta _{2}}$.
\end{remark}

In this section, we will only focus on the \emph{general case} when $\mu
_{1}\geq 0$ and $\lambda _{2}\neq 0$ since the former $\lambda _{2}=0,\mu
_{1}\geq 0$ is much easier to handle due the validity of the maximum
principle for $d$, cf. Proposition \ref{maxp}. The following proposition
establishes the existence of an absorbing ball in $\Upsilon _{\beta ,s}$ for
the dynamical system $\left( S_{\theta _{2}},\Upsilon _{\beta ,s}\right) $
in the case $\theta >0,$ $\theta _{2}\geq 1$. Here and everywhere else, $%
\mathcal{B}_{X}\left( R\right) $ denotes the ball in $X$ of radius $R,$
centered at $0$.

\begin{proposition}
\label{t:attr-exist2}Let $s\in \left( \frac{n}{4},1\right] ,$ $n=2,3$ and
consider the following nonempty interval%
\begin{equation*}
J_{n}:=\left( -\theta _{2},\min \left( \theta -\frac{n}{2},\theta -\theta
_{2}\right) \right] \cap \left[ 4s-\theta -3\theta _{2},+\infty \right) .
\end{equation*}%
For $\beta \in J_{n}\backslash \left\{ \theta -n/2\right\} \neq \varnothing $%
, let the following conditions hold.

(i) $b_{0}:V^{\alpha }\times V^{\alpha }\times V^{\theta -\beta }\rightarrow
\mathbb{R}$ is bounded, where $\alpha =\min \{\beta ,\theta -\theta _{2}\}$;

(ii) $b_{0}(v,w,Qw)=0$ for any $v,w\in \mathcal{V}$;

(iii) $g\in V^{\beta -\theta }$ is time independent.

Then for every $R>0$, there exists $t_{\ast }=t_{\ast }\left( R\right) >0$,
such that, for any $\varphi _{0}:=\left( u_{0},d_{0}\right) \in \mathcal{B}%
_{\Upsilon _{\beta ,s}}\left( R\right) ,$%
\begin{equation}
\sup_{t\geq t_{\ast }}\left( \left\Vert \left( u\left( t\right) ,d\left(
t\right) \right) \right\Vert _{\Upsilon _{\beta
,s}}^{2}+\int_{t}^{t+1}\left( \left\Vert u\left( s\right) \right\Vert
_{\theta +\beta }^{2}+\left\Vert d\left( s\right) \right\Vert
_{2s+1}^{2}\right) ds\right) \leq C,  \label{energy_ftris}
\end{equation}%
for some constant $C>0$, independent of time and the initial data.
\end{proposition}

\begin{proof}
By Propositions \ref{t:attr-exist}, for every $R>0$ with $\varphi _{0}\in
\mathcal{B}_{\Upsilon _{\beta ,s}}\left( R\right) $ there exists $%
t_{0}=t_{0}(R)>0$ such that%
\begin{equation}
\sup_{t\geq t_{0}}\left\Vert \left( u\left( t\right) ,d\left( t\right)
\right) \right\Vert _{\mathcal{Y}_{\theta _{2}}}^{2}+\int_{t}^{t+1}\left(
\left\Vert u\left( s\right) \right\Vert _{\theta -\theta
_{2}}^{2}+\left\Vert A_{1}d\left( s\right) \right\Vert _{L^{2}}^{2}\right)
ds\leq C_{\ast }.  \label{energy-fbis2}
\end{equation}%
Moreover, by Theorem \ref{t:reg} and application of the uniform Gronwall's
lemma \cite[Lemma III.1.1]{T}\ in (\ref{inn3}), by virtue of (\ref%
{energy-fbis2}), we infer the existence of a new time $t_{\ast }=t_{0}+1$
such that%
\begin{equation}
\sup_{t\geq t_{\ast }}\left( \left\Vert u\left( t\right) \right\Vert _{\beta
}^{2}+\left\Vert d\left( t\right) \right\Vert _{2s}^{2}\right) \leq C,
\label{comp-abs}
\end{equation}%
for some positive constant $C$ independent of time and the initial data.
Moreover, integration over $\left( t,t+1\right) $ of the inequality (\ref%
{inn3}) yields%
\begin{equation}
\sup_{t\geq t_{\ast }}\int_{t}^{t+1}\left( \left\Vert u\left( s\right)
\right\Vert _{\beta +\theta }^{2}+||d\left( s\right) ||_{2s+1}^{2}\right)
\leq C,  \label{comp-abs-bis}
\end{equation}%
owing once again to (\ref{comp-abs}). The claim (\ref{energy_ftris}) is then
immediate.
\end{proof}

Next we show the existence of finite dimensional global attractors for our
regularized family of models (\ref{e:op}) when $\theta >0$. However, due to
lack of compactness of the solutions in the space $\Upsilon _{\beta ,s}$ we
cannot proceed in a standard way. Indeed, the strong coupling in the full
Ericksen-Leslie system (\ref{e:op}) for $\left( u,d\right) $ prevents us
from establishing any additional smoothing properties of the solutions
without requiring more restrictive assumptions on the body force $g$ and the
other parameters of the problem. In fact, in what follows we shall prove
even more: the existence of an exponential attractor for $\left( S_{\theta
_{2}},\Upsilon _{\beta ,s}\right) .$ We recall that the exponential
attractor always contains the global attractor and also attracts bounded
subsets of the energy phase-space at an exponential rate, which makes it a
more useful object in numerical simulations than the global attractor.

We shall accomplish this program in a series of subsequent lemmas. First, we
have the basic statement.

\begin{proposition}
\label{basic-inv}Let $\mathcal{B}_{0}$ be a bounded absorbing ball whose
existence has been proven in Proposition \ref{t:attr-exist2}. The set%
\begin{equation*}
\mathcal{B}_{\ast }:=\bigcup_{t\geq t_{\ast }}S_{\theta _{2}}\left( t\right)
\mathcal{B}_{0}
\end{equation*}%
is bounded in $\Upsilon _{\beta ,s}$ and positively invariant for $S_{\theta
_{2}}.$
\end{proposition}

Clearly, $\mathcal{B}_{\ast }$ is also absorbing for the semigroup $%
S_{\theta _{2}}$. Thus, it is sufficient to construct the exponential
attractor for the restriction of this semigroup on $\mathcal{B}_{\ast }$
only. With this at hand, we can show the uniform H\"{o}lder continuity of $%
t\mapsto S_{\theta _{2}}\left( t\right) \varphi _{0}$ in the $V^{-\theta
_{2}}\times W^{1}$-norm, namely,

\begin{lemma}
\label{hctime} Let the assumptions of Proposition \ref{t:attr-exist2} be
satisfied. Consider $\varphi \left( t\right) =S_{\theta _{2}}\left( t\right)
\varphi _{0}$ with $\varphi _{0}\in \mathcal{B}_{\ast }$. Then, we have%
\begin{equation}
\left\Vert u\left( t\right) -u\left( \widetilde{t}\right) \right\Vert
_{-\theta _{2}}+\left\Vert d\left( t\right) -d\left( \widetilde{t}\right)
\right\Vert _{1}\lesssim (\left\vert t-\widetilde{t}\right\vert ^{\epsilon
_{1}}+\left\vert t-\widetilde{t}\right\vert ^{\epsilon _{2}}),
\label{hclinf2}
\end{equation}%
for all $t,\widetilde{t}\in \left[ 0,T\right] $, for some $0<\epsilon
_{1},\epsilon _{2}<1$ depending only on $\theta ,\theta _{2}$ and $s\in (%
\frac{n}{4},1].$
\end{lemma}

\begin{proof}
We can rely once again on the proof of Theorem \ref{t:exist} and Theorem \ref%
{t:stab}. Indeed, since the $V^{\beta }$-norm of $u$ and the $W^{2s}$-norm
of $d$ are globally bounded by Theorem \ref{t:reg} if $\varphi _{0}\in
\mathcal{B}_{\ast }$, then also%
\begin{equation*}
\partial _{t}\varphi =\left( \partial _{t}u,\partial _{t}d\right) \in L_{%
\text{loc}}^{2}\left( 0,\infty ;V^{-\theta -\theta _{2}}\times L^{2}\left(
\Omega \right) \right) .
\end{equation*}%
Finally, the simple relation%
\begin{equation*}
\varphi \left( t\right) -\varphi \left( \widetilde{t}\right) =\int_{t}^{%
\widetilde{t}}\partial _{s}\varphi \left( s\right) ds
\end{equation*}%
and proper interpolation inequalities in the spaces $V^{\beta }\subset
V^{-\theta _{2}}\subset V^{-\theta -\theta _{2}},$ $W^{2s}\subset
W^{1}\subset L^{2}\left( \Omega \right) ,$ $2s>n/2$, imply the desired
inequality (\ref{hclinf2}).
\end{proof}

The crucial step in order to establish the existence of an exponential
attractor is the validity of so-called smoothing property for the difference
of any two energy solutions $\varphi _{i},$ \thinspace $i=1,2$. In the
present case, such a property is a consequence of the following two lemmas.
The first result establishes that the semigroup $S_{\theta _{2}}\left(
t\right) $ is some kind of contraction map, up to the term $\left\Vert
\varphi _{1}-\varphi _{2}\right\Vert _{L^{2}(0,t;\mathcal{Y}_{\theta _{2}})}$%
.

\begin{lemma}
\label{dec} Let the assumptions of Proposition \ref{t:attr-exist2} hold. For
any two energy solutions $\varphi _{i}=\left( u_{i},d_{i}\right) $\
associated with the initial data $\varphi _{0i}\in \mathcal{B}_{\ast }$, the
following estimate holds:%
\begin{align}
& \left\Vert u_{1}\left( t\right) -u_{2}\left( t\right) \right\Vert
_{-\theta _{2}}^{2}+\left\Vert d_{1}\left( t\right) -d_{2}\left( t\right)
\right\Vert _{1}^{2}  \label{decay1} \\
& \lesssim e^{-\eta t}\left( \left\Vert u_{1}\left( 0\right) -u_{2}\left(
0\right) \right\Vert _{-\theta _{2}}^{2}+\left\Vert d_{1}\left( 0\right)
-d_{2}\left( 0\right) \right\Vert _{1}^{2}\right)  \notag \\
& +\int_{0}^{t}\left( \left\Vert u_{1}\left( s\right) -u_{2}\left( s\right)
\right\Vert _{-\theta _{2}}^{2}+\left\Vert d_{1}\left( s\right) -d_{2}\left(
s\right) \right\Vert _{1}^{2}\right) ds,  \notag
\end{align}%
for all $t\geq 0$, for some positive constant $\eta $ independent of time.
\end{lemma}

\begin{proof}
We rely mainly on the estimates exploited in the proof of Theorem \ref%
{t:stab} and Theorem \ref{t:stab2}. Indeed, each energy solution $\varphi
_{i}$ is globally bounded in $L^{\infty }\left( 0,\infty ;L^{\infty }\left(
\Omega \right) \right) $ by Remark \ref{rem-bd} if $\varphi _{0i}\in
\mathcal{B}_{\ast }$. It turns out that the main steps require nothing more
than what is already contained in the proof of Theorem \ref{t:stab} (or
Theorem \ref{t:stab2}).

Our starting point is the inequality (\ref{uniqest1bis}). With the exception
of $I_{2},I_{63},$ $I_{64},$ all the other terms can be estimated word by
word as in (\ref{uniqest2}), (\ref{uniqest61}), (\ref{uniqest62}), (\ref%
{uniqest66}) and (\ref{uniqest5}), respectively. For $I_{2},I_{63}$ and $%
I_{64}$, we need more refined estimates. We have the bounds:%
\begin{eqnarray}
\left\vert I_{2}\right\vert &=&\left\vert \left\langle
A_{1}^{1-s}B_{1}\left( u,d\right) ,A_{1}^{s}d_{2}\right\rangle \right\vert
\label{uniqest7} \\
&\lesssim &\left\Vert d_{2}\right\Vert _{2s}\left\Vert Qu\right\Vert
_{\theta _{2}}\left\Vert \nabla d\right\Vert _{1}  \notag \\
&\leq &\delta \left\Vert A_{1}d\right\Vert _{L^{2}}^{2}+C\left\Vert
u\right\Vert _{-\theta _{2}}^{2}\left\Vert d_{2}\right\Vert _{2s}^{2}  \notag
\end{eqnarray}%
and%
\begin{align}
\left\vert I_{63}\right\vert & =\left\vert \mu _{2}\left\langle \widehat{%
\mathcal{N}}\otimes d_{1}-\mathcal{N}_{2}\otimes d,\nabla Qu\right\rangle
\right\vert  \label{uniqest8} \\
& \lesssim \left\Vert \nabla Qu\right\Vert _{\theta _{2}-1}\left( ||\widehat{%
\mathcal{N}}\otimes d_{1}||_{1-\theta _{2}}+\left\Vert \mathcal{N}%
_{2}\otimes d\right\Vert _{1-\theta _{2}}\right)  \notag \\
& \leq C_{M}\left\Vert u\right\Vert _{-\theta _{2}}\left( \left\Vert
A_{1}d\right\Vert _{L^{2}}+\left\Vert u_{2}\right\Vert _{-\theta
_{2}}\left\Vert d_{2}\right\Vert _{2s}\left\Vert A_{1}d\right\Vert
_{L^{2}}\right)  \notag \\
& \leq \delta \left\Vert A_{1}d\right\Vert _{L^{2}}^{2}+C_{\delta
,M}\left\Vert u\right\Vert _{-\theta _{2}}^{2}\left( 1+\left\Vert
u_{2}\right\Vert _{-\theta _{2}}^{2}\left\Vert d_{2}\right\Vert
_{2s}^{2}\right) .  \notag
\end{align}%
In detail, these bounds are deduced using the definition of $\widehat{%
\mathcal{N}}$ and $\mathcal{N}_{2}$ together with the following crucial
properties:

\begin{itemize}
\item[(a1)] Each term $Qu_{i}\left( \partial _{i}d_{j}\right) $ is a product
of functions in $H^{\theta _{2}}$ and $H^{1}$ and therefore bounded in $%
H^{2-2s},$ by Lemma \ref{l:hole} (Appendix) since $2s>n/2$ and $\theta
_{2}\geq 1.$ This yields (\ref{uniqest7}).

\item[(a2)] For $I_{63}$, we observe that each term $\left[ \left( \partial
_{i}Qu_{2j}\right) d_{j}\right] d_{l}$ is a product of functions from $H^{0}$
and $H^{2}\subset L^{\infty }$ and therefore bounded in $H^{0}\subseteq
H^{1-\theta _{2}}.$ Finally, the terms $\left( \partial _{i}Qu_{2j}\right)
d_{j}$ are products of functions in $H^{\theta _{2}-1}$ and $H^{2s}\subset
L^{\infty },$ and therefore bounded in $L^{2}.$
\end{itemize}

\noindent A similar argument to the derivation of (\ref{uniqest8}) yields
the bound%
\begin{equation*}
\left\vert I_{64}\right\vert \leq \delta \left\Vert A_{1}d\right\Vert
_{L^{2}}^{2}+C_{\delta ,M}\left\Vert u\right\Vert _{-\theta _{2}}^{2}\left(
1+\left\Vert u_{2}\right\Vert _{-\theta _{2}}^{2}\left\Vert d_{2}\right\Vert
_{2s}^{2}\right) .
\end{equation*}%
Since also $\left( u_{i},d_{i}\right) \in L^{\infty }\left( 0,\infty
;V^{\beta }\times W^{2s}\right) $, from (\ref{uniqest1bis}) we finally see
that%
\begin{eqnarray*}
&&\frac{d}{dt}\left( \Vert u\left( t\right) \Vert _{-\theta
_{2}}^{2}+\left\Vert d\left( t\right) \right\Vert _{1}^{2}\right) +\left(
2c_{A_{0}}-6\delta \right) \Vert u\left( t\right) \Vert _{\theta -\theta
_{2}}^{2}+\left( -2\lambda _{1}^{-1}-10\delta \right) \left\Vert
A_{1}d\left( t\right) \right\Vert _{L^{2}}^{2} \\
&\leq &C_{M,\delta }\left( \Vert u\left( t\right) \Vert _{-\theta
_{2}}^{2}+\left\Vert d\left( t\right) \right\Vert _{1}^{2}\right) ,
\end{eqnarray*}%
for all $t\geq 0$, for a sufficiently small $\delta \in \left( 0,\min \left(
c_{A_{0}}/3,-\lambda _{1}^{-1}/5\right) \right) $. Thus, Gronwall's
inequality entails the desired estimate (\ref{decay1}).
\end{proof}

We now need some compactness for the term $\left\Vert \varphi _{1}-\varphi
_{2}\right\Vert _{L^{2}\left( 0,t;\mathcal{Y}_{\theta _{2}}\right) }$ on the
right-hand side of (\ref{decay1}). This is given by

\begin{lemma}
\label{lipdif} Let the assumptions of Proposition \ref{t:attr-exist2} hold.
Then, the following estimate holds:%
\begin{align}
& \left\Vert \partial _{t}u_{1}-\partial _{t}u_{2}\right\Vert _{L^{2}\left(
0,t;V^{-\theta -\theta _{2}}\right) }^{2}+\left\Vert \partial
_{t}d_{1}-\partial _{t}d_{2}\right\Vert _{L^{2}\left( 0,t;L^{2}\left( \Omega
\right) \right) }^{2}  \label{comp1} \\
& +\int_{0}^{t}\left( \left\Vert u_{1}\left( s\right) -u_{2}\left( s\right)
\right\Vert _{-\theta _{2}}^{2}+\left\Vert d_{1}\left( s\right) -d_{2}\left(
s\right) \right\Vert _{1}^{2}\right) ds  \notag \\
& \leq \varrho \left( t\right) \left( \left\Vert u_{1}\left( 0\right)
-u_{2}\left( 0\right) \right\Vert _{-\theta _{2}}^{2}+\left\Vert d_{1}\left(
0\right) -d_{2}\left( 0\right) \right\Vert _{1}^{2}\right) ,  \notag
\end{align}%
for all $t\geq 0$.
\end{lemma}

\begin{proof}
The required control of the integral term on the left-hand side of (\ref%
{comp1}) is readily provided by (\ref{c-dep}). It remains to gain some
control on the time derivative $\left( \partial _{t}u,\partial _{t}d\right)
, $ where $u:=u_{1}-u_{2},$ $d:=d_{1}-d_{2}.$ We recall the variational
formulation (\ref{diffuniq1})-(\ref{diffuniq2}) and rely once again on the
fact that%
\begin{equation}
\varphi _{i}=\left( u_{i},d_{i}\right) \in L^{\infty }\left( 0,\infty ;%
\mathcal{B}_{\ast }\right) ,\text{ }d_{i}\in L^{\infty }\left( 0,\infty
;L^{\infty }\left( \Omega \right) \right) ,  \label{very-g}
\end{equation}%
for each $i=1,2$. For any test functions $w\in V^{\theta +\theta _{2}}$ and $%
\eta \in L^{2}\left( \Omega \right) $, using the corresponding variational
formulation, one has%
\begin{equation}
\begin{array}{l}
\left\vert \left\langle \partial _{t}u,w\right\rangle \right\vert \leq
\left( \left\Vert \text{r.h.s.u}\right\Vert _{-\theta -\theta
_{2}}+\left\Vert \sigma _{Qu_{1}}-\sigma _{Qu_{2}}\right\Vert _{1-\theta
-\theta _{2}}\right) \left\Vert w\right\Vert _{\theta +\theta _{2}}, \\
\left\vert \left\langle \partial _{t}d,\eta \right\rangle \right\vert \leq
\left\Vert \text{r.h.s.d}\right\Vert _{L^{2}}\left\Vert \eta \right\Vert
_{L^{2}},%
\end{array}
\label{very-g2}
\end{equation}%
with%
\begin{equation*}
\text{r.h.s.u}:=-A_{0}u-B_{0}\left( u,u_{1}\right) -B_{0}\left(
u_{2},u\right) +R_{0}\left( A_{1}d_{2},d\right) +R_{0}\left(
A_{1}d,d_{1}\right)
\end{equation*}%
and%
\begin{align*}
\text{r.h.s.d}& :=\lambda _{1}^{-1}A_{1}d-B_{1}\left( u,d_{1}\right)
-B_{1}\left( u_{2},d\right) \\
& +\lambda _{1}^{-1}\left( f\left( d_{1}\right) -f\left( d_{2}\right)
\right) +\omega _{Qu}d_{1}+\omega _{Qu_{2}}d \\
& -\lambda _{2}\lambda _{1}^{-1}\left( A_{Qu}d_{1}+A_{Qu_{2}}d\right) .
\end{align*}%
Repeated use of (\ref{very-g}) and arguing as in the proof of Theorem \ref%
{t:stab}, owing to $\theta >0$, $\theta _{2}\geq 1,$ it is now
straightforward to show that%
\begin{align*}
& \left\Vert \text{r.h.s.u}\right\Vert _{-\theta -\theta _{2}}+\left\Vert
\text{r.h.s.d}\right\Vert _{L^{2}}+\left\Vert \sigma _{Qu_{1}}-\sigma
_{Qu_{2}}\right\Vert _{1-\theta -\theta _{2}} \\
& \leq C\left( \left\Vert u\right\Vert _{\theta -\theta _{2}}+\left\Vert
A_{1}d\right\Vert _{L^{2}}\right) ,
\end{align*}%
for some constant $C>0$ which depends on $\mathcal{B}_{\ast }$, but is
independent of time. This estimate together with (\ref{very-g2}) and (\ref%
{c-dep}) gives the desired estimate on the time derivatives in (\ref{comp1}%
). The proof is finished.
\end{proof}

The main result of this section is concerned with the existence of
exponential attractors for problem (\ref{e:op}) in the case $\theta >0$.

\begin{theorem}
\label{expo-thm}Let the assumptions of Proposition \ref{t:attr-exist2} be
satisfied. Then the dynamical system $\left( S_{\theta _{2}},\Upsilon
_{\beta ,s}\right) $ possesses an exponential attractor $\mathcal{M}_{\theta
_{2},\beta ,s}\subset \Upsilon _{\beta ,s}$ which is bounded in $\Upsilon
_{\beta ,s}$. More precisely by definition, we have

\begin{itemize}
\item[(a)] $\mathcal{M}_{\theta _{2},\beta ,s}$ is compact and
semi-invariant with respect $S_{\theta _{2}}\left( t\right) ,$ that is,%
\begin{equation*}
S_{\theta _{2}}\left( t\right) \left( \mathcal{M}_{\theta _{2},\beta
,s}\right) \subseteq \mathcal{M}_{\theta _{2},\beta ,s},\quad \forall
\,t\geq 0.
\end{equation*}

\item[(b)] The fractal dimension $\dim _{F}\left( \mathcal{M}_{\theta
_{2},\beta ,s},\mathcal{Y}_{\theta _{2}}\right) $ of $\mathcal{M}_{\theta
_{2},\beta ,s}$ is finite and an upper bound can be computed explicitly.

\item[(c)] $\mathcal{M}_{\theta _{2},\beta ,s}$ attracts exponentially fast
any bounded subset $B$ of $\Upsilon _{\beta ,s}$, that is, there exist a
positive nondecreasing function $\mathcal{L}$ and a constant $\tau >0$ such
that
\begin{equation*}
dist_{\mathcal{Y}_{\theta _{2}}}\left( S_{\theta _{2}}\left( t\right) B,%
\mathcal{M}_{\theta _{2},\beta ,s}\right) \leq \mathcal{L}(\Vert B\Vert
_{\Upsilon _{\beta ,s}})e^{-\tau t},\quad \forall \,t\geq 0.
\end{equation*}
\end{itemize}

\noindent Here $dist_{\mathcal{Y}_{\theta _{2}}}$ denotes the Hausdorff
semi-distance between sets in $\mathcal{Y}_{\theta _{2}}$ and $\Vert B\Vert
_{\Upsilon _{\beta ,s}}$ stands for the size of $B$ in $\Upsilon _{\beta
,s}. $ Both $\mathcal{L}$ and $\tau $ can be explicitly calculated.
\end{theorem}

\begin{proof}
We apply an abstract result stated in the Appendix, see Proposition \ref%
{abstract}. Recall that by Proposition \ref{basic-inv}, the ball $\mathcal{B}%
_{\ast }$ is absorbing and positively invariant for $S_{\theta _{2}}\left(
t\right) $. On the other hand, due to the results proven in this section, we
have%
\begin{equation*}
\sup_{t\geq 0}\left\Vert \left( u\left( t\right) ,d\left( t\right) \right)
\right\Vert _{\Upsilon _{\beta ,s}}\leq C_{\beta ,s},
\end{equation*}%
for every trajectory $\varphi =\left( u,d\right) $ originating from $\varphi
_{0}=\left( u_{0},d_{0}\right) \in \mathcal{B}_{\ast }$, for some positive
constant $C_{\beta ,s}$ which is independent of the choice of $\varphi
_{0}\in \mathcal{B}_{\ast }$. We can now apply the abstract result of
Proposition \ref{abstract} to the map%
\begin{equation*}
\mathbb{S}=S_{\theta _{2}}\left( T_{\ast }\right) :\mathbb{B\rightarrow B},
\end{equation*}%
where $\mathbb{B}=\mathcal{B}_{\ast }$ and $\mathcal{H}=V^{-\theta
_{2}}\times W^{1}$, for a fixed $T_{\ast }>0$ such that $e^{-\eta T_{\ast }}<%
\frac{1}{2}$, $\eta >0$ is the same as in Lemma \ref{dec}. To this end, we
introduce the functional spaces%
\begin{equation}
\begin{array}{l}
\mathcal{V}_{1}:=L^{2}\left( 0,T;V^{\theta -\theta _{2}}\times L^{2}\left(
\Omega \right) \right) \cap H^{1}\left( 0,T;V^{-\theta -\theta _{2}}\times
L^{2}\left( \Omega \right) \right) , \\
\mathcal{V}:=L^{2}\left( 0,T;V^{-\theta _{2}}\times W^{1}\right)%
\end{array}
\label{fsp}
\end{equation}%
and note that $\mathcal{V}_{1}$ is compactly embedded into $\mathcal{V}$ due
to the Aubin-Lions-Simon compactness lemma. Finally, we introduce the
operator $\mathbb{T}:\mathcal{B}_{\ast }\rightarrow \mathcal{V}_{1}$, by $%
\mathbb{T}\varphi _{0}:=\varphi \in \mathcal{V}_{1},$ where $\varphi $
solves (\ref{e:op}) with $\varphi \left( 0\right) =\varphi _{0}\in \mathcal{B%
}_{\ast }$. We claim that the maps $\mathbb{S}$, $\mathbb{T}$, the spaces $%
\mathcal{H}$,$\mathcal{V}$,$\mathcal{V}_{1}$ thus defined satisfy all the
assumptions of Proposition \ref{abstract}. Indeed, the global Lipschitz
continuity (\ref{gl1}) of $\mathbb{T}$ is an immediate corollary of Lemma %
\ref{lipdif}, and estimate (\ref{gl2}) follows from estimate (\ref{decay1}).
Therefore, due to Proposition \ref{abstract}, the semigroup $\mathbb{S}%
(n)=S_{\theta _{2}}\left( nT_{\ast }\right) $ generated by the iterations of
the operator $\mathbb{S}:\mathcal{B}_{\ast }\rightarrow \mathcal{B}_{\ast }$
possesses a (discrete) exponential attractor $\left( \mathcal{M}_{\theta
_{2},\beta ,s}\right) _{d}$ in $\mathcal{B}_{\ast }$ endowed with the
topology of $V^{-\theta _{2}}\times W^{1}$. In order to construct the
exponential attractor $\mathcal{M}_{\theta _{2},\beta ,s}$ for the semigroup
$S_{\theta _{2}}(t)$ with continuous time, we note that, due to Theorem \ref%
{t:stab2}, this semigroup is Lipschitz continuous with respect to the
initial data in the topology of $V^{-\theta _{2}}\times W^{1}$. Moreover, by
Lemma \ref{hctime} the map $\left( t,\varphi _{0}\right) \mapsto S_{\theta
_{2}}\left( t\right) \varphi _{0}$ is also uniformly H\"{o}lder continuous
on $\left[ 0,T\right] \times \mathcal{B}_{\ast }$, where $\mathcal{B}_{\ast
} $ is equipped with the metric topology of $V^{-\theta _{2}}\times W^{1}$.
Hence, the desired exponential attractor $\mathcal{M}_{\theta _{2},\beta ,s}$
for the continuous semigroup $S_{\theta _{2}}(t)$ can be obtained by the
standard formula%
\begin{equation}
\mathcal{M}_{\theta _{2},\beta ,s}=\bigcup_{t\in \left[ 0,T_{\ast }\right]
}S_{\theta _{2}}\left( t\right) \left( \mathcal{M}_{\theta _{2},\beta
,s}\right) _{d}.  \label{st}
\end{equation}%
Theorem \ref{expo-thm} is now proved.
\end{proof}

As a consequence of the above theorem, we have the following.

\begin{corollary}
\label{global-thm}Under the assumptions of Theorem \ref{expo-thm}, there
exists a global attractor $\mathcal{A}_{\theta _{2},\beta ,s}$ which
attracts the bounded sets of $\Upsilon _{\beta ,s}$. Moreover, $\mathcal{A}%
_{\theta _{2},\beta ,s}$ is connected, bounded in $\Upsilon _{\beta ,s}$ and
$\mathcal{A}_{\theta _{2},\beta ,s}$ has finite fractal dimension:%
\begin{equation*}
\dim _{F}\left( \mathcal{A}_{\theta _{2},\beta ,s},\mathcal{Y}_{\theta
_{2}}\right) <\infty .
\end{equation*}
\end{corollary}

\begin{remark}
In fact due to interpolation, Theorem \ref{expo-thm} also implies that the
fractal dimension of the global and exponential attractors is finite in $%
V^{\beta -\varepsilon _{1}}\times W^{2s-\varepsilon _{2}}$, for every $%
-\theta _{2}<\varepsilon _{1}<\beta $ and $n/2<\varepsilon _{2}<2s$. The
attraction property in (c) also holds in the stronger topology of $V^{\beta
-\varepsilon _{1}}\times W^{2s-\varepsilon _{2}}$.
\end{remark}

Note that Proposition \ref{t:attr-exist2} provides many examples where the
conclusion of Theorem \ref{expo-thm} is satisfied. For example, checking all
the requirements of Proposition \ref{t:attr-exist2} in the case when $\theta
=1$, the conclusions of Theorem \ref{expo-thm} and Corollary \ref{global-thm}
are satisfied for the modified 3D Leray-EL-$\alpha $ (ML-EL-$\alpha $)
model, the 3D SBM-EL model and the 3D NS-EL-$\alpha $ system. These results
were not reported anywhere in the literarure for the full Ericksen-Leslie
model.

\section{Convergence to steady states}

\label{s:convss}

In this section, we show that any global-in-time bounded energy solution to
the full regularized or nonregularized Ericksen-Leslie model (\ref{e:op})
converges (in a certain sense) to a single equilibrium as time tends to
infinity. The proof of the main statements are based on a suitable version
of the Lojasiewicz--Simon theorem and the results developed in the previous
sections. We emphasize that our subsequent results hold \emph{only} for the
energy weak solutions introduced through Definition \ref{weak-energy}, even
when uniqueness is not available. In particular, they hold for limit points
of the Galerkin approximation scheme exploited in Theorem \ref{t:exist}, as
well as for other approximation schemes in which the energy inequality (\ref%
{en-both})\ can be proven. Thus, the energy inequality is crucial for
investigating the long-time behavior as time goes to infinity. It will also
serve as a selection criterion in eliminating all those \emph{non-physical}
weak solutions in the framework of Definition \ref{weak}, which may not
necessarily satisfy the energy inequality (\ref{en-both}). Finally, in some
cases when the energy solutions become more regular, we can also prove
stronger convergence results.

We shall first prove that every energy solution given by Definition \ref%
{weak-energy} has a non-empty $\omega $-limit set.

\begin{lemma}
\label{ss:omega}Let the assumptions of Theorem \ref{t:exist} be satisfied,
and suppose that $g$ also obeys the following condition:%
\begin{equation}
\int_{t}^{\infty }\left\Vert g\left( s\right) \right\Vert _{-\theta -\theta
_{2}}^{2}ds\lesssim \left( 1+t\right) ^{-\left( 1+\delta \right) },\text{
for all }t\geq 0,  \label{g-assumpt}
\end{equation}%
for some constant $\delta \in \left( 0,1\right) $. Let $\left( u,d\right) $
be an energy solution in the sense of Definition \ref{weak-energy}. Then,
the $\omega $-limit set of $\left( u,d\right) $ is nonempty. More precisely,
we have%
\begin{equation}
\lim_{t\rightarrow \infty }u\left( t\right) =0\text{ weakly in }V^{-\theta
_{2}}  \label{o1}
\end{equation}%
and any divergent sequence $\left\{ t_{n}\right\} \subset \left[ 0,\infty
\right) $ admits a subsequence, denoted by $\left\{ t_{n_{k}}\right\} $,
such that%
\begin{equation}
\lim_{t_{n_{k}}\rightarrow \infty }d\left( t_{n_{k}}\right) =d_{\ast }\text{
weakly in }W^{1}\text{, strongly in }W^{0},  \label{o2}
\end{equation}%
for some $d_{\ast }\in D\left( A_{1}\right) $ which is a solution of%
\begin{equation}
A_{1}d_{\ast }+f\left( d_{\ast }\right) =0\text{ in }\Omega .  \label{o3}
\end{equation}
\end{lemma}

\begin{proof}
First, we recall that an energy solution in the sense of Definition \ref%
{weak-energy} exists by virtue of Theorem \ref{t:exist}. Our proof follows
the lines of the argument given in \cite[Theorem 2.6]{prslong} and our
arguments developed in Theorem \ref{t:exist}. We prove our subsequent
results in \textbf{Case 1 }(\textbf{Case 2 }is analogous and follows with
some minor modifications). To this end, let $\left\{ t_{n}\right\} \subset %
\left[ 0,\infty \right) $ be a divergent sequence. The energy inequality (%
\ref{en-both}) together with assumption (\ref{g-assumpt}) implies, at least
for a suitable subsequence of $\left\{ t_{n}\right\} $, still labelled as $%
\left\{ t_{n}\right\} $, that%
\begin{equation*}
u\left( t_{n}\right) \rightarrow u_{\ast }\text{ weakly in }V^{-\theta _{2}},%
\text{ }d\left( t_{n}\right) \rightarrow d_{\ast }\text{ weakly in }W^{1},
\end{equation*}%
for some $\left( u_{\ast },d_{\ast }\right) \in \mathcal{Y}_{\theta _{2}}$.
Consider now the initial value problem (\ref{e:op}) on the time interval $%
\left[ t_{n},t_{n+1}\right] $ with the initial values $\left( u\left(
t_{n}\right) ,d\left( t_{n}\right) \right) $ and observe that $\left(
u_{n}\left( t\right) ,d_{n}\left( t\right) \right) :=\left( u\left(
t+t_{n}\right) ,d\left( t+t_{n}\right) \right) $ are also weak solutions of (%
\ref{e:op}) for $t\in \left[ 0,1\right] $. Then, from the energy inequality (%
\ref{en-both}) and (\ref{apriori-tensor}), as $t_{n}\rightarrow \infty $ we
infer%
\begin{equation}
\begin{array}{l}
u_{n}\rightarrow 0\text{ strongly in }L^{2}\left( 0,1;V^{\theta -\theta
_{2}}\right) ,\text{ weakly star in }L^{\infty }\left( 0,1;V^{-\theta
_{2}}\right) , \\
d_{n}\rightarrow \overline{d}\text{ weakly in }L^{2}\left( 0,1;W^{2}\right) ,%
\text{ weakly star in }L^{\infty }\left( 0,1;W^{1}\right) ,%
\end{array}
\label{cs1}
\end{equation}%
for a suitable function $\overline{d}$. Repeating the comparison arguments
developed in the proof of Theorem \ref{t:exist}, we deduce%
\begin{equation}
\partial _{t}d_{n}\rightarrow \partial _{t}\overline{d}\text{ weakly in }%
L^{2}\left( 0,1;W^{-2}\right)  \label{cs2}
\end{equation}%
and, in particular, $\partial _{t}u_{n}$ is uniformly bounded in $%
L^{p}\left( 0,1;V^{-\gamma }\right) $ for $p>1$ and $\gamma \geq 0$ as given
by (\ref{thep}). In particular, by the Aubin-Lions-Simon compactness
criterion and (\ref{cs1})-(\ref{cs2}) we obtain%
\begin{equation}
d_{n}\rightarrow \overline{d}\text{ strongly in }L^{2}\left(
0,1;W^{1}\right) \cap C\left( 0,1;W^{-2}\right) ,  \label{cs3}
\end{equation}%
as well as%
\begin{equation}
u_{n}\rightarrow 0\text{ strongly in }C\left( 0,1;V^{-\gamma }\right)
\label{cs4}
\end{equation}%
for some $\gamma \geq 0$. Moreover, by the definition of $Q$, one has%
\begin{equation}
Qu_{n}\rightarrow 0\text{ strongly in }L^{2}\left( 0,1;V^{\theta +\theta
_{2}}\right) ,\text{ weakly star in }L^{\infty }\left( 0,1;V^{\theta
_{2}}\right) .  \label{cs4bis}
\end{equation}%
Thus, (\ref{cs4}) yields that $u_{\ast }=0$, which implies (\ref{o1}) in
view of (\ref{cs1}) and (\ref{cs4}). On the other hand, by the energy
inequality (\ref{en-both}) and (\ref{apriori-tensor}) we also have%
\begin{align}
A_{1}d_{n}+f\left( d_{n}\right) & \rightarrow 0\text{ strongly in }%
L^{2}\left( 0,1;L^{2}\left( \Omega \right) \right) ,  \label{cs5} \\
A_{Q}^{n}d_{n}& \rightarrow 0\text{ strongly in }L^{2}\left( 0,1;L^{2}\left(
\Omega \right) \right) ,  \notag \\
d_{n}^{T}A_{Q}^{n}d_{n}& \rightarrow 0\text{ strongly in }L^{2}\left(
0,1;L^{2}\left( \Omega \right) \right) .  \notag
\end{align}%
Next, by the estimates (\ref{cs1}), (\ref{cs3})-(\ref{cs4}), and observing
that terms like $\left( \partial _{i}Qu_{j}\right) d_{j}$ are a product of
functions in $H^{\theta +\theta _{2}-1}$ and $H^{1},$ and therefore bounded
in $H^{-2},$ for any $\theta ,\theta _{2}\geq 0$ by Lemma \ref{l:hole}, we
also infer%
\begin{equation}
\omega _{Q}^{n}d_{n}\rightarrow 0\text{ strongly in }L^{2}\left(
0,1;W^{-2}\right) .  \label{cs6}
\end{equation}%
Using (\ref{cs4bis}), a similar argument entails that%
\begin{equation}
B_{1}\left( u_{n},d_{n}\right) \rightarrow 0\text{ strongly in }L^{2}\left(
0,1;W^{-2}\right) .  \label{cs7}
\end{equation}%
Thus, comparing terms in the second equation of (\ref{e:op}), we also obtain%
\begin{equation}
\partial _{t}d_{n}\rightarrow 0\text{ strongly in }L^{2}\left(
0,1;W^{-2}\right) ;  \label{cs8}
\end{equation}%
henceforth, it follows that $\overline{d}=d_{\ast }$ for all $t\in \left[ 0,1%
\right] .$ Letting now $t_{n}\rightarrow \infty $ in the equation for the
director field $d_{n}$, satisfying (\ref{e:op}), we observe that $d_{\ast }$
is also a solution of $A_{1}d_{\ast }+f\left( d_{\ast }\right) =0$ in $%
\Omega $ and $d_{\ast }\in D\left( A_{1}\right) $, as claimed. Lemma \ref%
{ss:omega} is proved.
\end{proof}

Even though we are dealing with an asymptotically decaying force $g\left(
t\right) $ due to (\ref{g-assumpt}), in general we cannot conclude that each
energy solution of (\ref{e:op}) converges to a \emph{single} equilibrium, as
the set of steady states associated with (\ref{o3}) can be quite complicated
(see, e.g., \cite{Ha, WXL12}). This means that we are required to prove (\ref%
{o2}) for the whole sequence $\left\{ t_{n}\right\} $ and not only a
subsequence. The main tool is the same energy functional from (\ref{total
energy of the system}) $E_{Q}\left( t\right) =:\mathcal{E}_{Q}\left( u\left(
t\right) ,d\left( t\right) \right) $, that is,%
\begin{equation*}
\mathcal{E}_{Q}\left( u,d\right) :=\frac{1}{2}\left\langle u,Qu\right\rangle
+\widehat{\mathcal{E}}\left( d\right) ,\text{ }\widehat{\mathcal{E}}\left(
d\right) :=\frac{1}{2}||A_{1}^{1/2}d||_{L^{2}}^{2}+\int_{\Omega }W\left(
d\right) dx.
\end{equation*}%
We note that $d_{\ast }$ is a critical point of $\widehat{\mathcal{E}}$ over
$D(A_{1}^{1/2}).$

The version of the \L ojasiewicz-Simon inequality we need is given by the
following lemma, proved in \cite{CJ, JB}.

\begin{lemma}
\label{l3}There exist constants $\zeta \in (0,1/2)$ and $C_{L}>0,$ $\eta >0$
depending on $d_{\ast }$ such that, for any $d\in D(A_{1}^{1/2}),$ if $%
\left\Vert d-d_{\ast }\right\Vert _{1}\leq \eta ,$ denoting by $\widehat{%
\mathcal{E}}^{\prime }$ the Fr\'{e}chet derivative of $\widehat{\mathcal{E}}$%
, we have%
\begin{equation}
C_{L}||\widehat{\mathcal{E}}{^{\prime }}\left( d\right) ||_{-1}\geq |%
\widehat{\mathcal{E}}\left( d\right) -\widehat{\mathcal{E}}\left( d_{\ast
}\right) |^{1-\zeta }.  \label{3.2bis}
\end{equation}
\end{lemma}

The following statement is valid for any energy solution $\left( u,d\right) $
of Definition \ref{weak-energy}.

\begin{proposition}
\label{p:const-energ}There exists a constant $e_{\infty }\in \mathbb{R}$
such that $\widehat{\mathcal{E}}\left( d_{\ast }\right) =e_{\infty },$ for
all solutions $d_{\ast }$ satisfying (\ref{o3}), and we have%
\begin{equation}
\lim_{t\rightarrow \infty }\mathcal{E}_{Q}\left( u\left( t\right) ,d\left(
t\right) \right) =e_{\infty }.  \label{e_const}
\end{equation}%
Moreover, the functional $\Phi \left( t\right) $ is nonincreasing along all
energy solutions $\left( u\left( t\right) ,d\left( t\right) \right) $ and,
for all $t\geq 0,$%
\begin{equation}
\frac{d}{dt}\Phi \left( t\right) \leq -\left( \frac{c_{A_{0}}}{2}\left\Vert
u\left( t\right) \right\Vert _{\theta -\theta _{2}}^{2}-\lambda
_{1}^{-1}\left\Vert A_{1}d\left( s\right) +f\left( d\left( s\right) \right)
\right\Vert _{L^{2}}^{2}\right) ,  \label{big-en}
\end{equation}%
where%
\begin{equation}
\Phi \left( t\right) :=\mathcal{E}_{Q}\left( u\left( t\right) ,d\left(
t\right) \right) +\frac{2\Vert Q\Vert _{-\theta _{2};\theta _{2}}^{2}}{%
c_{A_{0}}}\int_{t}^{\infty }\left\Vert g\left( s\right) \right\Vert
_{-\theta -\theta _{2}}^{2}ds.  \label{big_phi}
\end{equation}
\end{proposition}

Our first result is concerned with the convergence of energy solutions of
problem (\ref{e:op}) to single equilibria, showing that their $\omega $%
-limit set is always a singleton.

\begin{theorem}
\label{conv-equil} Let the assumptions of Lemma \ref{ss:omega} hold. The $%
\omega $-limit set of the component $d$ of any weak energy solution $\left(
u,d\right) $ of problem (\ref{e:op}), as given by Definition \ref%
{weak-energy}, is a singleton. Further we also have%
\begin{equation}
\lim_{t\rightarrow \infty }||d(t)-d_{\ast }||_{L^{2}}=0,\text{ }%
\lim_{t\rightarrow \infty }\left\langle u\left( t\right) ,v\right\rangle =0,%
\text{ }\forall v\in V^{\theta _{2}}  \label{conv-l2}
\end{equation}%
and the following convergence rate:%
\begin{equation}
||d(t)-d_{\ast }||_{W^{1-}}\lesssim \left( 1+t\right) ^{-\chi },
\label{conv_rate}
\end{equation}%
for some $\chi \in \left( 0,1\right) $ depending on $d_{\ast }$.
\end{theorem}

\begin{proof}
The second claim of (\ref{conv-l2}) follows from (\ref{o1}). To prove the
first claim, we adapt the ideas of \cite{CJ, JB} and use an argument that we
applied in \cite[Section 5.3]{GM-JNS} for a simplified regularized
Ericksen-Leslie model. First, we observe that the conclusion of Lemma \ref%
{l3} also holds, provided that we choose even a smaller constant $\zeta \in
(0,\frac{1}{2})\cap $ $(0,\delta \left( 1+\delta \right) ^{-1}),$ where $%
\delta <1$ is the decay rate in (\ref{g-assumpt}). In order to see that, it
suffices to choose a constant $\eta >0$ in Lemma \ref{l3} so small that $|%
\widehat{\mathcal{E}}\left( d\right) -\widehat{\mathcal{E}}\left( d_{\ast
}\right) |\leq 1$ whenever $\left\Vert d-d_{\ast }\right\Vert _{1}\leq \eta $%
. Define further%
\begin{equation*}
\widehat{\Phi }\left( t\right) :=\mathcal{E}_{Q}\left( u\left( t\right)
,d\left( t\right) \right) -\widehat{\mathcal{E}}\left( d_{\ast }\right) +%
\frac{2\Vert Q\Vert _{-\theta _{2};\theta _{2}}^{2}}{c_{A_{0}}}%
\int_{t}^{\infty }\left\Vert g\left( s\right) \right\Vert _{-\theta -\theta
_{2}}^{2}ds
\end{equation*}%
and notice that $\widehat{\Phi }\left( t\right) $ differs from $\Phi \left(
t\right) $ in (\ref{big_phi}) only by a constant. Hence, setting%
\begin{equation*}
\Upsilon \left( t\right) :=\left\Vert u\left( t\right) \right\Vert _{\theta
-\theta _{2}}+\left\Vert A_{1}d\left( t\right) +f\left( d\left( t\right)
\right) \right\Vert _{L^{2}},
\end{equation*}%
for every $t\geq 0$, from (\ref{big-en}) we have%
\begin{equation}
\frac{d}{dt}\widehat{\Phi }\left( t\right) \lesssim -\Upsilon ^{2}\left(
t\right) \leq 0  \label{cs9}
\end{equation}%
so that $\widehat{\Phi }$ is also a nonincreasing function on $\left[
0,\infty \right) $. Furthermore, integrating this relation over $\left(
0,\infty \right) $ and recalling (\ref{e_const}), we also obtain%
\begin{equation*}
\int_{0}^{\infty }\Upsilon ^{2}\left( s\right) ds<\infty ,
\end{equation*}%
thanks to (\ref{o1})-(\ref{o3}). Together with the energy inequality (\ref%
{en-both}), then one has%
\begin{align}
u& \in L^{\infty }\left( 0,\infty ;V^{-\theta _{2}}\right) \cap L^{2}\left(
0,\infty ;V^{\theta -\theta _{2}}\right) ,  \label{cs10} \\
d& \in L^{\infty }\left( 0,\infty ;W^{1}\right) ,  \label{cs10bis} \\
A_{1}d+f\left( d\right) & \in L^{2}\left( 0,\infty ;L^{2}\left( \Omega
\right) \right) .  \label{cs10tris}
\end{align}%
As before, these bounds together with proper handling of the other nonlinear
terms in the director equation of (\ref{e:op}) imply%
\begin{align}
A_{Q}d& \in L^{2}\left( 0,\infty ;W^{-2}\right) ,  \label{cs11} \\
\omega _{Q}d& \in L^{2}\left( 0,\infty ;W^{-2}\right) ,  \label{cs11bis} \\
B_{1}\left( u,d\right) & \in L^{2}\left( 0,\infty ;W^{-2}\right) .
\label{cs11tris}
\end{align}%
In detail the estimates in (\ref{cs11})-(\ref{cs11tris}) are obtained by
application of Lemma \ref{l:hole} (Appendix), as follows:

\begin{itemize}
\item[(1)] The terms $\left( \partial _{i}Qu_{j}\right) d_{j}$ are a product
of functions in $H^{\theta +\theta _{2}-1}$ and $H^{1}$ and therefore
bounded in $H^{-2}$ since $\theta ,\theta _{2}\geq 0$.

\item[(2)] The terms $\left( Qu_{i}\right) \partial _{i}d_{j}$ are a product
of functions in $H^{\theta +\theta _{2}}$ and $L^{2},$ and therefore bounded
in $H^{-2}$ again, since $\theta ,\theta _{2}\geq 0.$
\end{itemize}

Therefore, once again comparing terms in the director equation for $d$ from (%
\ref{e:op}), the estimates (\ref{cs10tris})-(\ref{cs11tris}) entail%
\begin{equation}
\partial _{t}d\in L^{2}\left( 0,\infty ;W^{-2}\right) .  \label{cs12}
\end{equation}

Our next goal is to show that there exists $t_{0}>0$ sufficiently large,
such that $\partial _{t}d\in L^{1}\left( t_{0},\infty ;W^{-2}\right) $. Now,
define%
\begin{equation*}
\Sigma :=\{t\geq 1:||d\left( t\right) -d_{\ast }||_{L^{2}}\leq \eta /3\}
\end{equation*}%
and observe that $\Sigma $ is unbounded by Lemma \ref{ss:omega}. For every $%
t\in \Sigma ,$ we define%
\begin{equation*}
\tau \left( t\right) =\sup \{t^{^{\prime }}\geq t:\sup_{s\in \lbrack
t,t^{^{\prime }}]}||d\left( t\right) -d_{\ast }||_{L^{2}}<\eta \}.
\end{equation*}%
By continuity, $\tau \left( t\right) >t$ for every $t\in \Sigma $. Let now $%
t_{0}\in \Sigma $ and divide the interval $J:=[t_{0},\tau \left(
t_{0}\right) )$ into two subsets%
\begin{equation*}
\Sigma _{1}:=\left\{ t\in J:\Upsilon \left( t\right) \geq \left(
\int_{t}^{\tau \left( t_{0}\right) }\left\Vert g\left( s\right) \right\Vert
_{-\theta -\theta _{2}}^{2}ds\right) ^{1-\zeta }\right\} ,\text{ }\Sigma
_{2}:=J\backslash \Sigma _{1}.
\end{equation*}

Setting further, as above,%
\begin{equation*}
\widehat{\Phi }\left( t\right) :=\mathcal{E}_{Q}\left( u\left( t\right)
,d\left( t\right) \right) -\widehat{\mathcal{E}}\left( d_{\ast }\right) +%
\frac{2\Vert Q\Vert _{-\theta _{2};\theta _{2}}^{2}}{c_{A_{0}}}%
\int_{t}^{\tau \left( t_{0}\right) }\left\Vert g\left( s\right) \right\Vert
_{-\theta -\theta _{2}}^{2}ds
\end{equation*}%
we notice that $\widehat{\Phi }\left( t\right) $ again satisfies (\ref{cs9})
for every $t\in J,$ and $\widehat{\Phi }$ is a nonincreasing function on $J$%
. Moreover, for every $t\in J$ we have%
\begin{align}
\frac{d}{dt}\left( |\widehat{\Phi }\left( t\right) |^{\zeta }sgn(\widehat{%
\Phi }\left( t\right) \right) & =\zeta |\widehat{\Phi }\left( t\right)
|^{\zeta -1}\frac{d}{dt}\widehat{\Phi }\left( t\right)  \label{rel4.4} \\
& \lesssim -|\widehat{\Phi }\left( t\right) |^{\zeta -1}\Upsilon ^{2}\left(
t\right) ,  \notag
\end{align}%
which implies that the functional $sgn(\widehat{\Phi }\left( t\right) )|%
\widehat{\Phi }\left( t\right) |^{\zeta }$ is decreasing on $J$. By (\ref%
{3.2bis}) and Proposition \ref{p:const-energ}, for every $t\in \Sigma _{1}$
we can easily establish%
\begin{align}
|\widehat{\Phi }\left( t\right) |^{1-\zeta }& \leq \left\vert \mathcal{E}%
_{Q}\left( u\left( t\right) ,d\left( t\right) \right) -\widehat{\mathcal{E}}%
\left( d_{\ast }\right) \right\vert ^{1-\zeta }+\left( \frac{2\Vert Q\Vert
_{-\theta _{2};\theta _{2}}^{2}}{c_{A_{0}}}\int_{t}^{\tau \left(
t_{0}\right) }\left\Vert g\left( s\right) \right\Vert _{-\theta -\theta
_{2}}^{2}ds\right) ^{1-\zeta }  \label{rel4.5} \\
& \lesssim \Upsilon \left( t\right) ,  \notag
\end{align}%
owing to the basic inequality
\begin{equation*}
\left\Vert u\left( t\right) \right\Vert _{-\theta _{2}}^{2}\lesssim
\left\Vert u\left( t\right) \right\Vert _{\theta -\theta _{2}}^{\frac{1}{%
1-\zeta }}\text{, for a.e. }t>0,
\end{equation*}%
which holds thanks to the embedding $V^{\theta -\theta _{2}}\subseteq
V^{-\theta _{2}}$ and (\ref{cs10}). Combining now (\ref{rel4.5}) with (\ref%
{rel4.4}) yields%
\begin{equation}
-\frac{d}{dt}\left( |\widehat{\Phi }\left( t\right) |^{\zeta }sgn(\widehat{%
\Phi }\left( t\right) \right) \gtrsim \Upsilon \left( t\right) .
\label{rel4.6}
\end{equation}%
Moreover, exploiting (\ref{rel4.6}) we have%
\begin{align}
\int_{\Sigma _{1}}\Upsilon \left( s\right) ds& \lesssim -\int_{\Sigma _{1}}%
\frac{d}{ds}\left( |\widehat{\Phi }\left( s\right) |^{\zeta }sgn(\widehat{%
\Phi }\left( s\right) \right) ds  \label{rel4.7} \\
& \lesssim \left( |\widehat{\Phi }\left( t_{0}\right) |^{\zeta }+|\widehat{%
\Phi }\left( \tau \left( t_{0}\right) \right) |^{\zeta }\right) <\infty ,
\notag
\end{align}%
where we interpret the term involving $\tau \left( t_{0}\right) $ on the
right hand side of (\ref{rel4.7}) as $0$ if $\tau \left( t_{0}\right)
=\infty $ (recall (\ref{e_const})). On the other hand, if $t\in \Sigma _{2},$
using assumption (\ref{g-assumpt}) we obtain%
\begin{equation}
\Upsilon \left( t\right) \leq \left( \int_{t}^{\tau \left( t_{0}\right)
}\left\Vert g\left( s\right) \right\Vert _{-\theta -\theta
_{2}}^{2}ds\right) ^{1-\zeta }\lesssim \left( 1+t\right) ^{-\left( 1-\zeta
\right) \left( 1+\delta \right) },  \label{rel4.8}
\end{equation}%
so once again the function $\Upsilon $ is dominated by an integrable
function on $\Sigma _{2}$ since $\zeta \left( 1+\delta \right) <\delta $.
Combining the inequalities (\ref{rel4.7}), (\ref{rel4.8}), we deduce that $%
\Upsilon $ is absolutely integrable on $J$ and%
\begin{equation}
\lim_{t_{0}\rightarrow \infty ,t_{0}\in \Sigma }\int_{t_{0}}^{\tau \left(
t_{0}\right) }\Upsilon \left( s\right) ds=0.  \label{rel4.9}
\end{equation}%
On the other hand, recalling estimates (\ref{cs10})-(\ref{cs11tris}) and the
observations (1)-(2), from the second equation of (\ref{e:op}) it follows
that%
\begin{align}
\left\Vert \partial _{t}d\left( t\right) \right\Vert _{-2}& \lesssim
\left\Vert B_{1}\left( u\left( t\right) ,d\left( t\right) \right)
\right\Vert _{L^{2}}+\left\Vert A_{1}d\left( t\right) +f\left( d\left(
t\right) \right) \right\Vert _{L^{2}}  \label{rel4.10} \\
& +\left\Vert A_{Q}d\left( t\right) \right\Vert _{L^{2}}+\left\Vert \omega
_{Q}d\left( t\right) \right\Vert _{L^{2}}  \notag \\
& \lesssim \left\Vert u\left( t\right) \right\Vert _{\theta -\theta
_{2}}\left\Vert d\left( t\right) \right\Vert _{1}+\left\Vert A_{1}\phi
\left( t\right) +f\left( \phi \left( t\right) \right) \right\Vert _{L^{2}}
\notag \\
& \lesssim \Upsilon \left( t\right) .  \notag
\end{align}%
Consequently, we also have%
\begin{equation}
\lim_{t_{0}\rightarrow \infty ,t_{0}\in \Sigma }\int_{t_{0}}^{\tau \left(
t_{0}\right) }\left\Vert \partial _{t}d\left( s\right) \right\Vert _{-2}ds=0.
\label{rel4.11}
\end{equation}%
This fact, combined with a simple contradiction argument (see \cite{CJ},
\cite{GM-JNS}) yields that we must have $\tau \left( t_{0}\right) =\infty ,$
for some sufficiently large $t_{0}\in \Sigma $. Thus, $\partial _{t}d\in
L^{1}\left( t_{0},\infty ;W^{-2}\right) $ as desired. By compactness and a
basic interpolation inequality, we have $d\left( t\right) \rightarrow
d_{\ast }$ in the strong topology of $W^{1-}$. Hence, the $\omega $-limit
set of the $d$ component of any weak energy solution $\left( u,d\right) $ is
the singleton $d_{\ast }$, which is a solution of (\ref{o3}). The estimate
of the rate of convergence in (\ref{conv_rate}) is a straightforward
consequence of (\ref{rel4.4})-(\ref{rel4.5}), the definition of $\Phi $ and
basic interpolation results. We leave the details to the interested reader.
The proof of Theorem \ref{conv-equil} is complete.
\end{proof}

\begin{remark}
Our theorem covers all the special cases listed in Table \ref{t:spec2} and
many other models (see Remark \ref{rem:all-mod}). In particular, our result
yields convergence to a single steady state $\left( 0,d_{\ast }\right) $ of
any weak \emph{energy solution} of the full three dimensional NSE-EL,
Leray-EL-$\alpha $, ML-EL-$\alpha $, NSV-EL, SBM-EL, NS-EL-$\alpha $ models.
None of these results have been reported previously.
\end{remark}

We can derive a sufficient condition such that a stronger convergence result
holds in (\ref{conv-l2}).

\begin{theorem}
\label{t:strong-conv}Let $\left( u,d\right) \in L^{\infty }\left( 0,\infty ;%
\mathcal{Y}_{\theta _{2}}\right) $ be an energy solution in the sense of
Lemma \ref{ss:omega}, determined by the assumptions of Theorem \ref{t:exist}%
. In addition, assume%
\begin{equation}
\theta +\theta _{2}\geq 1\text{ and }d\in L^{\infty }\left( 0,\infty
;L^{\infty }\left( \Omega \right) \right) .  \label{criteria-conv}
\end{equation}%
Then, there holds%
\begin{equation}
\lim_{t\rightarrow +\infty }\left( \left\Vert u\left( t\right) \right\Vert
_{-\theta _{2}}+||A_{1}^{1/2}\left( d\left( t\right) -d_{\ast }\right)
||_{L^{2}}\right) =0,  \label{conv-st}
\end{equation}%
where $d_{\ast }\in D\left( A_{1}\right) $ is a solution of (\ref{o3}).
\end{theorem}

\begin{proof}
We recall that each energy solution $\left( u,d\right) $\ of Definition \ref%
{weak-energy} satisfies the bounds (\ref{cs10})-(\ref{cs10tris}) and that $%
y\in L^{1}\left( 0,\infty \right) ,$ owing to $V^{\theta -\theta
_{2}}\subseteq V^{-\theta _{2}}$, where we have set%
\begin{equation*}
y\left( t\right) :=\frac{1}{2}\left\langle u,Qu\right\rangle +\left\Vert
\rho \left( t\right) \right\Vert _{-1}^{2},\text{ }\rho \left( t\right)
:=A_{1}d\left( t\right) +f\left( d\left( t\right) \right) .
\end{equation*}%
In particular, by (\ref{en-case1})-(\ref{apriori-tensor}) we recall that
\begin{equation}
\int_{0}^{\infty }\left( \left\Vert u\left( s\right) \right\Vert _{\theta
-\theta _{2}}^{2}+\left\Vert \rho \left( s\right) \right\Vert
_{L^{2}}^{2}+\mu _{1}\left\Vert d^{T}A_{Q}d\right\Vert _{L^{2}}^{2}\right)
ds<\infty .  \label{cs16v0}
\end{equation}%
By the second condition in (\ref{criteria-conv}), there exists a constant $%
M>0$ independent of time such that
\begin{equation*}
\left\Vert d\right\Vert _{L^{\infty }\left( \mathbb{R}_{+};L^{\infty }\left(
\Omega \right) \right) }\leq M.
\end{equation*}%
Our goal is to show that $y$ satisfies the inequality%
\begin{equation}
\frac{dy}{dt}\left( t\right) \leq C+\Lambda \left( t\right) ,\text{ for all }%
t\geq 0,  \label{crucial-ineq}
\end{equation}%
for some constant $C>0$ independent of time and some function $\Lambda \in
L^{1}\left( 0,\infty \right) $. Then, the application of \cite[Lemma 6.2.1]%
{Zh} yields $y\left( t\right) \rightarrow 0$ as $t\rightarrow \infty $; the
convergence (\ref{conv-st}) is then an immediate consequence of this crucial
fact, owing to the basic inequality%
\begin{equation*}
||A_{1}^{1/2}\left( d\left( t\right) -d_{\ast }\right) ||_{L^{2}}\lesssim
||\rho \left( t\right) ||_{-1}+\left\Vert f\left( d\right) -f\left( d_{\ast
}\right) \right\Vert _{-1}+\left\Vert A_{1}d_{\ast }+f\left( d_{\ast
}\right) \right\Vert _{-1}
\end{equation*}%
and (\ref{conv-l2}), (\ref{o3}). Of course, (\ref{crucial-ineq}) can be
justified by employing a proper approximation scheme, such as the one used
in the proof of Theorem \ref{t:exist}.

To this end, we first pair the first equation of (\ref{e:op}) with $Qu$,
then use assumption (ii) of Theorem \ref{t:exist}, to deduce the identity%
\begin{equation}
\frac{1}{2}\frac{d}{dt}\left\langle u,Qu\right\rangle +\left\langle
A_{0}u,Qu\right\rangle =\left\langle g,Qu\right\rangle +\left\langle
R_{0}\left( \rho ,d\right) ,Qu\right\rangle +\left\langle \sigma _{Q},\nabla
\left( Qu\right) \right\rangle ,  \label{cs15}
\end{equation}%
where the tensor $\sigma _{Q}$ is given by (\ref{v5}). On the other hand, in
view of the second equation of (\ref{e:op}) we have%
\begin{align}
\frac{d}{dt}\left\Vert \rho \right\Vert _{-1}^{2}& =-\left\langle
B_{1}\left( u,d\right) ,\rho \right\rangle +\left\langle \omega _{Q}d,\rho
\right\rangle -\frac{\lambda _{2}}{\lambda _{1}}\left\langle A_{Q}d,\rho
\right\rangle +\frac{1}{\lambda _{1}}\left\Vert \rho \right\Vert _{L^{2}}^{2}
\label{cs16} \\
& -\left\langle f^{^{\prime }}\left( d\right) B_{1}\left( u,d\right)
,A_{1}^{-1}\rho \right\rangle +\left\langle f^{^{\prime }}\left( d\right)
\omega _{Q}d,A_{1}^{-1}\rho \right\rangle  \notag \\
& -\frac{\lambda _{2}}{\lambda _{1}}\left\langle f^{^{\prime }}\left(
d\right) A_{Q}d,A_{1}^{-1}\rho \right\rangle +\frac{1}{\lambda _{1}}%
\left\langle f^{^{\prime }}\left( d\right) \rho ,A_{1}^{-1}\rho
\right\rangle .  \notag
\end{align}%
Adding the relations (\ref{cs15})-(\ref{cs16}) together, then using (\ref%
{e:coercive-abis}) and noting that $\left\langle B_{1}\left( u,d\right)
,\rho \right\rangle =\left\langle R_{0}\left( \rho ,d\right)
,Qu\right\rangle $, one has%
\begin{align}
& \frac{dy}{dt}+c_{A_{0}}\left\Vert u\right\Vert _{\theta -\theta _{2}}^{2}-%
\frac{1}{\lambda _{1}}\left\Vert \rho \right\Vert _{L^{2}}^{2}  \label{cs17}
\\
& \leq \left\langle g,Qu\right\rangle +\left\langle \sigma _{Q},\nabla
\left( Qu\right) \right\rangle +\left\langle \omega _{Q}d,\rho \right\rangle
-\frac{\lambda _{2}}{\lambda _{1}}\left\langle A_{Q}d,\rho \right\rangle
\notag \\
& -\left\langle f^{^{\prime }}\left( d\right) B_{1}\left( u,d\right)
,A_{1}^{-1}\rho \right\rangle +\left\langle f^{^{\prime }}\left( d\right)
\omega _{Q}d,A_{1}^{-1}\rho \right\rangle  \notag \\
& -\frac{\lambda _{2}}{\lambda _{1}}\left\langle f^{^{\prime }}\left(
d\right) A_{Q}d,A_{1}^{-1}\rho \right\rangle +\frac{1}{\lambda _{1}}%
\left\langle f^{^{\prime }}\left( d\right) \rho ,A_{1}^{-1}\rho \right\rangle
\notag \\
& =:I_{1}+...+I_{8}.  \notag
\end{align}%
We now obtain proper bounds for the terms on the right-hand side in the
following manner:

\begin{itemize}
\item[(b1)] As usual for the first one, for every $\delta >0$ we have%
\begin{equation*}
\left\vert I_{1}\right\vert \leq \delta \left\Vert u\right\Vert _{\theta
-\theta _{2}}^{2}+C_{\delta }\left\Vert g\right\Vert _{-\theta -\theta
_{2}}^{2}.
\end{equation*}

\item[(b2)] For $I_{3},I_{4}$, exploiting (\ref{criteria-conv}) together
with the Sobolev embedding $V^{\theta -\theta _{2}}\subseteq V^{1-2\theta
_{2}}$, as $\theta +\theta _{2}\geq 1,$ yields%
\begin{equation*}
\left\vert I_{3}\right\vert +\left\vert I_{4}\right\vert \leq \left\Vert
\rho \right\Vert _{L^{2}}\left\Vert \nabla \left( Qu\right) \right\Vert
_{L^{2}}\left\Vert d\right\Vert _{L^{\infty }}\leq \delta \left\Vert
u\right\Vert _{\theta -\theta _{2}}^{2}+C_{\delta ,M}\left\Vert \rho
\right\Vert _{L^{2}}^{2}.
\end{equation*}

\item[(b3)] Using the definition for $\sigma _{Q}$, we further split the
second term $I_{2}$ into three more terms $I_{21},I_{22},I_{23}$. The first
one we bound as follows:%
\begin{align*}
\left\vert I_{21}\right\vert & =\left\langle \mu _{1}(d^{T}A_{Q}d)d\otimes
d,\nabla Qu\right\rangle \\
& \leq \mu _{1}\left\Vert \nabla Qu\right\Vert _{L^{2}}\left\Vert
d^{T}A_{Q}d\right\Vert _{L^{2}}\left\Vert d\otimes d\right\Vert _{L^{\infty
}} \\
& \lesssim \mu _{1}\left\Vert u\right\Vert _{\theta -\theta _{2}}\left\Vert
d^{T}A_{Q}d\right\Vert _{L^{2}}\left\Vert d\right\Vert _{L^{\infty }}^{2} \\
& \leq \delta \left\Vert u\right\Vert _{\theta -\theta _{2}}^{2}+C_{\delta
,M}\left( \mu _{1}\left\Vert d^{T}A_{Q}d\right\Vert _{L^{2}}\right) ^{2}.
\end{align*}%
Next, since by definition $\mathcal{N}_{Q}=\lambda _{1}^{-1}\rho -\lambda
_{2}/\lambda _{1}A_{Q}d$, proceeding as for the estimate for $I_{3}$, we get%
\begin{align*}
\left\vert I_{22}\right\vert & =\left\vert \left\langle \mu _{2}\mathcal{N}%
_{Q}\otimes d+\mu _{3}d\otimes \mathcal{N}_{Q},\nabla Qu\right\rangle
\right\vert \\
& \lesssim \left\Vert \nabla Qu\right\Vert _{L^{2}}\left\Vert \mathcal{N}%
_{Q}\right\Vert _{L^{2}}\left\Vert d\right\Vert _{L^{\infty }} \\
& \leq C_{M}\left\Vert u\right\Vert _{\theta -\theta _{2}}\left( \left\Vert
\rho \right\Vert _{L^{2}}+\left\Vert A_{Q}d\right\Vert _{L^{2}}\right) \\
& \leq C_{M}\left\Vert u\right\Vert _{\theta -\theta _{2}}\left\Vert \rho
\right\Vert _{L^{2}}+C_{M}\left\Vert u\right\Vert _{\theta -\theta _{2}}^{2},
\end{align*}%
owing once more to the boundedness of $d$. Finally, a similar argument gives
the same bound:%
\begin{align*}
\left\vert I_{23}\right\vert & =\left\vert \left\langle \mu
_{5}(A_{Q}d)\otimes d+\mu _{6}d\otimes (A_{Q}d),\nabla Qu\right\rangle
\right\vert \\
& \leq C_{M}\left\Vert u\right\Vert _{\theta -\theta _{2}}\left( \left\Vert
\rho \right\Vert _{L^{2}}+1\right) +C_{M}\left\Vert u\right\Vert _{\theta
-\theta _{2}}^{2}.
\end{align*}

\item[(b4)] Since the term $f^{^{\prime }}\left( d\right) \left( B_{1}\left(
u,d\right) \right) $ is bounded in $W^{-1}$ as a product of vector-valued
functions in $W^{1}$ and $L^{2}$, using $V^{\theta -\theta _{2}}\subseteq
V^{-\theta _{2}}$ we derive%
\begin{align*}
\left\vert I_{3}\right\vert & =\left\vert \left\langle f^{^{\prime }}\left(
d\right) B_{1}\left( u,d\right) ,A_{1}^{-1}\rho \right\rangle \right\vert \\
& \lesssim \left\Vert f^{^{\prime }}\left( d\right) \right\Vert
_{1}\left\Vert B_{1}\left( u,d\right) \right\Vert _{L^{2}}\left\Vert \rho
\right\Vert _{-1} \\
& \lesssim \left\Vert f^{^{\prime }}\left( d\right) \right\Vert
_{1}\left\Vert \nabla Qu\right\Vert _{1}\left\Vert A_{1}d\right\Vert
_{L^{2}}\left\Vert \rho \right\Vert _{-1} \\
& \leq \delta \left\Vert u\right\Vert _{\theta -\theta _{2}}^{2}+C_{\delta
}\left( \left\Vert \rho \right\Vert _{L^{2}}+\left\Vert f\left( d\right)
\right\Vert _{L^{2}}\right) ^{2}\left\Vert f^{^{\prime }}\left( d\right)
\right\Vert _{1}^{2}\left\Vert \rho \right\Vert _{-1}^{2}.
\end{align*}

\item[(b5)] Deriving a bound for $I_{4}$ and $I_{5}$, one argues exactly in
the same fashion using the definition of the tensors $A_{Q}$ and $\omega
_{Q} $. One has%
\begin{align*}
\left\vert I_{4}\right\vert & =\left\vert \left\langle f^{^{\prime }}\left(
d\right) \omega _{Q}d,A_{1}^{-1}\rho \right\rangle \right\vert \\
& \lesssim \left\Vert f^{^{\prime }}\left( d\right) \right\Vert
_{1}\left\Vert \omega _{Q}d\right\Vert _{L^{2}}\left\Vert \rho \right\Vert
_{-1} \\
& \lesssim \left\Vert f^{^{\prime }}\left( d\right) \right\Vert
_{1}\left\Vert \nabla Qu\right\Vert _{L^{2}}\left\Vert d\right\Vert
_{L^{\infty }}\left\Vert \rho \right\Vert _{-1} \\
& \leq \delta \left\Vert u\right\Vert _{\theta -\theta _{2}}^{2}+C_{\delta
,M}\left\Vert f^{^{\prime }}\left( d\right) \right\Vert _{1}^{2}\left\Vert
\rho \right\Vert _{-1}^{2}.
\end{align*}%
We also get the same upper bound for $I_{5}.$

\item[(b6)] For the final term $I_{2}$, we use the fact that $f^{^{\prime
}}\left( d\right) \rho $ is bounded in $W^{-1}$ as a product of
vector-valued functions in $W^{1}\times L^{2}$. One readily obtains the
estimate:%
\begin{equation*}
\left\vert I_{8}\right\vert =\frac{1}{\lambda _{1}}\left\vert \left\langle
f^{^{\prime }}\left( d\right) \rho ,A_{1}^{-1}\rho \right\rangle \right\vert
\leq \delta \left\Vert \rho \right\Vert _{L^{2}}^{2}+C_{\delta }\left\Vert
f^{^{\prime }}\left( d\right) \right\Vert _{1}^{2}\left\Vert \rho
\right\Vert _{-1}^{2}.
\end{equation*}
\end{itemize}

Collect now all the previous estimates from (b1)-(b6) and insert them on the
right-hand side of (\ref{cs17}). Choosing a sufficiently small $\delta \sim
\min \left( c_{A_{0}},-\lambda _{1}^{-1}\right) >0$, taking into account the
uniform in-time bounds $\left\Vert \rho \left( t\right) \right\Vert
_{-1}\leq C,$ $||f^{^{\prime }}\left( d\left( t\right) \right) ||_{1}\leq C$
and $\left\Vert \left( u\left( t\right) ,d\left( t\right) \right)
\right\Vert _{\mathcal{Y}_{\theta _{2}}}\leq C$, which hold for all times $%
t\geq 0$, and the basic controls (\ref{cs16v0}), (\ref{g-assumpt}), we
readily infer the validity of inequality (\ref{crucial-ineq}) for a proper
function $\Lambda \in L^{1}\left( 0,\infty \right) $. Hence, we have proved
our claim and the proof is finished.
\end{proof}

\begin{remark}
The second assumption of (\ref{criteria-conv}) can be slightly weakened if $%
\mu _{1}=0$ ($\lambda _{2}\neq 0$). We recall that when $\lambda _{2}=0,$ $%
\mu _{1}\geq 0$, the second of (\ref{criteria-conv}) is already satisfied on
account of Proposition \ref{maxp}.
\end{remark}

For example, setting $\theta =1$, $\theta _{1}=\theta _{2}=0$ and $\mathcal{V%
}=\mathcal{V}_{per},$ $\Omega =\mathbb{T}^{2}$, global existence of a unique
strong solution in the class%
\begin{equation}
\left( u,d\right) \in L^{\infty }\left( 0,\infty ;V^{1}\times D\left(
A_{1}\right) \right)  \label{extra}
\end{equation}%
for the two dimensional 2D Navier-Stokes-Ericksen-Leslie model was
established in \cite[Lemma 4.3 and Remark 4]{WXL12} in either of the
following cases (a) $\mu _{1}=0,$ $\lambda _{2}\neq 0$, (b) $\mu _{1}\geq 0,$
$\lambda _{2}=0$. We observe that due to the additional regularity (\ref%
{extra}) and the embedding $D\left( A_{1}\right) \subset L^{\infty }$, the
second of (\ref{criteria-conv}) is automatically satisfied for this model
and so the conclusion of Theorem \ref{t:strong-conv} holds. In this case,
the convergence result can also be found in \cite[Theorem 4.2]{WXL12}. On
the other hand, the convergence to a single steady state $\left( 0,d_{\ast
}\right) $ of any (regular) energy solution, satisfying the assumptions of
Theorem \ref{t:reg} is also ensured on account of Remark \ref{rem-bd} and
Proposition \ref{t:attr-exist2}. In particular, this is true for the 3D
modified Leray-EL-$\alpha $ (ML-EL-$\alpha $) model, the 3D SBM-EL model and
the 3D NS-EL-$\alpha $ system in the general case when $\mu _{1}\geq 0$ and $%
\lambda _{2}\neq 0$. Besides, in the case (b) above, the 3D NSV-EL model
possesses (unique)\ energy solutions that converge to single steady states
as concluded by (\ref{conv-st}). These convergence results were not
previously reported in the literature for any of these models.

\section{Concluding remarks}

\label{s:cr}

In this article, we consider a general family of regularized Ericksen-Leslie
models which captures some specifics and variants of the models that have
not been considered or analyzed anywhere in the literature before. We give a
unified analysis of the Ericksen-Leslie system using tools in nonlinear
analysis and Sobolev function theory together with energy methods, and then
use them to obtain sharp results. In particular, in Section \ref{s:well} we
develop some well-posedness results for our family of nonlinear models,
which include existence results (Section \ref{ss:exist}), regularity results
(Section \ref{ss:reg}), and uniqueness and stability results (Section \ref%
{ss:stab}). In Section \ref{ss:global}, we show the existence of a
finite-dimensional global attractor in the case $\theta >0$ and give some
further properties, by first establishing the existence of an exponential
attractor. In Section \ref{s:convss}, we prove the asymptotic stabilization
as time goes to infinity of any energy solution for our problem (\ref{e:op})
to a single steady state. The present unified analysis can be exploited to
extend and establish existence, regularity and existence of finite
dimensional attractor results also in the case $\theta =0$; this case is
more delicate and requires a more refined analysis which lies beyond the
scope of the present article. Indeed, problem (\ref{e:op}) with $\theta =0$
can be seen as a non-dissipative system in which the fluid equation looses
its parabolic character and behaves more like a hyperbolic equation. For
instance, this is the case when the velocity component satisfies the 3D
Navier-Stokes-Voigt equation. For a simplified regularized Ericksen-Leslie
model ($\sigma _{Q}\equiv 0,$ $\omega _{Q}\equiv 0$ and $\lambda _{2}=0$),
such results have already been established in \cite{GM-JNS}. For the full
regularized Ericksen-Leslie model (\ref{e:op}), we will consider such
questions in a forthcoming contribution.

We conclude this section with some remarks on the assumption about the space
$\mathcal{V}$ in (\ref{select-p}), and the precise connections between the
models as introduced in Table \ref{t:spec} and their equivalent formulations
which are most recognizable in the physics literature. To this end, let us
assume that the fluid velocity satisfies either one of the following
regularized versions of the 3D Navier-Stokes equations ($M=I$, $Q=I$, $%
A_{0}=-\mu _{4}\Delta $):

\begin{itemize}
\item[(1)] The 3D Leray-$\alpha $ system with $\theta =1,$ $\theta _{1}=1,$ $%
\theta _{2}=0$:%
\begin{equation}
\partial _{t}u-\mu _{4}\Delta u+\left( \Pi u\cdot \nabla \right) u+\nabla p=%
\overrightarrow{F}\left( d\right) ,  \label{f1}
\end{equation}%
with $\Pi =\left( I-\alpha ^{2}\Delta \right) ^{-1}$ and $Q=I.$ In this
case, $u=Q^{-1}v\left( =v\right) $ is the fluid velocity as in (\ref{e:rel}).

\item[(2)] The modified 3D Leray-$\alpha $ (ML-$\alpha $)\ system with $%
\theta =1,$ $\theta _{1}=0,$ $\theta _{2}=1$:%
\begin{equation}
\partial _{t}\left( v-\alpha ^{2}\Delta v\right) -\mu _{4}\Delta \left(
v-\alpha ^{2}\Delta v\right) +\left( v-\alpha ^{2}\Delta v\right) \cdot
\nabla v+\nabla p=\overrightarrow{F}\left( d\right) ,  \label{f2}
\end{equation}%
with $Q=\Pi $ and $M=I$. In this case, $v=\Pi u$ is the (regularized) fluid
velocity as in (\ref{e:rel}).

\item[(3)] The 3D simplified Bardina model with $\theta =\theta _{1}=\theta
_{2}=1$:%
\begin{equation}
\partial _{t}\left( v-\alpha ^{2}\Delta v\right) -\mu _{4}\Delta \left(
v-\alpha ^{2}\Delta v\right) +v\cdot \nabla v+\nabla p=\overrightarrow{F}%
\left( d\right) ,  \label{f3}
\end{equation}%
with $Q=\Pi $ and $M=\Pi $. As above in (3), $v=\Pi u$ is the (regularized)
fluid velocity for the system (\ref{e:rel}).

\item[(4)] The 3D Navier-Stokes-Voigt equation with $\theta =0,$ $\theta
_{1}=\theta _{2}=1$:%
\begin{equation}
\partial _{t}\left( v-\alpha ^{2}\Delta v\right) -\mu _{4}\Delta v+v\cdot
\nabla v+\nabla p=\overrightarrow{F}\left( d\right) ,  \label{f4}
\end{equation}%
with $Q=M=\Pi $. Again $v=\Pi u$ (i.e., $u=\Pi ^{-1}v$) corresponds to the
(regularized) fluid velocity.

\item[(5)] The 3D Lagrangian averaged Navier-Stokes-$\alpha $ equation with $%
\theta =1,$ $\theta _{1}=0,$ $\theta _{2}=1$:%
\begin{equation}
\partial _{t}\left( v-\alpha ^{2}\Delta v\right) -\mu _{4}\Delta \left(
v-\alpha ^{2}\Delta v\right) +v\cdot \nabla \left( v-\alpha ^{2}\Delta
v\right) +\nabla \left( v^{T}\right) \cdot \left( v-\alpha ^{2}\Delta
v\right) +\nabla p=\overrightarrow{F}\left( d\right) ,  \label{f5}
\end{equation}%
with $M=\Pi $ and $Q=I$. In this case $v=Q^{-1}u=u$ corresponds to the
(regularized) fluid velocity.
\end{itemize}

\noindent Above in (1)-(5), the flow is incompressible (i.e., div$\left(
v\right) =0$), $p$ denotes pressure and $\overrightarrow{F}$ consists of a
body force $g\left( t\right) $ acting on the fluid as well as the
stresses/forces due to the coupling of the fluid velocity with the director
field $d$, i.e.,
\begin{equation*}
\overrightarrow{F}\left( d\right) \overset{\text{def}}{=}A_{1}d\cdot \nabla
d+\text{div}\left( \boldsymbol{\sigma }_{Q}\right) +g,
\end{equation*}%
such that $d$ obeys the equation%
\begin{equation}
\partial _{t}d+v\cdot \nabla d-\omega _{Q}d+\frac{\lambda _{2}}{\lambda _{1}}%
A_{Q}d=\frac{1}{\lambda _{1}}\left( A_{1}d+\nabla _{d}W(d)\right) .
\label{df}
\end{equation}%
Here, $A_{Q},\omega _{Q},\sigma _{Q}$ are given in (\ref{v4}) and (\ref{v5}%
), respectively. In particular, we have $A_{Q}=(\nabla v+\nabla ^{T}v)/2$
and $\omega _{Q}=(\nabla v-\nabla ^{T}v)/2.$

We emphasize that problem (\ref{e:rel})-(\ref{example}), for any of the
choices of the parameters $\left( \theta ,\theta _{1},\theta _{2}\right) $
in (1)-(5) above, is in fact equivalent to any regularized 3D
Ericksen-Leslie model (oREL) in which the fluid velocity $v$ satisfies
either one of the equations (\ref{f1})-(\ref{f5}) above and the director
field $d$ satisfies (\ref{df}). Indeed this is the case when the operators $%
Q:V^{-\theta _{2}}\rightarrow V^{\theta _{2}}$ and $Q^{-1}:V^{\theta
_{2}}\rightarrow V^{-\theta _{2}}$ are isometries. Furthermore, according to
the statements proven in Section \ref{s:well}, the transformed problem (\ref%
{e:rel})-(\ref{example}) is well posed in $V^{\beta }\times W^{l}$ for some $%
\beta \geq -\theta _{2},$ $l\geq 1$; this makes any of the (oREL) problems
for (\ref{f1})-(\ref{f5}), (\ref{df}) well-posed in $V^{\beta +2\theta
_{2}}\times W^{l},$ thus generating a solution semigroup of operators%
\begin{align*}
\widehat{S}_{\theta _{2}}\left( t\right) & :V^{\beta +2\theta _{2}}\times
W^{l}\rightarrow V^{\beta +2\theta _{2}}\times W^{l}, \\
\left( v_{0},d_{0}\right) & \mapsto \left( v\left( t\right) ,d\left(
t\right) \right) .
\end{align*}%
More precisely, if the system (\ref{e:rel})-(\ref{example}) generates a
semigroup of solution operators $S_{\theta _{2}}$, as determined by the
conditions of Section \ref{s:well},
\begin{equation*}
S_{\theta _{2}}\left( t\right) :V^{\beta }\times W^{l}\rightarrow V^{\beta
}\times W^{l},\text{ }\left( u_{0},d_{0}\right) \mapsto \left( u\left(
t\right) ,d\left( t\right) \right)
\end{equation*}%
then this semigroup is linked through the corresponding semigroup $\widehat{S%
}_{\theta _{2}}\left( t\right) $ of any of the (oREL) problems above by the
relation%
\begin{equation*}
u\left( t\right) =Q^{-1}v\left( t\right) ,\text{ }\forall t\geq 0.
\end{equation*}%
Concerning the longtime behavior of $\widehat{S}_{\theta _{2}}$, then $%
\widehat{S}_{\theta _{2}}$ possesses a global attractor $\widehat{\mathcal{A}%
}$, which can be seen to satisfy%
\begin{equation*}
\widehat{\mathcal{A}}=\left\{ \left( v,d\right) \in V^{\beta +2\theta
_{2}}\times W^{l}:v=Qu,\text{ }\left( u,d\right) \in \mathcal{A}\right\} ,
\end{equation*}%
where $\mathcal{A}$ is a global attractor associated with any dynamical
system for the solution operator $S_{\theta _{2}}$.

Finally, our last comment is about assumption (\ref{select-p}). For this,
let us now consider $\Omega $ as a compact Riemannian manifold with boundary
$\Gamma $ and take again $E=T\Omega $ the tangent bundle. We observe that
the assumption on $\mathcal{V}$ in (\ref{select-p}) is satisfied in a more
general setting than suggested by the example given in Section \ref{ss:ee}.
Indeed, in the context of the specified regularized models of (1)-(5), it
suffices to consider $\mathcal{V}$ as a closed subspace of%
\begin{equation*}
\mathcal{V}_{\mathrm{ns}}=\{v\in C^{\infty }(T\Omega ):\mathrm{div}\left(
v\right) =0,\text{ }v\cdot n=0\text{ on }\Gamma \}.
\end{equation*}%
In this case, (\ref{select-p}) is clearly satisfied by the velocity $v\in
\mathcal{V}\subseteq \mathcal{V}_{\mathrm{ns}}$ of any of the problems
(1)-(5).

\section{Appendix}

\label{ss:app}

In this section, we include some supporting material on Gr\"{o}nwall-type
inequalities, Sobolev inequalities and abstract results. The first lemma is
a slight generalization of the usual Gr\"{o}nwall-type inequality.

\begin{lemma}
\label{Gineq}Let $\mathcal{E}:\mathbb{R}_{+}\rightarrow \mathbb{R}_{+}$ be
an absolutely continuous function satisfying
\begin{equation*}
\frac{d}{dt}\mathcal{E}(t)+2\eta \mathcal{E}(t)\leq h(t)\mathcal{E}%
(t)+l\left( t\right) +k,
\end{equation*}%
where $\eta >0$, $k\geq 0$ and $\int_{s}^{t}h\left( \tau \right) d\tau \leq
\eta (t-s)+m$, for all $t\geq s\geq 0$ and some $m\in \mathbb{R}$, and $%
\int_{t}^{t+1}l\left( \tau \right) d\tau \leq \gamma <\infty $. Then, for
all $t\geq 0$,
\begin{equation*}
\mathcal{E}(t)\leq \mathcal{E}(0)e^{m}e^{-\eta t}+\frac{2\gamma e^{m+\eta }}{%
e^{\eta }-1}+\frac{ke^{m}}{\eta }.
\end{equation*}
\end{lemma}

With $s,p\in \mathbb{R}_{+}$, let \thinspace $W^{s,p}$ be the standard
Sobolev space on an $n$-dimensional compact Riemannian manifold with $n\geq
2 $. The following result states the classical Gagliardo-Nirenberg-Sobolev
inequality (cf. \cite{BCL, He} and \cite{CM1, CM2}).

\begin{lemma}
\label{GNS}Let $0\leq k<m$ with $k,m\in \mathbb{N}$ and numbers $p,q,q\in %
\left[ 1,\infty \right] $ satisfy%
\begin{equation*}
k-\frac{n}{p}=\tau \left( m-\frac{n}{q}\right) -\left( 1-\tau \right) \frac{n%
}{r}.
\end{equation*}%
Then there exists a positive constant $C$ independent of $u$ such that%
\begin{equation*}
\left\Vert u\right\Vert _{W^{k,p}}\leq C\left\Vert u\right\Vert
_{W^{m,q}}^{\tau }\left\Vert u\right\Vert _{L^{r}}^{1-\tau },
\end{equation*}%
with $\tau \in \left[ \frac{k}{m},1\right] $ provided that $m-k-\frac{n}{r}%
\notin \mathbb{N}_{0},$ and $\tau =\frac{k}{m}$ provided that $m-k-\frac{n}{r%
}\in \mathbb{N}_{0}.$
\end{lemma}

We state here a standard result on pointwise multiplication of functions in
the Sobolev spaces $H^{k}=W^{k,2}$ (see~\cite{Ma04a}; cf. also \cite{HLT}).

\begin{lemma}
\label{l:hole} Let $s$, $s_{1}$, and $s_{2}$ be real numbers satisfying
\begin{equation*}
s_{1}+s_{2}\geq 0,\qquad min(s_{1},s_{2})\geq s,\qquad \text{and}\qquad
s_{1}+s_{2}-s>\frac{n}{2},
\end{equation*}%
where the strictness of the last two inequalities can be interchanged if $%
s\in \mathbb{N}_{0}$. Then, the pointwise multiplication of functions
extends uniquely to a continuous bilinear map
\begin{equation*}
H^{s_{1}}\otimes H^{s_{2}}\rightarrow H^{s}.
\end{equation*}
\end{lemma}

Our construction of an exponential attractor is based on the following
abstract result \cite[Proposition 4.1]{EZ}.

\begin{proposition}
\label{abstract}Let $\mathcal{H}$,$\mathcal{V}$,$\mathcal{V}_{1}$ be Banach
spaces such that the embedding $\mathcal{V}_{1}\hookrightarrow \mathcal{V}$
is compact. Let $\mathbb{B}$ be a closed bounded subset of $\mathcal{H}$ and
let $\mathbb{S}:\mathbb{B}\rightarrow \mathbb{B}$ be a map. Assume also that
there exists a uniformly Lipschitz continuous map $\mathbb{T}:\mathbb{B}%
\rightarrow \mathcal{V}_{1}$, i.e.,%
\begin{equation}
\left\Vert \mathbb{T}b_{1}-\mathbb{T}b_{2}\right\Vert _{\mathcal{V}_{1}}\leq
L\left\Vert b_{1}-b_{2}\right\Vert _{\mathcal{H}},\quad \forall
b_{1},b_{2}\in \mathbb{B},  \label{gl1}
\end{equation}%
for some $L\geq 0$, such that%
\begin{equation}
\left\Vert \mathbb{S}b_{1}-\mathbb{S}b_{2}\right\Vert _{\mathcal{H}}\leq
\gamma \left\Vert b_{1}-b_{2}\right\Vert _{\mathcal{H}}+K\left\Vert \mathbb{T%
}b_{1}-\mathbb{T}b_{2}\right\Vert _{\mathcal{V}},\quad \forall
b_{1},b_{2}\in \mathbb{B},  \label{gl2}
\end{equation}%
for some constant $0\le \gamma <\frac{1}{2}$ and $K\geq 0$. Then, there
exists a (discrete) exponential attractor $\mathcal{M}_{d}\subset \mathbb{B}$
of the semigroup $\{\mathbb{S}(n):=\mathbb{S}^{n},n\in \mathbb{Z}_{+}\}$
with discrete time in the phase space $\mathcal{H}$, which satisfies the
following properties:

\begin{itemize}
\item semi-invariance: $\mathbb{S}\left( \mathcal{M}_{d}\right) \subset
\mathcal{M}_{d}$;

\item compactness: $\mathcal{M}_{d}$ is compact in $\mathcal{H}$;

\item exponential attraction: $dist_{\mathcal{H}}(\mathbb{S}^{n}\mathbb{B},%
\mathcal{M}_{d})\leq C_{0}e^{-\chi n},$ for all $n\in \mathbb{N}$ and for
some $\chi >0$ and $C_{0}\geq 0$, where $dist_{\mathcal{H}}$ denotes the
standard Hausdorff semidistance between sets in $\mathcal{H}$;

\item finite-dimensionality: $\mathcal{M}_{d}$ has finite fractal dimension
in $\mathcal{H}$.
\end{itemize}
\end{proposition}

Moreover, the constants $C_{0}$ and $\chi ,$ and the fractal dimension of $%
\mathcal{M}_{d}$ can be explicitly expressed in terms of $L$, $K$, $\gamma $%
, $\left\Vert \mathbb{B}\right\Vert _{\mathcal{H}}$ and Kolmogorov's $\kappa
$-entropy of the compact embedding $\mathcal{V}_{1}\hookrightarrow \mathcal{V%
},$ for some $\kappa =\kappa \left( L,K,\gamma \right) $. We recall that the
Kolmogorov $\kappa $-entropy of the compact embedding $\mathcal{V}%
_{1}\hookrightarrow \mathcal{V}$ is the logarithm of the minimum number of
balls of radius $\kappa $ in $\mathcal{V}$ necessary to cover the unit ball
of $\mathcal{V}_{1}$.

\end{document}